\documentclass[hidelinks]{article}
\usepackage[utf8]{inputenc}
\usepackage[utf8]{inputenc}
\usepackage{mathrsfs}
\usepackage{amssymb}
\usepackage[margin=1in]{geometry}
\usepackage{changepage}
\usepackage{stmaryrd}
\usepackage{xcolor}
\usepackage{bbold}
\usepackage{amsmath,amsthm}
\usepackage{graphicx}
\usepackage{hyperref}
\usepackage{tikz}

\newtheorem{definition}{Definition}
\newtheorem{propo}{Proposition}
\newtheorem{remark}{Remark}
\newtheorem{lemma}{Lemma}

\newtheorem{theorem}{Theorem}
\newtheorem*{theorem*}{Theorem}

\numberwithin{equation}{section}
\newcommand{\vertiii}[1]{{\left\vert\kern-0.25ex\left\vert\kern-0.25ex\left\vert #1 
		\right\vert\kern-0.25ex\right\vert\kern-0.25ex\right\vert}}

\def\N{\mathbb{N}}
\def\R{\mathbb{R}}

\def\d{\,\mathrm{d}}

\def\B{\mathbb{B}}
\def\E{\mathbb{E}}
\def\P{\mathbb{P}}
\def\ssi{\Longleftrightarrow}

\title{\textbf{Weak well-posedness and weak discretization error for stable-driven SDEs with Lebesgue drift}}
\author{Mathis Fitoussi\footnotemark[1], \, Benjamin Jourdain\footnotemark[2], \, St\'ephane Menozzi\footnotemark[1]  }
\date{\today}

\begin{document}
	\maketitle
	
		\begin{abstract}
		We are interested in the discretization of stable driven SDEs with additive noise for $\alpha\in (1,2) $ and $L^q-L^p $ drift under the Serrin type condition $\frac{\alpha}{q}+\frac dp<\alpha-1 $. We show weak existence and uniqueness as well as heat kernel estimates for the SDE and obtain a convergence rate of order $\frac{1}{\alpha}\big(\alpha-1-(\frac{\alpha}{q}+\frac dp)\big) $ for the difference of the densities for the Euler scheme approximation involving suitably cutoffed and time randomized drifts.
		\end{abstract}

	\footnotetext[1]{Universit\'e Paris-Saclay, Laboratoire de Math\'ematiques et Mod\'elisation d’\'Evry (LaMME), 23 boulevard de France, 91 037 \'Evry, France\\}
    \footnotetext[2]{Cermics, \'Ecole des Ponts, INRIA, Marne-la-Vall\'ee, France.}
	\section{Introduction}\label{sec-intro}
	For a fixed time horizon $T>0$, we are interested in the weak well-posedness and the Euler-Maruyama dicretization of the SDE
	\begin{equation}\label{sde}
	\d X_t = b(t,X_t)\d t + \d Z_t, \qquad X_0 = x, \qquad \forall t \in [0,T],
	\end{equation}
	where $b$ belongs to the Lebesgue space $L^q ([0,T],L^p (\R^d)):= \left\{ f :[0,T] \times \R^d : \left\Vert t\mapsto \Vert f(t,\cdot) \Vert_{L^p} \right\Vert_{L^q([0,T])}  < \infty \right\}$ $=:L^q-L^p$ and $Z_t$ is a symmetric non-degenerate $d$-dimensional $\alpha$-stable process, whose spectral measure is equivalent to the Lebesgue measure on the unit sphere $\mathbb{S}^{d-1}$ (see Subsection \ref{subsec-density-noise} for detailed assumptions on the noise).\\ 

	 We will work under the integrability condition
	\begin{equation}\label{serrin}
		\frac{d}{p} + \frac{\alpha}{q}<\alpha -1, \qquad \alpha \in (1,2).
	\end{equation}
	
	This condition can be seen as the $\alpha$-stable extension of the Krylov-R\"ockner condition for Brownian-driven SDEs (see \cite{KR05}), although not guaranteeing \textit{strong} well-posedness in the strictly stable setting ($\alpha<2$). To this end, some additional smoothness conditions on the drift naturally appear, expressed in terms of Bessel potential spaces (see \cite{XZ20}).\\
	
	In this paper, we first establish well-posedness of \eqref{sde} through the study of a suitably associated Euler scheme, for which we prove heat kernel estimates. These then allow to follow the usual route to derive well-posedness: tightness, identification of a martingale problem solution and stability.	As a consequence of this approach, we derive Duhamel-type expansions for the densities of the Euler scheme and the diffusion, which paves the way for an error analysis.\\
	
	There has recently been a growing interest for SDEs of the type \eqref{sde} which involve a singular drift, both from the theoretical and numerical points of view. Drifts of the above form indeed appear in some physics-related models, having in mind, for example, the Biot-Savart kernel or Keller-Segel-type equations.\\
	
	This paper can be viewed as a stable-driven extension of \cite{JM23}, in which the corresponding Brownian case was addressed for the weak error. Stable processes naturally appear when modelling anomalous diffusion phenomena (see \cite{Esc06} for the fractional Keller-Segel model and \cite{MS12} for general fractional models). It is therefore important to be able to quantify how discretization schemes can approximate \eqref{sde}.
	
	\subsection{Definition of the Euler scheme}\label{subsec-def-euler-scheme}

	Since we consider a potentially unbounded drift coefficient, it is natural to introduce a cutoff for the discretization scheme. For a time step size $h$, the two cutoffs we consider are the following:
	\begin{itemize}
		\item If $p=q=\infty,$ we simply take $\forall (t,y)\in [0,T]\times\R^d,b_h(t,y)=\bar{b}_h(t,y)=b(t,y)$.
		\item Otherwise, we set
		\begin{align}
			b_h (t,y) &:= \frac{\min \left\{ |b(t,y)|,Bh^{-\frac{d}{\alpha p}-\frac{1}{q}}\right\}}{|b(t,y)|} b(t,y)\mathbb{1}_{|b(t,y)|>0} ,\qquad &(t,y)\in [0,T]\times\R^d,  \label{cutoff-1} \\
			\bar{b}_h (t,y) &= \frac{\min \left\{ |b(t,y)|,Bh^{\frac{1}{\alpha}-1}\right\}}{|b(t,y)|} b(t,y)\mathbb{1}_{t\geq h,|b(t,y)|>0}, &(t,y)\in [0,T] \times \R^d, \label{cutoff-2}
		\end{align}
		for some constant $B>0$ which can be chosen freely as long as it does not depend on $h$ nor $T$.
			\end{itemize}
		The first option has a cutoff level related to the integrability condition \eqref{serrin}, while the second one is related to the auto-similarity index of the driving noise. The latter also artificially sets the drift to 0 on the first step (we will see later that this allows in particular to compute estimates on the gradient of the density of the Euler scheme). The idea behind this cutoff level is to make sure the contribution of the drift does not dominate over that of the stable noise on each time step of the scheme.\\
			
		We then define a discretization scheme with $n$ time steps over $[0,T]$, with constant step size $h:=T/n$. For the rest of this paper, we denote, $\forall k \in \{1,...,n\}, t_k := kh$ and $\forall s>0, \tau_s^h := h \lfloor \frac{s}{h} \rfloor \in (s-h,s]$, which is the last grid point before time $s$. Namely, if $s\in [t_k,t_{k+1}), \tau_s^h = t_k$. \\
		
		In order to avoid assumptions on the drift $b$ beyond integrability and measurability, we are led to randomize the evaluations of $b_h$ (resp. $\bar b_h$) in the time variable. For each $k\in \{0,...,n-1\}$, we will draw a random variable $U_k$ according to the uniform law on $[kh,(k+1)h]$, independently of each other and the noise $(Z_t)_{t\geq 0}$.	We can then define a step of the Euler scheme as
		\begin{equation}\label{euler-scheme}
			X_{t_{k+1}}^h = X_{t_k}^h + (Z_{t_{k+1}}-Z_{t_{k}})+ hb_h(U_k,X_{t_k}^h),
		\end{equation}
		and its time interpolation as the solution to
		\begin{equation}\label{scheme-interpo}
		\d X^h_t=\d Z_t+b_h(U_{\lfloor \frac t h\rfloor},X^h_{\tau^h_t}) \d t.
	\end{equation}
	Similarly, for the alternative cutoff, we define
	\begin{equation}\label{euler-scheme-bar}
		\bar X_{t_{k+1}}^h = \bar X_{t_k}^h + (Z_{t_{k+1}}-Z_{t_{k}})+ h\bar b_h(U_k,X_{t_k}^h),
	\end{equation}
	and its time interpolation as the solution to
	\begin{equation}\label{scheme-interpo-bar}
		\d \bar X^h_t=\d Z_t+\bar b_h(U_{\lfloor \frac t h\rfloor},\bar X^h_{\tau^h_t}) \d t.
	\end{equation}
	As $b_h$ and $\bar b_h$ are bounded, the schemes \eqref{scheme-interpo} and \eqref{scheme-interpo-bar} are well defined and admit densities in positive times. We will denote by $\Gamma^h (0,x,t,\cdot)$ and $\bar \Gamma^h (0,x,t,\cdot)$ their respective densities at time $t\in (0,T]$ when starting from $x$ at time 0.
	\subsection{Well-posedness - state of the art}
	
	\color{red}

	\color{black}Let us recall that weak well-posedness is often investigated through the parabolic PDE which is naturally associated with the SDE \eqref{sde}
	\begin{equation}\label{pde}
		\left( \partial_s + b(s,x)\cdot \nabla_x + \mathcal{L}^\alpha \right)u(s,x)=f(s,x) \; \mathrm{on} \;[0,\textcolor{red}{T})\times \R^d, \qquad u(t,\cdot) = g \; \mathrm{on} \; \R^d,
	\end{equation}
	where $\mathcal{L}^\alpha$ is the generator of the noise and $f$ and $g$ are suitable functions. Bearing in mind that, in the $\beta$-H\"older setting, the associated parabolic gain is $\beta+\alpha$, the regularity condition $\beta+\alpha >1$ naturally appears to define the gradient of the solution. Let us point out that this condition already appeared in the seminal work of \cite{TTW74}. For weak and strong well-posedness in the H\"older setting, we can e.g. refer to \cite{MP14} and \cite{CZZ21}, which also includes the super-critical case.	Since we do not have any regularity available on the drift $b$, we are naturally led to consider sub-critical regimes for the stability index (i.e. $\alpha>1$).\\
	
	Establishing estimates on the gradient of the solution to the PDE naturally leads to weak uniqueness in the multidimensional setting for \eqref{sde} through the martingale problem. In this paper, under \eqref{serrin}, we obtain such estimates exploiting heat kernel estimates for the density $\bar \Gamma^h$ of $\bar X^h $ and taking the limit as $h$ goes to 0. Keep in mind that some additional smoothness is required to derive strong well-posedness in the multidimensional case. \\
	
	In the strictly stable and time-homogeneous setting with mere integrability assumptions on the drift, weak existence and uniqueness of a solution to \eqref{sde} was first investigated in \cite{Por94} in $\R$ and extended to the multidimensional case in \cite{PP95} under the condition $\frac{d}{p}<\alpha-1$ \textcolor{black}{by constructing the density using its parametrix expansion.} When considering the embedding $L^p (\R^d) \hookrightarrow \B_{\infty,\infty}^{-\frac{d}{p}}(\R^d)$ (the latter being the Besov space with regularity $-\frac{d}{p}$), the previous condition is then consistent with the condition $\alpha+\beta>1$ appearing in the H\"older case.

	 Let us also mention the work \cite{CdRM22}, in which weak well-posedness is proved for distributional drifts in the Besov-Lebesgue space $L^q-\B_{p,r}^\beta$ under the condition $\beta > \frac{1-\alpha+\frac{\alpha}{q}+\frac{d}{p}}{2}$. In view of this threshold, our well-posedness result can be seen as an extension of this work for $\beta=0$.\\
	
	Our approach to well-posedness naturally provides heat kernel estimates for both the discretization scheme and the limit SDE that quantify the behavior of their time marginal laws. Namely, as detailed in the seminal work by Kolokoltsov \cite{Kol00}, for a smooth bounded drift, the time marginals of the solution \eqref{sde} and the noise behave alike. This work was then extended in various directions, mostly for H\"older continuous drifts (see \cite{KK18}, \cite{Kul19}, \cite{CHZ20} and \cite{MZ22}), and more recently for distributional drifts (see \cite{PvZ22} in the Brownian setting and \cite{Fit23} in the strictly stable case). In those works, the authors again establish that the time marginal laws of the process have a density which is ``equivalent" (i.e. bounded from above and below) to the density of the noise, and that the spatial gradients exhibit the same time singularities and decay rates (see Theorem \ref{thm-wp-diffusion} below in the current Lebesgue setting).
	
	\subsection{Euler scheme - state of the art}
 	
 	For the discretization of singular drift diffusions, a rather vast literature exists, although it mostly focuses on the Brownian setting for an additive noise. A first approach consists in using the sewing lemma (see \cite{Le20}) in order to obtain results on the strong error rate, which is defined as the convergence rate of
 	\begin{equation}
 		\left\Vert \sup_{t\in (0,T)} |X_t-X_t^h| \right\Vert_{L^r}
 	\end{equation}
 	for some $r>1$. This was done in the work of L\^e and Ling (\cite{LL22}), who obtain a convergence in $h^{\frac{1}{2}} |\ln (h)|$ under the Krylov-R\"ockner condition $\frac{d}{p}+\frac{2}{q}<1$ (see also \cite{DGI22}) even with multiplicative noise (when the corresponding coefficient is Lipschitz in the spatial variable) for the semi-discrete scheme where the time-variable of the coefficients is not discretized. This is a remarkable result since, up to the logarithmic factor, this corresponds to the convergence rate for the strong error associated with a Brownian SDE with Lipschitz coefficients with non-trivial diffusion term. 
        It remains open to understand whether the strong convergence rate improves in terms of the gap to criticity $1-\big(\frac{d}{p}+\frac{2}{q}\big)$ in the additive noise case. \\
 	
 	The main contribution of the sewing lemma consists in bounding $L^r$ norms of the form
	\begin{equation}\label{sewing-contribution}
		\E \left[ \left| \int_{0}^{t}b(s,X_s^h)-b(s,X_{\tau_s^h}^h) \d s \right|^r \right],
	\end{equation}
 	that is, the strong error associated with local differences of the path along an irregular function with suitable integrability properties.\\

 	On the other hand, deriving weak error rates usually involves  studying the PDE \eqref{pde} or the associated Duhamel representation. Indeed, the weak error is related to the difference between the density of the SDE \eqref{sde} and that of the corresponding Euler scheme \eqref{scheme-interpo}. Using the Duhamel representations satisfied by the respective transition densities $\Gamma$ and $\Gamma^h$ of the diffusion and its Euler scheme, we will estimate
 	$$|\Gamma (0,x,t,y) - \Gamma^h(0,x,t,y)|.$$ This approach allows to integrate against any type of irregular test functions enjoying suitable integrability properties.\\
 	
 	When the coefficients of \eqref{sde} are smooth, the seminal paper of Talay and Tubaro (\cite{TT90}) gives a convergence of order 1 in $h$. With $\beta$-H\"older coefficients, the work of Mikulevicius and Platen (\cite{MP91}) proves a convergence in $h^\frac{\beta}{2}$ in the Brownian case. In these works, for $u$ solving \eqref{pde} with smooth terminal condition $g$ and no source term $f$, applying Itô's formula, the error writes
 	\begin{align*}
 		\mathcal{E}(g,t,x,h)&=\E_{0,x}[ g(X_t^h)-g(X_t)]=\E_{0,x}[ u(t,X_t^h)-u(0,x)]\\
 		&=\E_{0,x}\left[\int_0^t \left( b(r,X_r^h)-b(\tau_r^h,X_{\tau_r^h}^h) \right) \cdot \nabla u (r,X_r^h)\right] \d r,
 	\end{align*}
     where the index $0,x$ of the expectation sign means that the scheme is started from $X^h_0=x$ at time $0$.
 	The authors then use classic Schauder type estimates, see e.g. \cite{AF64}, to control $\Vert \nabla u\Vert_{L^\infty}$. From the $\beta$-H\"older continuity of the drift, the following bound is then derived
 	\begin{equation}\label{rate-MP91}
 		\mathcal{E}(g,t,x,h)\leq C \Vert \nabla u \Vert_{L^\infty} \int_0^t \E _{0,x}\left[|X_r^h -X_{\tau_r^h}^h|^\beta \right]\d r\leq C \Vert \nabla u \Vert_{L^\infty} h^{\frac{\beta}{2}}.
 	\end{equation}
 	 The above final rate then comes from the magnitude of the increment of the Euler scheme on one time step in the $L^\beta (\mathbb{P})$ norm. However, one can see that this essentially consists in using strong error analysis techniques to derive a weak error rate, which does not seem adequate.\\
 	 
 	 In the current work, we want to investigate errors of the form $\mathcal{E}(\delta_y,t,x,h)$ (where $\delta_y$ is the dirac mass at point $y$). This formally writes
 	 \begin{align}\label{errdelt}
			\mathcal{E}(\delta_y,t,x,h)=\E _{0,x}\left[\int_0^t \left( b(r,X_r^h)-b(\tau_r^h,X_{\tau_r^h}^h)) \right) \cdot \nabla_z \Gamma (r,t,z,y)|_{z=X_r^h} \d r\right],
 	 \end{align}
 	 where $\Gamma$ is the density of \eqref{sde}. When comparing this equation to \eqref{sewing-contribution}, we see that, in the weak setting, an additional gradient term appears. Whenever this term is not regular enough, which is the case in the current Lebesgue setting, it lowers the time integrability properties of the irregular function that we want to investigate. However, in the specific case of a H\"older continuous drift and terminal condition $g$, this additional term can be handled using sewing techniques. Doing so in \cite{Hol22}, the author improves the convergence rate in \eqref{rate-MP91} to $h^{\frac{\beta+1}{2}}$. The study of the weak error for H\"older coefficients and a final Dirac mass will be the topic of an upcoming work.\\
 	
 	In the irregular setting, for the weak error associated with the densities, the randomization of the time variable permits to replace $b(\tau_r^h,X_{\tau_r^h}^h)$ by $b(r,X_{\tau_r^h}^h) \d u$ in \eqref{errdelt} (up to some error term on $[\tau^h_t,t]$) and another new idea was introduced in \cite{BJ20} in order to tackle mere bounded drifts, which consists in writing
 	\begin{align}
 		\E_{0,x} [b(r,X_r^h)\cdot\nabla \Gamma (r,t,X_r^h,y)&-b(r,X_{\tau_r^h}^h)\cdot\nabla \Gamma (r,t,X_{\tau_r^h}^h,y) ]\notag\\&= \int [ \Gamma^h (0,x,r,z) - \Gamma^h (0,x,\tau_r^h,z) ]b(r,z) \cdot\nabla \Gamma (r,t,z,y)\d z\label{bj20}
 	\end{align}
 	and exploiting the regularity in the forward time variable of $\Gamma^h$ instead of that of $b\cdot \nabla \Gamma$. In \cite{JM23}, authors use this technique with a drift in $L^q-L^p$ to derive a rate of order $\frac{\gamma}{2}$, where $\gamma:= 1-\frac{d}{p}-\frac{2}{q}$ is the so-called ``gap to singularity" or ``gap to criticity". Note that, with respect to the rate obtained in \cite{LL22}, due to the additional \textit{gradient} term in $\nabla \Gamma$ in \eqref{bj20} (as opposed to \eqref{sewing-contribution}), an order $\frac{1}{2}$ is lost on the convergence rate. However, the techniques developed therein allow to take advantage of the gap to singularity. \\
  
    As mentioned, the rate for the strong error under the Krylov-R\"ockner condition is (at least) $\frac{1}{2}$, up to a logarithmic factor. Since we expect the weak error rate to be at least as good, it remains to understand how to improve it beyond $\frac{1}{2}$.\\
 	
 	
 	In Theorem \ref{thm-main}, we obtain a weak error rate in $h^{\frac{\gamma}{\alpha}}$, where our ``gap to singularity" is now defined as 
	$\gamma := \alpha-1-\Big(\frac{d}{p}+\frac{\alpha}{q}\Big)>0$. Importantly, there is continuity w.r.t. the stability index for the associated error rates.
 	
	\subsection{Driving noise and related density properties}\label{subsec-density-noise}

	Let us denote by $\mathcal{L}^\alpha$ the generator of the driving noise $Z$. In the case $\alpha = 2$, $\mathcal{L}^\alpha$ is the usual Laplacian $\frac{1}{2}\Delta$. When $\alpha \in (1,2)$, in whole generality, the generator of a symmetric stable process writes, $\forall \phi \in C_0^\infty (\R^d,\R)$ (smooth compactly supported functions),
	\begin{align*}
		\mathcal{L}^\alpha \phi (x)&= \mathrm{p.v.} \int_{\R^d} \left[ \phi(x+z) - \phi(x)\right]\nu(\mathrm{d} z)\\
		&=\int_{\R_+}\int_{\mathbb{S}^{d-1}}\left[ \phi(x+\rho \xi) - \phi(x)\right]\mu(\mathrm{d} \xi)\frac{\d \rho}{\rho^{1+\alpha}}
	\end{align*}
	(see \cite{Sat99} for the polar decomposition of the spectral measure) where $\mu$ is a finite measure on the unit sphere $\mathbb{S}^{d-1}$ such that $\mu(A)=\mu(\{\xi\in\mathbb{S}^{d-1}:-\xi\in A\})$ for each Borel subset $A$ of $\mathbb{S}^{d-1}$.
        
	This general setting will not allow us to derive heat kernel estimates, because it does not lead to global estimates of the noise density. In \cite{Wa07}, Watanabe investigates the behavior of the density of an $\alpha$-stable process in terms of properties fulfilled by the support of its spectral measure. From this work, we know that whenever the measure $\mu$ is not equivalent to the Lebesgue measure $m$ on the unit sphere, accurate estimates on the density of the stable process are delicate to obtain. However, Watanabe (see \cite{Wa07}, Theorem 1.5) and Kolokoltsov (\cite{Kol00}, Propositions 2.1--2.5) showed that if \begin{equation}
         c^{-1} m(\mathrm{d} \xi) \leq \mu (\mathrm{d} \xi) \leq c m (\mathrm{d} \xi) \mbox{ for some }c\in [1,+\infty),\label{minmajleb}
        \end{equation} the following estimates hold for the density $z\mapsto p_\alpha(v,z)$ of $Z_v$ with respect to the Lebesgue measure on $\R^d$ when $v>0$ : there exists a constant $C$ depending only on $\alpha,d$, s.t. $\forall v\in \R_+^*, z\in \R^d$,
    \begin{equation}\label{ARONSON_STABLE}
    C^{-1}v^{-\frac{d}{\alpha}}\left( 1+ \frac{|z|}{v^{\frac{1}{\alpha}}} \right)^{-(d+\alpha)}\leq p_\alpha(v,z)\leq Cv^{-\frac{d}{\alpha}}\left( 1+ \frac{|z|}{v^{\frac{1}{\alpha}}} \right)^{-(d+\alpha)}.
    \end{equation}
	As our approach heavily relies on these global bounds, we assume that $\mu$ satisfies \eqref{minmajleb}.   
	
	Note that in Section \ref{sec-lemmas} and Appendix \ref{APP_HK} which are dedicated to technical lemmas, we will be using the proxy notation 
	\begin{equation}
		\bar{p}_\alpha (v,z) := \frac{C_\alpha}{v^{\frac{d}{\alpha}}} \left(1+\frac{|z|}{v^{\frac{1}{\alpha}}} \right)^{-(d+\alpha)},\qquad v>0,z\in \R^d,\label{DEF_P_BAR}
	\end{equation}
	where $C_\alpha$ is chosen so that $\forall v>0, \int \bar{p}_\alpha (v,y) \d y = 1$, because we therein explicitly rely on the global bounds provided by $\bar{p}_\alpha$. In the rest of the paper, we will prefer the notation $p_\alpha$, directly referring to the density of the noise.\\
	
	Further properties related to the density of the driving noise are stated in Lemmas \ref{lemma-stable-sensitivities} and \ref{lemFK} below.
	
		\subsection{Main results}\label{sec-main-results}
We are now in position to state the main results of the current work.	The first result concerns the well-posedness of \eqref{sde}. 

	\begin{theorem}[Weak existence and density estimates for the diffusion]\label{thm-wp-diffusion}
		Assume \eqref{serrin}. The stochastic differential equation \eqref{sde} admits a weak solution such that for each $t\in (0,T]$, $X_t$ admits a density $y\mapsto \Gamma (0,t,x,y)$ w.r.t. the Lebesgue measure such that $\exists C:=C(b,T)<\infty : \forall t \in (0,T], \forall (x,y)\in (\R^d)^2$,
		\begin{equation}\label{ineq-density-diffusion}
                        \Gamma (0,x,t,y) \leq C p_\alpha (t,y-x),
		\end{equation}
and this density is the unique solution to the following Duhamel representation among functions of $(t,y)\in[0,T]\times\R^d$ satisfying \eqref{ineq-density-diffusion}: 
		\begin{equation}\label{duhamel-diffusion}
		\forall t\in (0,T],\;\forall y\in \R^d,\;	\Gamma(0,x,t,y) =  p_\alpha (t,y-x) - \int_0^t\int_{\R^d}\Gamma(0,x,r,z)b(r,z).\nabla_yp_\alpha(t-r,y-z)\d z \d r. 
		\end{equation}
		\textcolor{black}{Furthermore, there exists a unique solution to the martingale problem related to $b\cdot \nabla + \mathcal{L}^\alpha$ starting from $x$ at time $0$ in the sense of Definition \ref{DEF_MART_PB} (see page \pageref{DEF_MART_PB} below).}\\
		
		\noindent Finally, let us define the ``gap to singularity" as
		\begin{equation}
			\gamma:= \alpha -1 -\Big(\frac{d}{p}+\frac{\alpha}{q}\Big) >0.\label{gap}
		\end{equation}
		Then, $\Gamma$ has the following regularity in the forward spatial variable: $\forall t \in (0,T], \forall (x,y,y')\in (\R^d)^3,$
		\begin{equation}\label{holder-space-gamma}
			|\Gamma (0,x,t,y)-\Gamma (0,x,t,y')| \leq C \frac{|y-y'|^{\gamma}\wedge t^{\frac{\gamma}{\alpha}}}{ t^{\frac{\gamma}{\alpha}}} \left( p_\alpha(t,y-x)+p_\alpha(t,y'-x)\right).
		\end{equation}
	\end{theorem}
The proof of the heat kernel estimates for the diffusion heavily relies on the following heat kernel estimates for the density the Euler schemes  \eqref{scheme-interpo} and \eqref{scheme-interpo-bar}.
\begin{propo}[Density estimates for the Euler scheme]\label{prop-main-estimates}
		Assume \eqref{serrin}. Set $h=\frac{T}{n},\ n\in \mathbb N^*$. Let $X^h$ be the scheme defined  in \eqref{scheme-interpo} (resp. $\bar X^h$ the scheme defined in \eqref{scheme-interpo-bar}) starting from $x_0\in \R^d$ at time $0$.\\
		Then, for all $0\le t_k:=kh<t\le T ,\ k\in \llbracket 0,n-1 \rrbracket, (x,y)\in (\R^d)^2$, the random variable $X_t^h$ admits, conditionally to $X_{t_k}^h=x$, a 
                density $\Gamma^h(t_k,x,t,\cdot) $, which enjoys the following Duhamel representation: for all $y\in \R^d $,
		\begin{align}
			\Gamma^h(t_k,x,t,y)
			=p_\alpha(t-t_k,y-x)-\int_{t_k}^{ t}\E_{t_k,x}\left[b_h(U_{\lfloor \frac rh\rfloor },X^h_{\tau_r^h})\cdot\nabla_y p_\alpha(t-r,y-X^h_r)\right]\d r,\label{duhamel-scheme}
		\end{align}
		where the index $t_k,x$ of the expectation sign means that the scheme $(X^h_r)_{r\in[t_k,T]}$ is started from $X^h_{t_k}=x$ at time $t_k$.
		Similarly, the random variable $\bar X_t^h$ admits, conditionally to $\bar X_{t_k}^h=x$,  a transition density $\bar \Gamma^h(t_k,x,t,\cdot) $, which enjoys the following Duhamel representation: for all $y\in \R^d $,
		\begin{align}
			\bar\Gamma^h(t_k,x,t,y)	=p_\alpha(t-t_k,y-x)-\int_{t_k}^{ t}\E_{t_k,x}\left[\bar b_h(U_{\lfloor \frac rh\rfloor },\bar X^h_{\tau_r^h})\cdot\nabla_y p_\alpha(t-r,y-\bar X^h_r)\right]\d r.\label{duhamel-scheme-bar}
		\end{align}
		Furthermore, there exists a finite constant $C$ not depending on $h=\frac{T}{n}$ such that for all $k\in \llbracket 0,n-1\rrbracket,
		t\in (t_k,T], x,y,y'\in\R^d$,
		\begin{align}\label{ineq-density-scheme}
			&\Gamma^h(t_k,x,t,y)\le C p_\alpha (t-t_k,y-x)
		\end{align}
		and 
		\begin{align}\label{holder-space-gammah}
			|\Gamma^h(t_k,x,t,y') &-\Gamma^h(t_k,x,t,y)|\notag\\
			\le& 
			C \frac{|y-y'|^{\gamma}\wedge (t-t_k)^{\frac{\gamma}{\alpha}}}{(t-t_k)^{\frac{\gamma}{\alpha}}}\big(p_\alpha(t-t_k,y-x)+p_\alpha(t-t_k,y'-x)\big),
		\end{align}
for $\gamma  $ defined in \eqref{gap}.
		Also, for all $0\le k<\ell<n,\; t\in[t_\ell,t_{\ell+1}], x,y\in\R^d$,
		\begin{align}\label{holder-time-gammah}
			|\Gamma^h(t_k,x,t,y)-\Gamma^h(t_k,x,t_\ell,y)|
			\le C \frac{(t-t_\ell)^{\frac{\gamma}{\alpha}}}{(t_\ell-t_k)^{\frac {\gamma}{\alpha}}}
			p_\alpha(t-t_k,y-x),
		\end{align}
		and the same estimates hold with $\bar\Gamma^h$ replacing $\Gamma^h$. 
	\end{propo}	Our second main result states a weak convergence rate bound for the Euler schemes \eqref{scheme-interpo} and \eqref{scheme-interpo-bar} :

	\begin{theorem}[Convergence Rate for the stable-driven Euler-Maruyama scheme with $L^q-L^p $ drift]\label{thm-main}
		Assume that \eqref{serrin} holds. 
                There exists a constant $C<\infty$ s.t. for  all $h=T/n$ with $n\in\N^*$, and all $t\in(0,T]$, $x,y\in \R^d $
		\begin{align*}
			|\Gamma^h(0,x,t,y)-\Gamma(0,x,t,y)| &\leq C  h^{\frac{\gamma}{\alpha}} p_\alpha(t,y-x),\\
			\mathrm{resp.} \qquad 	|\bar \Gamma^h(0,x,t,y)-\Gamma(0,x,t,y)| &\leq C  h^{\frac{\gamma}{\alpha}} p_\alpha(t,y-x).
		\end{align*}
		
	\end{theorem}

	\subsection{Notations}\label{subsec-def}
	We will use the following notations : 
	\begin{itemize}
		\item $A\lesssim B $ if there exists a constant $C\ge 1$, which depends only on $\alpha,d,p,q,b,T$, such that $A\leq CB$.
		\item $A \asymp B$ if there exists a constant $C\ge 1$, which depends only on $\alpha,d,p,q,b,T$, such that $ C^{-1} B \leq A\leq CB$.
		\item For $\ell\in [1,+\infty]$, we always denote by $\ell'\in [1,+\infty]$ its conjugate exponent, i.e. $\frac{1}{\ell}+ \frac{1}{\ell'}=1$.
	\end{itemize}
	The article is organized as follows. The proof of Theorem \ref{thm-main}
	is developed in Section \ref{SEC_PROOF_THM_ERROR} (assuming that the controls of Theorem \ref{thm-wp-diffusion} hold). Section  \ref{subsec-proof-prop1} is dedicated to the proof of the estimates for the schemes. The proof of Theorem \ref{thm-wp-diffusion} is presented in Section \ref{SEC_MART_PB}.
The proof of some technical results are gathered in Appendix \ref{APP_HK}. 	

\section{Proof of the convergence rate for the error (Theorem \ref{thm-main})}
\label{SEC_PROOF_THM_ERROR}

	\subsection{Technical tools}\label{sec-lemmas}
	We will profusely use the following technical lemmas which hold for any stability index $\alpha\in (1,2) $ and are proved in Appendix \ref{APP_HK}:
	\begin{lemma}[Stable sensitivities - Estimates on the $\alpha $-stable kernel]\label{lemma-stable-sensitivities}
		For each multi-index  $\zeta$ with length $|\zeta|\leq2$, and for all $0<u\leq u'\leq T$, $(x,x')\in (\R^d)^2$,
		\begin{itemize} 
			\item Bounds for space and time derivatives: for all $\beta\in \{0,1\} $,
			\begin{align}
				|\partial_u^\beta \nabla_x^\zeta p_\alpha (u,x)| &\lesssim \frac{p_\alpha (u,x)}{u^{\beta+\frac{|\zeta|}{\alpha}}}.\label{derivatives-palpha}
			\end{align}
			\item Spatial H\"older regularity: for all $\theta \in (0,1]$,
			\begin{align}\label{holder-space-palpha}
				\left|\nabla_x^\zeta p_\alpha (u,x)-\nabla_x^\zeta p_\alpha (u,x')\right| \lesssim \left(\frac{|x-x'|^\theta}{u^{\frac{\theta}{\alpha}}} \wedge 1\right)\frac{1}{u^{\frac{|\zeta|}{\alpha}}}\left(p_\alpha (u,x)+p_\alpha (u,x')\right).
			\end{align}
			\item Time H\"older regularity: for all $\theta\in (0,1] $,
			\begin{align}\label{holder-time-palpha}
				\left|\nabla_x^\zeta p_\alpha (u,x)-\nabla_x^\zeta p_\alpha (u',x)\right| \lesssim \frac{|u-u'|^\theta}{u^{\theta+\frac{|\zeta|}{\alpha}}} \left(p_\alpha (u,x)+p_\alpha (u',x)\right).
			\end{align}
			\item Time scales for spatial moments: for all $\ell\in [1,+\infty] $ and  $\delta \in [0,\frac{d}{\ell}+\alpha)$, 			
			\begin{align}\label{spatial-moments}
				\Vert p_\alpha (u,\cdot)|\cdot|^\delta \Vert_{L^{\ell '}}\leq C u^{-\frac{d }{\alpha\ell}+\frac{\delta }{\alpha}}.
			\end{align}
   
			\item Convolution: for all $(x,y) \in (\R^d)^2$, $ 0 < s < u < t\le T$, $\ell\geq 1$,
			\begin{equation}\label{p-q-convo}
				\Vert {p}_\alpha (t-u,\cdot-y) {p}_\alpha (u-s,x-\cdot)\Vert_{L^{\ell'}} \lesssim \left[ \frac{1}{(t-u)^{\frac{d}{\alpha  \ell}}} +\frac{1}{(u-s)^{\frac{d}{\alpha  \ell}}}\right] {p}_\alpha (t-s,x-y).
			\end{equation}
			\item Integration of an $L^{\ell} $ function in a spatial stable convolution: for all $(x,y) \in (\R^d)^2$, $ 0 \leq s <u < t\le T$, $\ell \geq 1, \phi \in L^\ell(\R^d,\R)$,
			\begin{equation}\label{convo-space-sing}
				\int p_\alpha (t-u,z-x)|\phi(z)| p_\alpha(u-s,y-z) \d z \lesssim \left[\frac{1}{(t-u)^{\frac{d}{\alpha \ell}}} + \frac{1}{(u-s)^{\frac{d}{\alpha \ell}}}\right] p_\alpha(t-s,y-x) \Vert \phi \Vert_{L^\ell}.
			\end{equation}
		\end{itemize}
	\end{lemma}
	
    \begin{lemma}[Feynman-Kac partial differential equation]\label{lemFK}
          Let $t>0$ and $\phi:\R^d\to\R$ be a ${\cal C}^2$ function with bounded derivatives. Then the function $v(s,y)={\mathbb 1}_{s<t}p_\alpha (t-s,\cdot)\star \phi(y)+{\mathbb 1}_{s=t}\phi(y)$ is ${\cal C}^{1,2}$ on $[0,t]\times \R^d$ and satisfies the Feynman-Kac partial differential equation
$$\forall (s,y)\in[0,t)\times\R^d,\;\partial_s v(s,y)+\mathcal{L}^\alpha v(s,y)=0 .$$\end{lemma}
        
	\begin{lemma}[Integration of the drift in a spatial stable time-space convolution]\label{lemma-convo-bulk} Let $0 \leq u < v \leq  t \leq T$ and $\beta_1,\beta_2 \in \R_+$. Let $b\in L^q([0,T], L^p(\R^d))$ with $p,q$ such that \eqref{serrin} holds.
		\begin{itemize}
			\item Singular case. If $v<t$ and 
			\begin{align*}
				q' \left(\frac{d}{\alpha p}+\beta_1\right)>1\qquad \mathrm{and} \qquad	q' \left(\frac{d}{\alpha p}+\beta_2\right)<1,
			\end{align*}
			then, 
			\begin{align}\label{lemma-convo-bulk-1}
				\int_{u}^{v}\int & p_\alpha (r,z-x)|b(r,z)|p_\alpha (t-r,y-z) \frac{1}{(t-r)^{\beta_1}}\frac{1}{r^{\beta_2}} \d r \notag\\
				&\lesssim  		p_\alpha(t,y-x)\Big((v-u)^{\frac{\gamma+1}\alpha
				-(\beta_1+\beta_2)}+(v-u)^{-\beta_2}(t-v)^{
				\frac{\gamma+1}\alpha-\beta_1
				}\Big).
			\end{align}
			\item Integrable case. If
			\begin{align*}
								q' \left(\frac{d}{\alpha p}+\beta_1\right)<1 \qquad \mathrm{and} \qquad	q' \left(\frac{d}{\alpha p}+\beta_2\right)<1,
			\end{align*}
			then, 
			\begin{align}\label{lemma-convo-bulk-2}
				\int_{u}^{v}\int p_\alpha (r,z-x)|b(r,z)|p_\alpha (t-r,y-z) \frac{1}{(t-r)^{\beta_1}}\frac{1}{r^{\beta_2}} \d r 
				\lesssim p_\alpha (t,y-x) [(v-u)^{\frac{\gamma+1}\alpha-(\beta_1+\beta_2)}].
			\end{align}
		\end{itemize}
	\end{lemma}

    The previous lemma will be used to treat the main error terms in the analysis of the error. The most common use case is when $\beta_2=0$ and $\beta_1=\frac{1}{\alpha}$ (we are thus in case \eqref{lemma-convo-bulk-2}) and $u=h, v=\tau_t^h-h$. This configuration appears when we previously used \eqref{derivatives-palpha} to bound the gradient of $p_\alpha (t-r,y-z)$ and that no other singularities come into play. 
    The case $\beta_2>0$ with an additional singular in $r$ factor is needed for the proof of Theorem \ref{thm-main} 
    (which will \textcolor{black}{require setting} $ \beta_2=\frac \gamma\alpha$).\\
    We will also use \eqref{lemma-convo-bulk-1} whenever there is an additional singularity in $(t-r)$ which makes the previous integral non-convergent (see e.g. \eqref{def-delta-4-i}). This will actually happen in order to obtain exactly the gap $\gamma $ defined \eqref{gap} in the convergence rate or in the Hölder exponents for the density, see e.g. \textcolor{black}{Section \ref{SUB_H_REG_IN_F_TIME} for the proof of the Hölder regularity of the density of the scheme stated in Proposition \ref{prop-main-estimates}}.
 
	\begin{remark}
		From the definition of $\bar{p}_\alpha(u,x)=\frac{C_\alpha}{u^{\frac{d}{\alpha}}}\left( 1+ \frac{|x|}{u^{\frac{1}{\alpha}}} \right)^{-(d+\alpha)}$, one can gather the following:\\
		Let $x\in \R^d$ and $u>0$.
		\begin{itemize}
			\item If $|x|\ge u^{\frac{1}{\alpha}}$ (off-diagonal regime),
			\begin{equation}\label{off-diag}
				\bar{p}_\alpha (u,x) \asymp \frac{u}{|x|^{d+\alpha}}.
			\end{equation}
			\item If $|x|\leq u^{\frac{1}{\alpha}}$ (diagonal regime),
			\begin{equation}\label{diag}
				\bar{p}_\alpha (u,x) \asymp \frac{1}{u^{\frac{d}{\alpha}}}.
			\end{equation}
		\end{itemize}
	\end{remark}
	
	Those two regimes will be central in our proofs. The scales which we consider for these regimes derive from the self-similarity of the noise.  Let us as well point out that the diagonal bound in \eqref{diag} is also a global upper bound for both $\bar p_\alpha $ and $p_\alpha $.\\
	
	The next lemma is very important since it precisely emphasizes that the drift $b_h$ (resp. $\bar b_h $) is actually a \textit{negligible} term w.r.t. the scale of the underlying noise for a one-step transition of the corresponding scheme.
	\begin{lemma}[About the cutoff on a one-step transition]\label{lemma-ignore-bh}
		\;Here, $\mathfrak b_h\in \{b_h,\bar b_h \} $ stands for one of the two drifts considered for the schemes.
		\begin{itemize} 
			\item For all $(u,r) \in (0,T]^2,  s\leq \min (u,h),  (x,y)\in (\R^d)^2$, and each multi-index  $\zeta\in\N^d$ with length $|\zeta|\leq 1$, 
			\begin{align}
				|\nabla^\zeta p_\alpha (u,y-s\mathfrak b_h(r,x))| &\lesssim \frac{p_\alpha(u,y)}{u^{\frac{|\zeta|}{\alpha}}}.\label{ignore-bh-eq1}
			\end{align}
			\item For all $(u,r) \in (0,T]^2,  s\leq \min (u,h), (x,y,y')\in (\R^d)^3$, for each multi-index  $\zeta\in\N^d$ with length $|\zeta|\leq 1$, and for all $\delta \in (0,1]$,
			\begin{align}
				\left| \nabla^\zeta p_\alpha (u,y-s\mathfrak b_h(r,x))-\nabla^\zeta p_\alpha (u,y'-s\mathfrak b_h(r,x)) \right| &\lesssim \left( \frac{|y-y'|^\delta}{u^{\frac{\delta}{\alpha}}} \wedge 1\right) \frac{1}{u^{\frac{|\zeta|}{\alpha}}} \left(p_\alpha(u,y) +p_\alpha (u,y')\right).\label{ignore-bh-eq2}
			\end{align}
		\end{itemize}
		
	\end{lemma}
	\subsection{Proof of the error bounds of Theorem \ref{thm-main}}
	 Comparing the Duhamel formula of the scheme, \eqref{duhamel-scheme}, to that of the diffusion, \eqref{duhamel-diffusion}, we get
	\begin{align*}
	&	\Gamma^h(0,x,t,y)-\Gamma(0,x,t,y)\\ & =\int_0^t\int_{\R^d}\Gamma(0,x,s,z)b(s,z).\nabla_yp_\alpha(t-s,y-z)\d z \d s-\E_{0,x} \left[ \int_0^t b_h(U_{\lfloor \frac{s}{h}\rfloor},X_{\tau_s^h}^h)\cdot\nabla_y p_\alpha (t-s,y-X_s^h) \d s \right].
	\end{align*}
	Respectively, for the alternative scheme involving $\bar b_h$,
	\begin{align*}
		&	\bar \Gamma^h(0,x,t,y)-\Gamma(0,x,t,y)\\ & =\int_0^t\int_{\R^d}\Gamma(0,x,s,z)b(s,z).\nabla_yp_\alpha(t-s,y-z)\d z\d s-\E_{0,x} \left[ \int_0^t \bar b_h(U_{\lfloor \frac{s}{h}\rfloor},X_{\tau_s^h}^h)\cdot\nabla_y p_\alpha (t-s,y-X_s^h) \d s \right].
	\end{align*}
	
	The error admits the following decomposition:
	\begin{align}\label{splitting-error}
			&\nonumber\Gamma^h(0,x,t,y)-\Gamma(0,x,t,y)= \int_0^t  \int[\Gamma(0,x,s,z)-\Gamma^h(0,x,s,z)] b(s,z) \cdot \nabla_y p_\alpha(t-s,y-z) \d z \d s \\
			&\nonumber\qquad+\mathbb{1}_{\{t\ge 3h\}}\int_{t_1}^{\tau_t^h-h}\int\Gamma^h(0,x,s,z)(b(s,z)-b_h(s,z))\cdot\nabla_y p_\alpha(t-s,y-z) \d z \d s\\
			&\nonumber\qquad +\mathbb{1}_{\{t\ge 3h\}}\int_{t_1}^{\tau_t^{{h}}-h} \int[\Gamma^h(0,x,s,z)-\Gamma^h(0,x,\tau^h_s,z)]b_h(s,z)\cdot\nabla_y p_\alpha(t-s,y-z) \d z \d s\\
			&\nonumber\qquad +\mathbb{1}_{\{t\ge 3h\}}\int_{t_1}^{\tau_t^{{h}}-h}\E_{0,x}\bigg[b_h(U_{\lfloor s/h\rfloor},X_{\tau^h_s}^h)\cdot (\nabla_y p_\alpha(t-U_{\lfloor s/h\rfloor},y-X_{\tau^h_s}^h)-\nabla_y p_\alpha(t-s,y-X_s^h))\bigg] \d s\\
			&\nonumber\qquad+\frac{1}{h}\int_{0}^{t_1\wedge t}\int_{0}^{h}\int p_\alpha(s,z-x-b_h(r,x)s)\left(b(s,z)-b_h(r,x)\right)\cdot\nabla_y p_\alpha(t-s,y-z) \d z \d r \d s\\
			&\nonumber\qquad +\mathbb{1}_{\{t\ge h\}} \frac{1}{h}\int_{(\tau_t^{{h}}-h)\vee t_1}^{t}\int_{\tau^h_s}^{\tau^h_s+h}\int \int \Gamma^h(0,x,\tau^h_s,w) p_\alpha \left(s-\tau^h_s,z-w-b_h(r,w)(s-\tau^h_s)\right)\\
			&\nonumber\qquad\qquad\qquad\qquad\qquad\qquad\qquad\qquad\qquad\qquad \times (b(s,z)-b_h(r,w))\cdot\nabla_y p_\alpha (t-s,y-z) \d z \d w \d r \d s\\
			&=:\Delta_1+\Delta_2+\Delta_3+\Delta_4+\Delta_5+\Delta_6,
	\end{align}
	where, for the last term, we use that for $s\in(t_1,T]$, not belonging to the discretization grid and $\phi:\R^d\times \R^d\times \R \rightarrow \R$ measurable and bounded, since $X_s^h=X_{\tau_s^h}^h+Z_s-Z_{\tau_s^h}+b_h(U_{\lfloor \frac{s}{h}\rfloor},X_{\tau_s^h}^h)(s-\tau_s^h)$ with $X_{\tau_s^h}^h$, $Z_s-Z_{\tau_s^h}$ and $U_{\lfloor \frac{s}{h}\rfloor}$ independent, we can write
	\begin{align}
		&\E_{0,x} \left[\phi(X_{\tau_s^h}^h,X_s^h,U_{\lfloor \frac{s}{h}\rfloor})\right]\notag\\ & \qquad\qquad=\frac{1}{h}\int_{\tau_s^h}^{\tau_s^h+h}\int \int \phi(w,z,r)\Gamma^h(0,x,\tau_s^h,w)p_\alpha(s-\tau_s^h,z-w-b_h(r,w)(s-\tau_s^h)) \d z \d w \d r.\label{espxxu}
	\end{align}
	Similarly, we define  
	\begin{align*}
		\bar\Gamma^h(0,x,t,y)-\Gamma(0,x,t,y) =: \bar \Delta_1+\bar\Delta_2+\bar\Delta_3+\bar\Delta_4+\bar\Delta_5+\bar\Delta_6,
	\end{align*}
	where $b_h$ is replaced by $\bar{b}_h$.\\
	
	    For $\Delta_{2},\Delta_3,\Delta_4$, we suppose that $t\ge 3h$ (otherwise these contributions vanish) and rely on the fact that the current integration time is distinct from 0 and from $t$, meaning that we can rely on the smoothness properties of the integrands on the considered time intervals. For $\Delta_{5},\Delta_{6}$, on the opposite, we rely on the smallness of the considered time intervals.\\
	    
	\newpage
 
	\begin{figure}[h]

		\begin{centering}
		\begin{tikzpicture}[scale=6]
			\draw (1.5,0)--(1.5,0) node[anchor=north] {$t$};
			\draw (0,0.02)--(0,-0.02) node[anchor=north] {0};
			\draw (1.48,0.02)--(1.5,0) ;
			\draw (1.48,-0.02)--(1.5,0) ;		
			
			\node at (0.25,0.05) {$t_1$};
			\draw (0.25,0.01)--(0.25,-0.01) {};
			
			\node at (1.15,0.05) {$\tau_t^h-h$};
			\draw (1.15,0.01)--(1.15,-0.01) {};
			
			\node at (1.4,0.05) {$\tau_t^h$};
			\draw (1.4,0.01)--(1.4,-0.01) {};
			
			\draw (0,0)--(1.5,0);
			\draw (0,-0.2)--(1.5,-0.2);
			\draw (0,-0.2)--(0.02,-0.18);
			\draw (0,-0.2)--(0.02,-0.22);
			\draw (1.5,-0.2)--(1.48,-0.18);
			\draw (1.5,-0.2)--(1.48,-0.22);
			\node [draw=none] at (0.75,-0.15) {$\Delta_1$} ;
			\node [draw=none,anchor=west] at (1.6,-0.2) {Grönwall lemma} ;
			
			\draw (0.25,-0.4)--(1.15,-0.4);
			\draw (0.25,-0.4)--(0.27,-0.38);
			\draw (0.25,-0.4)--(0.27,-0.42);
			\draw (1.15,-0.4)--(1.13,-0.38);
			\draw (1.15,-0.4)--(1.13,-0.42);
			\node [draw=none] at (0.7,-0.35) {$\Delta_2$};
			\node [draw=none,anchor=west] at (1.6,-0.4) {Cutoff error terms} ;
			
			\draw (0.25,-0.6)--(1.15,-0.6);
			\draw (0.25,-0.6)--(0.27,-0.58);
			\draw (0.25,-0.6)--(0.27,-0.62);
			\draw (1.15,-0.6)--(1.13,-0.58);
			\draw (1.15,-0.6)--(1.13,-0.62);
			\node [draw=none] at (0.7,-0.55) {$\Delta_3$};
			\node [draw=none,anchor=west] at (1.6,-0.6) {Forward time regularity of $\Gamma^h$} ;
			
			\draw (0.25,-0.8)--(1.15,-0.8);
			\draw (0.25,-0.8)--(0.27,-0.78);
			\draw (0.25,-0.8)--(0.27,-0.82);
			\draw (1.15,-0.8)--(1.13,-0.78);
			\draw (1.15,-0.8)--(1.13,-0.82);
			\node [draw=none] at (0.7,-0.75) {$\Delta_4$};
			\node [draw=none,anchor=west] at (1.6,-0.8) {Stable sensitivities} ;
			
			\draw (0,-1)--(0.25,-1);
			\draw (0,-1)--(0.02,-0.98);
			\draw (0,-1)--(0.02,-1.02);
			\draw (0.25,-1)--(0.23,-0.98);
			\draw (0.25,-1)--(0.23,-1.02);
			\node [draw=none] at (0.125,-0.95) {$\Delta_5$};
			\node [draw=none,anchor=west] at (1.6,-1) {Overall error on the first full time step} ;
			
			\draw (1.15,-1.2)--(1.4,-1.2);
			\draw (1.15,-1.2)--(1.17,-1.18);
			\draw (1.15,-1.2)--(1.17,-1.22);
			\draw (1.4,-1.2)--(1.38,-1.18);
			\draw (1.4,-1.2)--(1.38,-1.22);
			\node [draw=none] at (1.275,-1.15) {$\Delta_6$};
			\node [draw=none,anchor=west] at (1.6,-1.2) {Overall error on the last full time step} ;
			\end{tikzpicture}\\[0.5cm]
		\end{centering}
		\caption{Splitting of the error}
		\end{figure}
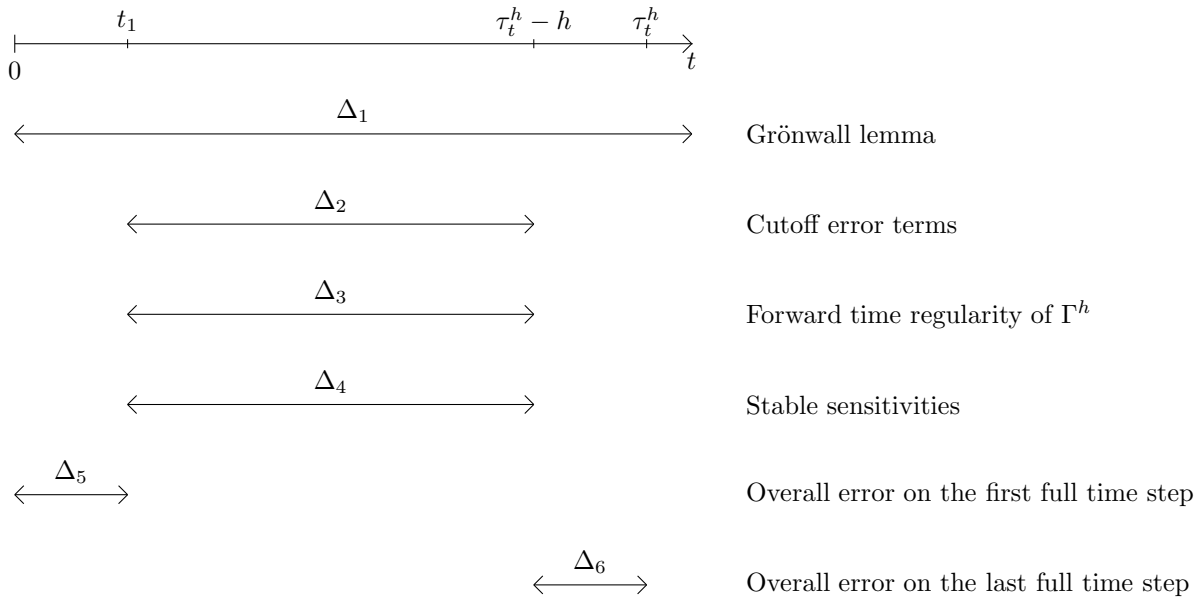

	 Let us first deal with $\Delta_2$. Since this term vanishes when $p=q=\infty$ (the same is true for $\bar\Delta_2$ when $h\le (\|b\|_{L^\infty-L^\infty}/B)^{\frac\alpha{1-\alpha}}$), we assume that either $p<\infty$ or $q<\infty$. Let $\lambda\ge 1$. Using the fact that $\forall y \in \R_+, \mathbb{1}_{\{ y\geq 1 \}} \leq y^{\lambda-1}$, we obtain that $\forall f:\R \rightarrow \R_+$, $\forall C>0$,
	$$f\mathbb{1}_{\{f\geq C\}} \leq  f^\lambda C^{1-\lambda} .$$
	This allows us to control the cutoff error in the following way:
	\begin{align}\label{maj-delta2-main}
		|b-b^h|=\left( |b|-Bh^{-\frac{d}{\alpha p}- \frac{1}{q}}\right)_+ \leq |b|\mathbb{1}_{|b|\geq Bh^{-\frac{d}{\alpha p}- \frac{1}{q}}} \leq |b|^\lambda B^{1-\lambda} h^{(\frac{d}{\alpha p}+\frac{1}{q})(\lambda-1)}.
	\end{align}
	Respectively
	\begin{align*}
		|b-\bar b^h|=\left( |b|-Bh^{\frac{1}{\alpha}-1}\right)_+ \leq |b|\mathbb{1}_{|b|\geq Bh^{\frac{1}{\alpha}-1}} \leq |b|^\lambda B^{1-\lambda} h^{(1-\frac{1}{\alpha})(\lambda-1)}.
	\end{align*}
	Along with the use of \eqref{ineq-density-scheme} and \eqref{derivatives-palpha}, we obtain
	\begin{align*}
					|\Delta_2| &= \left| \mathbb{1}_{\{t\ge 3h\}}\int_{t_1}^{\tau_t^h-h}\int\Gamma^h(0,x,s,z)(b(s,z)-b_h(s,z))\cdot\nabla_y p_\alpha(t-s,y-z) \d z \d s \right| \\
					&\lesssim  h^{-(\frac{d}{\alpha p}+\frac{1}{q})(1-\lambda)}\int_{t_1}^{\tau_t^h-h}\int p_\alpha (s,z-x)   |b(s,z)|^\lambda  \frac{p_\alpha(t-s,y-z)}{(t-s)^{\frac{1}{\alpha}}} \d z \d s.
	\end{align*}
	Let us check that we can choose 
	\begin{equation}\label{DEF_LAMBDA}
		{\lambda} = 1+ \frac{\frac{\gamma_1}{\alpha}}{\frac{d}{\alpha p}+\frac{1}{q}}=1+\frac{\gamma_1}{\frac{d}{p}+\frac{\alpha}{q}}\mbox{ with }\gamma_1 \in (\gamma,1],
              \end{equation}small enough so that 
              $\tilde p=\frac p\lambda$ and $\tilde q=\frac q\lambda$ satisfy $\tilde p\ge 1$ and $\tilde q\ge 1$. This is indeed possible since, by the definition \eqref{gap} of $\gamma$ and \eqref{serrin},
              $$\frac{p}{1+\frac{\gamma}{\frac{d}{p}+\frac{\alpha}{q}}}=\frac{d+\frac{\alpha p}q}{\alpha-1}>1\mbox{ and }\frac{q}{1+\frac{1}{\frac{d}{p}+\frac{\alpha}{q}}}>\frac{\frac{dq}p+\alpha}{\alpha}\ge 1.$$
              Morever, in order to estimate the time integrals that will appear below after the application of H\"older's inequality, let us observe that $\lambda>1$ and since $$\frac{q}{q-1-\frac{\gamma}{\frac{d}{p}+\frac{\alpha}{q}}}\left[\frac{1}{\alpha} + \frac{d}{\alpha p}\left(1+\frac{\gamma}{\frac{d}{p}+\frac{\alpha}{q}}\right)\right]=1
              ,$$
              $$\mbox{we have }-\left(\frac{q}{\lambda}\right)' \left[\frac{1}{\alpha} + \frac{d\lambda}{\alpha p}\right]=-\left(1-\frac{\lambda}{q}\right)^{-1}\left[\frac{1}{\alpha} + \frac{d\lambda}{\alpha p}\right]<-1.$$
    Using the identity $\forall f:\R \rightarrow \R, \forall \mu \in \R_+, \forall \mathfrak p\geq \mu, \Vert f^\mu \Vert_{L^{\frac{\mathfrak p}{\mu}}}=\Vert f \Vert_{L^{\mathfrak p}}^\mu$ and \eqref{convo-space-sing} with $\ell=\tilde p$ then Young's inequality,  H\"older's inequality in time and the last inequality combined with $(t-\tau_t^h+h)\geq h$, we get
	\begin{align}
		|\Delta_2| 	&\lesssim h^{(\frac{d}{\alpha p}+\frac{1}{q})({\lambda}-1)} \int_{t_1}^{\tau_t^h-h} \int p_\alpha (s,z-x)   |b(s,z)|^{{\lambda} } \frac{p_\alpha(t-s,y-z)}{(t-s)^{\frac{1}{\alpha}}} \d z \d s\notag\\
		&\lesssim  h^{(\frac{d}{\alpha p}+\frac{1}{q})({\lambda}-1)}p_\alpha (t,y-x)\int_{t_1}^{\tau_t^h-h}\frac{1}{(t-s)^{\frac{1}{\alpha}}} \Vert b(s,\cdot) \Vert_{L^{p}}^{{\lambda}} \left[\frac{1}{s^{\frac{d {\lambda}}{\alpha p}}}  +\frac{1}{(t-s)^{\frac{d {\lambda}}{\alpha p}}} \right] \d s\notag\\
	&\lesssim  h^{(\frac{d}{\alpha p}+\frac{1}{q})({\lambda}-1)}p_\alpha (t,y-x)\int_{t_1}^{\tau_t^h-h}\Vert b(s,\cdot) \Vert_{L^{p}}^{{\lambda}} \left[\frac{1}{s^{\frac{1}{\alpha}+\frac{d {\lambda}}{\alpha p}}}  +\frac{1}{(t-s)^{\frac{1}{\alpha}+\frac{d {\lambda}}{\alpha p}}} \right] \d s\notag\\	&\lesssim  h^{(\frac{d}{\alpha p}+\frac{1}{q})({\lambda}-1)}p_\alpha (t,y-x) \Vert s\mapsto \Vert b(s,\cdot) \Vert_{L^{p}}^{{\lambda}}  \Vert_{L^{\frac{q}{{\lambda}}}} \notag\\ & \qquad\qquad\qquad \times \left(\int_{t_1}^{\tau_t^h-h}\left[\frac{1}{s^{(\frac{q}{{\lambda}})'\left(\frac{1}{\alpha}+\frac{d {\lambda}}{\alpha p}\right)}}  +\frac{1}{(t-s)^{(\frac{q}{{\lambda}})'\left(\frac{1}{\alpha}+\frac{d {\lambda}}{\alpha p}\right)}} \right] \d s \right)^{\frac{1}{(\frac{q}{{\lambda}})'}}
 \notag\\&\lesssim  h^{(\frac{d}{\alpha p}+\frac{1}{q})({\lambda}-1)}p_\alpha (t,y-x)\left(h^{\frac{1}{(\frac{q}{{\lambda}})'}-\frac{1}{\alpha}-\frac{d {\lambda}}{\alpha p}}+(t-\tau_t^h+h)^{\frac{1}{(\frac{q}{{\lambda}})'}-\frac{1}{\alpha}-\frac{d {\lambda}}{\alpha p}}\right)\notag\\&\lesssim  h^{(\frac{d}{\alpha p}+\frac{1}{q})({\lambda}-1)+1-\frac{1}{\alpha}-{\lambda}\left(\frac{d}{\alpha p} + \frac{1}{q}\right)} p_\alpha (t,y-x)\notag\\&\lesssim h^{\frac{\gamma}{\alpha}}p_\alpha(t,y-x).\label{maj-D2-mainthm}	\end{align}
The same computations with the same choice of $\lambda$ yield $|\bar \Delta_2|\lesssim h^{(1-\frac{1}{\alpha})(\lambda-1)}h^{1-\frac{1}{\alpha}-{\lambda}\left(\frac{d}{\alpha p} + \frac{1}{q}\right)}p_\alpha(t,y-x) \lesssim h^{\frac{\lambda\gamma}{\alpha}}p_\alpha(t,y-x)\lesssim h^{\frac{\gamma}{\alpha}}p_\alpha(t,y-x)$.\\
	
	We now turn our attention to $\Delta_3$, for which we mainly rely on  the H\"older regularity of $\Gamma^h$ in time (equation \eqref{holder-time-gammah} of Proposition \ref{prop-main-estimates}). We assume $t\geq 3h$, otherwise this contribution vanishes. Using \eqref{holder-time-gammah}, we can write
	$$\forall s\ge t_1,\;|\Gamma^h(0,x,s,z)-\Gamma^h(0,x,\tau^h_s,z)| \lesssim \frac{(s-\tau_s^h)^{\frac{\gamma}{\alpha}}}{(\tau_s^h)^{\frac{\gamma}{\alpha}}}p_\alpha (s,z-x).$$
    We plug this into the definition of $\Delta_3$, using as well $|b_h|\leq |b|$ (resp. $|\bar b_h|\leq |b|$):
	\begin{align*}
		|\Delta_3 |&\lesssim \int_{t_1}^{\tau_t^{{h}}-h} \int\frac{(s-\tau_s^h)^{\frac{\gamma}{\alpha}}}{(\tau_s^h)^{\frac{\gamma}{\alpha}}}p_\alpha (s,z-x) |b(s,z)| |\nabla_y p_\alpha(t-s,y-z) |\d z \d s\\
		&\lesssim h^{\frac{\gamma}{\alpha}}\int_{t_1}^{\tau_t^{{h}}-h} \int\frac{1}{(\tau_s^h)^{\frac{\gamma}{\alpha}}}p_\alpha (s,z-x) |b(s,z)| |\nabla_y p_\alpha(t-s,y-z) |\d z \d s.
	\end{align*}
	We now deal with the gradient using \eqref{derivatives-palpha} and notice that since $s\ge t_1\geq h$ then $(\tau_s^h)^{-1} \leq 2 s^{-1}$ to write
	\begin{align*}
		|\Delta_3 |&\lesssim h^{\frac{\gamma}{\alpha}} \int_{t_1}^{\tau_t^{{h}}-h} \int\frac{1}{(t-s)^{\frac{1}{\alpha}}s^{\frac{\gamma}{\alpha}}} p_\alpha (s,z-x) |b(s,z)| p_\alpha(t-s,y-z) \d z \d s\\
		&\lesssim h^{\frac{\gamma}{\alpha}} \int_{0}^{t} \int\frac{1}{(t-s)^{\frac{1}{\alpha}}s^{\frac{\gamma}{\alpha}}} p_\alpha (s,z-x) |b(s,z)| p_\alpha(t-s,y-z) \d z \d s.
	\end{align*}
	Using \eqref{lemma-convo-bulk-2} (with $u=0,v=t, \beta_1=1/\alpha, \beta_2=\gamma/\alpha$ so that $\frac{\gamma+1}{\alpha}-(\beta_1+\beta_2)=0$), we obtain
	\begin{align}\label{maj-D3-mainthm}
		|\Delta_3 |\lesssim h^{\frac{\gamma}{\alpha}} p_\alpha (t,y-x).
	\end{align}

 For $\Delta_4$, we first expand the expectation with the known densities using \eqref{espxxu}:
	\begin{align*}
		\Delta_4 &= \int_{t_1}^{\tau_t^{{h}}-h}\E_{0,x}\bigg[b_h(U_{\lfloor s/h\rfloor},X_{\tau^h_s}^h)\cdot (\nabla p_\alpha(t-U_{\lfloor s/h\rfloor},y-X_{\tau^h_s}^h)-\nabla p_\alpha(t-s,y-X_s^h))\bigg] \d s\\
		&= \sum_{j=1}^{\lfloor \frac{t}{h}\rfloor -2} \int_{t_j}^{t_{j+1}}  \Gamma^h (0,x,t_j,z)  \frac{1}{h}\int_{t_j}^{t_{j+1}} \int \int p_\alpha (s-t_j,w-z-b_h(r,z)(s-t_j)) \\
		&\qquad\qquad\qquad\qquad\qquad\qquad\qquad \times b_h(r,z)\cdot (\nabla p_\alpha(t-r,y-z)-\nabla p_\alpha(t-s,y-w)) \d z \d w \d r \d s.
	\end{align*}
	We then derive, using \eqref{ineq-density-scheme} and \eqref{ignore-bh-eq1},
	\begin{align}\label{delta4-main-proof}
		&|\Delta_4| \lesssim \sum_{j=1}^{\lfloor \frac{t}{h}\rfloor -2}  \int_{t_j}^{t_{j+1}}    \frac{1}{h}\int_{t_j}^{t_{j+1}} \int \int \Gamma^h (0,x,t_j,z) p_\alpha (s-t_j,w-z-b_h(r,z)(s-t_j))|b(r,z)| \nonumber  \\
		&\qquad \times \left( \left|\nabla p_\alpha(t-r,y-z) -\nabla p_\alpha(t-s,y-z) \right| + \left| \nabla p_\alpha(t-s,y-z) -\nabla p_\alpha(t-s,y-w) \right| \right) \d z \d w \d r \d s \nonumber\\
		&\qquad \lesssim \sum_{j=1}^{\lfloor \frac{t}{h}\rfloor -2} \int_{t_j}^{t_{j+1}}    \frac{1}{h}\int_{t_j}^{t_{j+1}} \int \int p_\alpha (t_j,z-x)  p_\alpha (s-t_j,w-z)|b(r,z)| \nonumber \\
		&\qquad \times \left( \left|\nabla p_\alpha(t-r,y-z) -\nabla p_\alpha(t-s,y-z) \right| + \left| \nabla p_\alpha(t-s,y-z) -\nabla p_\alpha(t-s,y-w) \right| \right) \d z \d w \d r \d s.
	\end{align}
	Next, we use \eqref{holder-time-palpha} to write 
	\begin{align*}
		\left|\nabla p_\alpha(t-r,y-z) -\nabla p_\alpha(t-s,y-z) \right| & \lesssim \frac{|r-s|}{(t-r\vee s)^{1+\frac{1}{\alpha}}} \left(p_\alpha (t-r,y-z) +p_\alpha (t-s,y-z)\right).
	\end{align*}
	Since $r-s<h$, $t\geq 3h$ and $t-r\vee s \ge h$, we can use \eqref{ARONSON_STABLE} to deduce that
	\begin{align*}
	\left|\nabla p_\alpha(t-r,y-z) -\nabla p_\alpha(t-s,y-z) \right| & \lesssim \frac{|r-s|}{(t-r)^{1+\frac{1}{\alpha}}} p_\alpha (t-r,y-z),
	\end{align*}
	which also yields, for any $\gamma_1\in (\gamma,\alpha]$, recalling that $r-s<t-r$,
	\begin{align}
		\left|\nabla p_\alpha(t-r,y-z) -\nabla p_\alpha(t-s,y-z) \right| & \lesssim \frac{|r-s|^{\frac{\gamma_1}{\alpha}}}{(t-r)^{\frac{\gamma_1}{\alpha}+\frac{{1}}{\alpha}}} p_\alpha (t-r,y-z). \label{double-maj-2}
	\end{align}
	\color{black} For the second term in \eqref{delta4-main-proof}, assuming that $\gamma_1\in (\gamma,1]$ \color{black}, we deduce from \eqref{holder-space-palpha} that
	\begin{align}
		\left| \nabla p_\alpha(t-s,y-z) -\nabla p_\alpha(t-s,y-w) \right| & \lesssim \left(\frac{|z-w|^{\gamma_1}}{(t-s)^{\frac{\gamma_1}{\alpha}}}\wedge 1\right) \frac{1}{(t-s)^{\frac{1}{\alpha}}}\left( p_\alpha(t-s,y-z) + p_\alpha(t-s,y-w)\right). \label{double-maj-1}
	\end{align}
	Plugging \eqref{double-maj-1} and \eqref{double-maj-2} into \eqref{delta4-main-proof}, we can write 
	\begin{align}
		& \Delta_4 \lesssim\sum_{j=1}^{\lfloor \frac{t}{h}\rfloor -2} \int_{t_j}^{t_{j+1}}    \frac{1}{h}\int_{t_j}^{t_{j+1}} \int \int p_\alpha (t_j,z-x)  p_\alpha (s-t_j,w-z)|b(r,z)| p_\alpha (t-r,y-z) \frac{ |r-s|^{\frac{\gamma_1}{\alpha}}}{(t-r)^{\frac{{1+\gamma_1 }}{\alpha}}} \d z \d w \d r \d s \nonumber\\
		& \qquad+ \sum_{j=1}^{\lfloor \frac{t}{h}\rfloor -2} \int_{t_j}^{t_{j+1}}    \frac{1}{h}\int_{t_j}^{t_{j+1}} \int \int p_\alpha (t_j,z-x)  p_\alpha (s-t_j,w-z)|b(r,z)| p_\alpha (t-s,y-z) \frac{|z-w|^{\gamma_1}}{(t-s)^{\frac{1+{\gamma_1}}{\alpha}}}\d z \d w \d r \d s\nonumber\\
		& \qquad+ \sum_{j=1}^{\lfloor \frac{t}{h}\rfloor -2} \int_{t_j}^{t_{j+1}}    \frac{1}{h}\int_{t_j}^{t_{j+1}} \int \int p_\alpha (t_j,z-x)  p_\alpha (s-t_j,w-z)|b(r,z)|\nonumber\\ & \qquad\qquad\qquad\qquad\qquad \times p_\alpha (t-s,y-w) \left(\frac{|z-w|^{\gamma_1}}{(t-s)^{\frac{\gamma_1}{\alpha}}}\wedge 1\right) \frac{1}{(t-s)^{\frac{1}{\alpha}}} \d z \d w \d r \d s\nonumber\\
		&=:\Delta_4^1+\Delta_4^2+\Delta_4^3. \label{def-delta-4-i}
	\end{align}
	Let us treat $\Delta_4^1$. From the Fubini theorem, we integrate first in $w$ using the fact that $p_\alpha$ is a probability density:
	\begin{align}\label{d4-inter}
		\Delta_4^1 =\sum_{j=1}^{\lfloor \frac{t}{h}\rfloor -2} \int_{t_j}^{t_{j+1}}    \frac{1}{h}\int_{t_j}^{t_{j+1}} \int  p_\alpha (t_j,z-x)  |b(r,z)| p_\alpha (t-r,y-z) \frac{ |r-s|^{\frac{\gamma_1}{\alpha}}}{(t-r)^{\frac{1+\gamma_1}{\alpha}}} \d z \d r \d s.
	\end{align}
 \color{black}
	Then, using $t_j^{-1}\le 2r^{-1} $ and \eqref{ARONSON_STABLE} along with the fact that $|r-s|<h$, we get
		\begin{align*}
		\Delta_4^1& \lesssim\sum_{j=1}^{\lfloor \frac{t}{h}\rfloor -2} \int_{t_j}^{t_{j+1}}   h^{\frac{\gamma_1}{\alpha}-1} \int_{t_j}^{t_{j+1}} \int  p_\alpha (r,z-x)  |b(r,z)| p_\alpha (t-r,y-z) \frac{ 1}{(t-r)^{\frac{{1}+\gamma_1}{\alpha}}} \d z \d r \d s\\
		& \lesssim\int_{t_1}^{t_{\lfloor \frac{t}{h}\rfloor -1}}   h^{\frac{\gamma_1}{\alpha}} \int  p_\alpha (r,z-x)  |b(r,z)| p_\alpha (t-r,y-z) \frac{ 1}{(t-r)^{\frac{{1}+\gamma_1}{\alpha}}} \d z \d r.
		\end{align*}
		Using \eqref{lemma-convo-bulk-1} (singular case with $u=t_1,v=\tau_t^h-h, \beta_1=(1+\gamma_1)/\alpha,\beta_2=0$ and noting that since $t\ge 3h$, $v-u\ge h$), we get
			\begin{align}
			\Delta_4^1	& \lesssim p_\alpha(t,y-x)  h^{\frac{\gamma_1}{\alpha}} \left(h^{\frac{\gamma-\gamma_1}{\alpha}}+(t-t_{\lfloor \frac{t}{h}\rfloor -1})^{\frac{\gamma-\gamma_1}{\alpha}}\right)\notag\\
			&\lesssim p_\alpha(t,y-x) h^{\frac{\gamma}{\alpha}}.	\label{MAJ_DELTA_41}
		\end{align}
	Let us treat $\Delta_4^2$ defined in \eqref{def-delta-4-i}. We integrate in $w$ using \eqref{spatial-moments} and use the fact that $s-t_j\leq h$:
	\begin{align*}
		\Delta_4^2&=\sum_{j=1}^{\lfloor \frac{t}{h}\rfloor -2} \int_{t_j}^{t_{j+1}}    \frac{1}{h}\int_{t_j}^{t_{j+1}} \int \int p_\alpha (t_j,z-x)  p_\alpha (s-t_j,w-z)|b(r,z)| p_\alpha (t-s,y-z) \frac{|z-w|^{\gamma_1}}{(t-s)^{\frac{1+{\gamma_1}}{\alpha}}}\d z \d w \d r \d s\\
		&\lesssim\sum_{j=1}^{\lfloor \frac{t}{h}\rfloor -2} \int_{t_j}^{t_{j+1}}    \frac{1}{h}\int_{t_j}^{t_{j+1}}  \int p_\alpha (t_j,z-x) |b(r,z)| p_\alpha (t-s,y-z) \frac{(s-t_j)^{\frac{\gamma_1}{\alpha}}}{(t-s)^{\frac{1+\gamma_1}{\alpha}}}\d z  \d r \d s\\
		&\lesssim h^{\frac{\gamma_1}{\alpha}}\sum_{j=1}^{\lfloor \frac{t}{h}\rfloor -2} \int_{t_j}^{t_{j+1}}    \frac{1}{h}\int_{t_j}^{t_{j+1}}  \int p_\alpha (t_j,z-x) |b(r,z)| p_\alpha (t-s,y-z) \frac{1}{(t-s)^{\frac{1+\gamma_1}{\alpha}}}\d z  \d r \d s.
	\end{align*}
	Notice that in the previous integral, as above $p_\alpha (t_j,z-x)\lesssim p_\alpha (r,z-x)$ and $p_\alpha (t-s,y-z)\lesssim p_\alpha (t-r,y-z)$, $\frac{1}{(t-s)^{\frac{1+\gamma_1}{\alpha}}}\lesssim \frac{1}{(t-r)^{\frac{1+\gamma_1}{\alpha}}}$. This yields
		\begin{align*}
		\Delta_4^2&\lesssim h^{\frac{\gamma_1}{\alpha}} \int_{t_1}^{t_{\lfloor \frac{t}{h}\rfloor -1}}    \int p_\alpha (r,z-x) |b(r,z)| p_\alpha (t-r,y-z) \frac{1}{(t-r)^{\frac{1+\gamma_1}{\alpha}}}\d z  \d r,
	\end{align*}
	which is the right form to use \eqref{lemma-convo-bulk-1} with the same parameters as for $\Delta_4^1$. Doing so, we obtain similarly
	\begin{align}
		\Delta_4^2	&\lesssim p_\alpha (t,y-x)h^{\frac{\gamma_1}{\alpha}} (t-t_{\lfloor \frac{t}{h}\rfloor -1})^{\frac{\gamma-\gamma_1}{\alpha}}\notag\\
		&\lesssim p_\alpha (t,y-x)h^{\frac{\gamma}{\alpha}}.\label{MAJ_DELTA_42}
	\end{align}
Let us now turn to the term  $\Delta_4^3$ in \eqref{def-delta-4-i}:
	\begin{itemize}
		\item \textbf{Global off-diagonal case:} $|x-y|>t^{\frac{1}{\alpha}}$, then, since $|y-x|\leq |y-w|+|z-w|+|z-x|$, at least one of the stable transitions in $\Delta_4^3$ will be off-diagonal as well. In this case, we will actually manage to retrieve the global final regime for $p_\alpha(t,y-x) $ from the inner densities in $\Delta_4^3$.
		\begin{itemize}
			\item 		If  $|z-x|>\frac{1}{3}|x-y|>\frac{1}{3}t^{\frac{1}{\alpha}}\gtrsim t_j^{\frac{1}{\alpha}}$, we can write
			\begin{align*}
				p_\alpha(t_j,z-x)\lesssim \frac{t_j}{|z-x|^{d+\alpha}} \lesssim\frac{t}{|x-y|^{d+\alpha}} \lesssim p_\alpha (t,x-y).
			\end{align*}
			We can then compute
			\begin{align}\label{maj-d4-od-1}
				\Delta_4^{3,1}&:=\mathbb{1}_{|x-y|>t^{\frac{1}{\alpha}}}\sum_{j=1}^{\lfloor \frac{t}{h}\rfloor -2} \int_{t_j}^{t_{j+1}}    \frac{1}{h}\int_{t_j}^{t_{j+1}} \int \int \mathbb{1}_{|z-x|>\frac{1}{3}|x-y|}p_\alpha (t_j,z-x) \nonumber\\  & \qquad \times p_\alpha (s-t_j,w-z)|b(r,z)|p_\alpha (t-s,y-w) \left(\frac{|z-w|^{\gamma_1}}{(t-s)^{\frac{\gamma_1}{\alpha}}}\wedge 1\right) \frac{1}{(t-s)^{\frac{1}{\alpha}}}\d z \d w \d r \d s\nonumber\\
				&\lesssim p_\alpha(t,y-x)\sum_{j=1}^{\lfloor \frac{t}{h}\rfloor -2} \int_{t_j}^{t_{j+1}}    \frac{1}{h}\int_{t_j}^{t_{j+1}} \int \int p_\alpha (s-t_j,w-z)\\  & \qquad\qquad \times|b(r,z)|p_\alpha (t-s,y-w) \frac{|z-w|^{\gamma_1}}{(t-s)^{\frac{1+{\gamma_1}}{\alpha}}} \d z \d w \d r \d s \nonumber \\
				&\lesssim p_\alpha(t,y-x)\sum_{j=1}^{\lfloor \frac{t}{h}\rfloor -2} \int_{t_j}^{t_{j+1}}    \frac{1}{h}\int_{t_j}^{t_{j+1}} \Vert b(r,\cdot)\Vert_{L^p} \d r \nonumber\\  & \qquad\qquad \times\int \Vert p_\alpha (s-t_j,w-\cdot)|\cdot-w|^{\gamma_1}\Vert_{L^{p'}}\frac{p_\alpha (t-s,y-w) }{(t-s)^{\frac{1+{\gamma_1}}{\alpha}}} \d w \d s.\nonumber
			\end{align}	
			Note that $\gamma_1\leq d+\alpha-\frac{d}{p'}$, allowing us to use \eqref{spatial-moments}:
			$$\Vert p_\alpha (s-t_j,w-\cdot)|\cdot-w|^{\gamma_1}\Vert_{L^{p'}} \lesssim (s-t_j)^{\frac{\gamma_1}{\alpha}-\frac{d}{\alpha p}},$$
			yielding, once integrating in $w$,
			\begin{align}\label{maj-d4-31}
				\Delta_4^{3,1}&\lesssim p_\alpha(t,y-x)\sum_{j=1}^{\lfloor \frac{t}{h}\rfloor -2} \int_{t_j}^{t_{j+1}}    \frac{1}{h}\int_{t_j}^{t_{j+1}} \Vert b(r,\cdot)\Vert_{L^p} \d r \frac{(s-t_j)^{\frac{\gamma_1}{\alpha}-\frac{d}{\alpha p}} }{(t-s)^{\frac{1+{\gamma_1}}{\alpha}}}  \d s\nonumber\\
				&\lesssim p_\alpha(t,y-x) \int_{t_1}^{t_{\lfloor \frac{t}{h}\rfloor -1}}  h^{\frac{\gamma_1}{\alpha}-\frac{d}{\alpha p}-\frac{1}{q}} \frac{1 }{(t-s)^{\frac{1+{\gamma_1}}{\alpha}}}  \d s\nonumber\\
				&\lesssim p_\alpha(t,y-x) h^{\frac{\gamma_1}{\alpha}-\frac{d}{\alpha p}-\frac{1}{q}}h^{1-\frac{1}{\alpha}-\frac{\gamma_1}{\alpha}}\nonumber\\
				&\lesssim p_\alpha(t,y-x) h^{\frac{\gamma}{\alpha}},
			\end{align}	
			the last inequality being true only if $1-\frac{1}{\alpha}-\frac{\gamma_1}{\alpha}<0$, which is always possible to satisfy since the choice of $\gamma_1\in (\gamma,1]$ is free.
			\item If $|z-w|>\frac{1}{3}|x-y|>\frac{1}{3}t^{\frac{1}{\alpha}}$, remarking that $s-t_j\leq h$ and $0<\frac{\gamma}{\alpha}+\frac{1}{q}<1$, we can write
			$$p_\alpha(s-t_j,w-z)\lesssim \frac{s-t_j}{|w-z|^{d+\alpha}} \lesssim \frac{s-t_j}{t}\times\frac{t}{|x-y|^{d+\alpha}} \lesssim \frac{s-t_j}{t}p_\alpha (t,x-y)\lesssim  \frac{h^{\frac{\gamma}{\alpha}+\frac{1}{q}}}{t^{\frac{\gamma}{\alpha}+\frac{1}{q}}}p_\alpha (t,x-y),$$
			and then compute
			\begin{align*}
				\Delta_4^{3,2}&:=\mathbb{1}_{|x-y|>t^{\frac{1}{\alpha}}}\sum_{j=1}^{\lfloor \frac{t}{h}\rfloor -2} \int_{t_j}^{t_{j+1}}    \frac{1}{h}\int_{t_j}^{t_{j+1}} \int \int \mathbb{1}_{|z-w|>\frac{1}{3}|x-y|}p_\alpha (t_j,z-x) \\  & \qquad \times p_\alpha (s-t_j,w-z)|b(r,z)|p_\alpha (t-s,y-w) \left(\frac{|z-w|^{\gamma_1}}{(t-s)^{\frac{\gamma_1}{\alpha}}}\wedge 1\right) \frac{1}{(t-s)^{\frac{1}{\alpha}}} \d z \d w \d r \d s\\
				&\lesssim p_\alpha(t,y-x)\frac{h^{\frac{\gamma}{\alpha}+\frac{1}{q}-1}}{t^{\frac{\gamma}{\alpha}+\frac{1}{q}}}\sum_{j=1}^{\lfloor \frac{t}{h}\rfloor -2} \int_{t_j}^{t_{j+1}}   \int_{t_j}^{t_{j+1}} \int \int p_\alpha (t_j,z-x)\\  & \qquad \times |b(r,z)|p_\alpha (t-s,y-w) \frac{1}{(t-s)^{\frac{1}{\alpha}}} \d z \d w \d r \d s\\
				&\lesssim p_\alpha(t,y-x)\frac{h^{\frac{\gamma}{\alpha}+\frac{1}{q}-1}}{t^{\frac{\gamma}{\alpha}+\frac{1}{q}}}\sum_{j=1}^{\lfloor \frac{t}{h}\rfloor -2}  \int_{t_j}^{t_{j+1}} \int_{t_j}^{t_{j+1}} \int  p_\alpha (t_j,z-x) |b(r,z)|\frac{1}{(t-s)^{\frac{1}{\alpha}}} \d z  \d r \d s\\
				&\lesssim p_\alpha(t,y-x)\frac{h^{\frac{\gamma}{\alpha}+\frac{1}{q}-1}}{t^{\frac{\gamma}{\alpha}+\frac{1}{q}}}\sum_{j=1}^{\lfloor \frac{t}{h}\rfloor -2}   \int_{t_j}^{t_{j+1}}\Vert b(r,\cdot)\Vert_{L^p} \d r \int_{t_j}^{t_{j+1}}  \Vert  p_\alpha (t_j,\cdot-x)\Vert_{L^{p'}}\frac{1}{(t-s)^{\frac{1}{\alpha}}}  \d s\\
				&\lesssim p_\alpha(t,y-x)\frac{h^{\frac{\gamma}{\alpha}}}{t^{\frac{\gamma}{\alpha}+\frac{1}{q}}}\sum_{j=1}^{\lfloor \frac{t}{h}\rfloor -2}   \int_{t_j}^{t_{j+1}}  \frac{1}{t_j^{\frac{d}{\alpha p}}}\frac{1}{(t-s)^{\frac{1}{\alpha}}}   \d s,
			\end{align*}
			using \eqref{spatial-moments} (with $\delta=0$) and the H\"older inequality in space for the \textcolor{black}{antepenultimate} inequality and the H\"older inequality in time for the last one.\\
			Next, remarking that $h\leq t_j$ and therefore $ t_j^{-1}\le 2s^{-1}$, we can write
			\begin{align*}
				\Delta_4^{3,2}&\lesssim p_\alpha(t,y-x)\frac{h^{\frac{\gamma}{\alpha}}}{t^{\frac{\gamma}{\alpha}+\frac{1}{q}}}\int_{0}^{t}  \frac{1}{s^{\frac{d}{\alpha p}}}\frac{1}{(t-s)^{\frac{1}{\alpha}}}   \d s.
			\end{align*}
			Hence,
			\begin{align}\label{maj-d4-32}
				\Delta_4^{3,2}&\lesssim p_\alpha(t,y-x)\frac{h^{\frac{\gamma}{\alpha}}}{t^{\frac{\gamma}{\alpha}+\frac{1}{q}}}t^{1-\frac{d}{\alpha p}-\frac{1}{\alpha}}\int_{0}^{1}  \frac{1}{\lambda^{\frac{d}{\alpha p}}}\frac{1}{(1-\lambda)^{\frac{1}{\alpha}}}   \d \lambda\nonumber\\
				&\lesssim p_\alpha(t,y-x)h^{\frac{\gamma}{\alpha}},
			\end{align}	
            recalling the definition of $\gamma$ in \eqref{gap} for the last inequality.
			\item If $|y-w|>\frac{1}{3}|x-y|>\frac{1}{3}t^{\frac{1}{\alpha}}\gtrsim (t-s)^{\frac{1}{\alpha}}$, we can write
			$$p_\alpha(t-s,y-w)\lesssim \frac{t-s}{|y-w|^{d+\alpha}} \lesssim \frac{t}{|x-y|^{d+\alpha}} \lesssim p_\alpha (t,x-y).$$
			This yields, using \eqref{spatial-moments} to bound $\Vert p_\alpha (s-t_j,z-\cdot)|z-\cdot|^{\gamma_1}\Vert_{L^1}$,
			\begin{align}\label{maj-d4-od-3}
				\nonumber\Delta_4^{3,3}&:= \mathbb{1}_{|x-y|>t^{\frac{1}{\alpha}}}\sum_{j=1}^{\lfloor \frac{t}{h}\rfloor -2} \int_{t_j}^{t_{j+1}}    \frac{1}{h}\int_{t_j}^{t_{j+1}} \int \int \mathbb{1}_{|y-w|>\frac{1}{3}|x-y|}p_\alpha (t_j,z-x) \\ \nonumber  & \qquad \times p_\alpha (s-t_j,w-z)|b(r,z)|p_\alpha (t-s,y-w) \left(\frac{|z-w|^{\gamma_1}}{(t-s)^{\frac{\gamma_1}{\alpha}}}\wedge 1\right) \frac{1}{(t-s)^{\frac{1}{\alpha}}} \d z \d w \d r \d s\\ \nonumber
				&\lesssim p_\alpha(t,y-x)\sum_{j=1}^{\lfloor \frac{t}{h}\rfloor -2} \int_{t_j}^{t_{j+1}}    \frac{1}{h}\int_{t_j}^{t_{j+1}} \int \int p_\alpha (t_j,z-x)p_\alpha (s-t_j,w-z) \\ & \qquad\qquad\qquad\qquad\qquad\qquad\qquad\qquad\qquad\qquad \times |b(r,z)| \frac{|z-w|^{\gamma_1}}{(t-s)^{\frac{1+{\gamma_1}}{\alpha}}} \d z \d w \d r \d s\\
				&\lesssim p_\alpha(t,y-x)\sum_{j=1}^{\lfloor \frac{t}{h}\rfloor -2} \int_{t_j}^{t_{j+1}}    \frac{1}{h}\int_{t_j}^{t_{j+1}} \frac{(s-t_j)^{\frac{\gamma_1}{\alpha}}}{(t-s)^{\frac{1+{\gamma_1}}{\alpha}}} \int  p_\alpha (t_j,z-x)|b(r,z)| \d z \d r \d s\nonumber.
			\end{align}
			Next, we use $s-t_j \leq h$, a H\"older inequality in $z$ and \eqref{spatial-moments}:
			\begin{align*}
				\Delta_4^{3,3}&\lesssim p_\alpha(t,y-x)h^{\frac{\gamma_1}{\alpha}-1}\sum_{j=1}^{\lfloor \frac{t}{h}\rfloor -2} \int_{t_j}^{t_{j+1}}\Vert b(r,\cdot)\Vert_{L^p} \d r \int_{t_j}^{t_{j+1}}\Vert p_\alpha (t_j,\cdot-x)\Vert_{L^{p'}} \frac{1}{(t-s)^{\frac{1+{\gamma_1}}{\alpha}}}  \d s\\
				&\lesssim p_\alpha(t,y-x)h^{\frac{\gamma_1}{\alpha}-\frac{1}{q}}\sum_{j=1}^{\lfloor \frac{t}{h}\rfloor -2}  \int_{t_j}^{t_{j+1}} \frac{1}{t_j^{\frac{d}{\alpha p}}}\frac{1}{(t-s)^{\frac{1+{\gamma_1}}{\alpha}}}  \d s\\&\lesssim p_\alpha(t,y-x)h^{\frac{\gamma_1}{\alpha}-\frac{1}{q}-\frac{d}{\alpha p}}\int_{h}^{\tau^h_t-h} \frac{1}{(t-s)^{\frac{1+{\gamma_1}}{\alpha}}}  \d s.
			\end{align*}
			Choosing $\gamma_1\in (\gamma,1]$ such that $\frac{1+{\gamma_1}}{\alpha}>1$ we conclude that
			\begin{equation}
   				\Delta_4^{3,3}\lesssim p_\alpha(t,y-x)h^{\frac{\gamma_1}{\alpha}-\frac{1}{q}-\frac{d}{\alpha p}}(t-\tau^h_t+h)^{1-\frac{1+{\gamma_1}}{\alpha}}\lesssim p_\alpha(t,y-x)h^{\frac{\gamma}{\alpha}}.\label{maj-d4-33}
			\end{equation}	
		\end{itemize}
   \color{black}

		\item \textbf{Global diagonal case:}  $|x-y|<t^{\frac{1}{\alpha}}$. We will use the fact that $p_\alpha (t,y-x)\asymp t^{-\frac{d}{\alpha}}$ to replace one of the local transitions with $p_\alpha(t,y-x)$, and then the computations will be the same as in the global off-diagonal case:
	\begin{itemize}
		\item if $t_j<t/2$, $p_\alpha (t-s,y-w)\lesssim(t-s)^{-\frac{d}{\alpha}} \lesssim t^{-\frac{d}{\alpha}} \asymp p_\alpha (t,y-x)$, and the computations are the same as from \eqref{maj-d4-od-3},
		\item if $t_j\geq t/2$, $p_\alpha(t_j,z-x)\lesssim t_j^{-\frac{d}{\alpha}} \lesssim t^{-\frac{d}{\alpha}} \asymp p_\alpha (t,y-x)$, and the computations are the same as from \eqref{maj-d4-od-1}.
	\end{itemize}
	\end{itemize}
   \color{black}
	Overall, gathering the estimates \eqref{maj-d4-31}, \eqref{maj-d4-32} and \eqref{maj-d4-33} as well as the estimates from the global diagonal case, we obtain $	\Delta_4^3 \lesssim p_\alpha(t,y-x)h^{\frac{\gamma}{\alpha}},$ which together with \eqref{MAJ_DELTA_42}, \eqref{MAJ_DELTA_41} and \eqref{def-delta-4-i} eventually yields
	\begin{align}
	\label{MAJ_DELTA_4}
	\Delta_4 \lesssim p_\alpha(t,y-x)h^{\frac{\gamma}{\alpha}},
	\end{align}
	 as intended. As we only used $|b_h|\lesssim b$ for $\Delta_3$ and $\Delta_4$, the estimations remain valid for $\bar \Delta_3$ and $\bar \Delta_4$.\\
	
	Let us turn our attention to $\Delta_5$ in \eqref{splitting-error} (first time step). Note that, even though a term $b(s,z)-b_h(r,x)$ appears in $\Delta_5$, its smallness actually follows from the fact that it only covers the first time step (over $(0,t_1\wedge t)$). Thus, we will bound $\Delta_5$ using the triangular inequality $|b(s,z)-b_h(r,x)|\lesssim |b(s,z)|+h^{-\frac{d}{\alpha p}-\frac{1}{q}}$ (resp. using $|\bar b_h (r,x)|=0$ for $r\leq h$), and then compute a bound for each term. Namely,
	\begin{align*}
		|\Delta_5| &= \left| \frac{1}{h}\int_{0}^{t_1\wedge t}\int_{0}^{h}\int p_\alpha(s,z-x-b_h(r,x)s) ( b(s,z)-b_h(r,x)) \cdot\nabla_y p_\alpha(t-s,y-z) \d z \d r \d s  \right|\\
		& \lesssim \frac{1}{h}\int_{0}^{t_1\wedge t}\int_{0}^{h}\int p_\alpha(s,z-x-b_h(r,x)s)\left( |b(s,z)|+ h^{-\frac{d}{\alpha p}-\frac{1}{q}} \right) \frac{1}{(t-s)^{\frac{1}{\alpha}}} p_\alpha(t-s,y-z) \d z \d r \d s.
	\end{align*}
	Since in our current integrals, using \eqref{ignore-bh-eq1}, $p_\alpha(s,z-x-b_h(r,x)s) \lesssim p_\alpha(s,z-x)$, we can write 
	\begin{align*}
		|\Delta_5| 	& \lesssim \int_{0}^{t_1\wedge t}\int p_\alpha(s,z-x) |b(s,z)| \frac{1}{(t-s)^{\frac{1}{\alpha}}} p_\alpha(t-s,y-z) \d z \d s\\
		& \qquad + \int_{0}^{t_1\wedge t}\int p_\alpha(s,z-x) h^{-\frac{d}{\alpha p}- \frac{1}{q}}\frac{1}{(t-s)^{\frac{1}{\alpha}}} p_\alpha(t-s,y-z) \d z  \d s.
	\end{align*}
	We then use \eqref{lemma-convo-bulk-2} with $u=0$, $v=t_1\wedge t$, $\beta_1=\frac 1\alpha$ and $\beta_2=0$ for the first term in the right-hand side and the convolution properties of the stable kernel for the second one to conclude that :
	\begin{align}\label{maj-D5-mainthm}
		|\Delta_5|	
		& \lesssim p_\alpha(t,y-x) h^{\frac{\gamma}{\alpha}}.
	\end{align}
	\begin{align}\label{maj-BAR-D5-mainthm}
		\mathrm{resp.} \qquad |\bar \Delta_5|\lesssim p_\alpha(t,y-x) h^{\frac{\gamma}{\alpha}}.
	\end{align}
	
	Let us now turn to $\Delta_6$ in \eqref{splitting-error}, for which the same reasoning as for $\Delta_5$ applies, although this time we are working on the last time step, over $((\tau_t^{{h}}-h)\vee t_1,t)$. Let $t\geq h$ (otherwise, $\Delta_6$ vanishes). Using \eqref{ineq-density-scheme}, \eqref{ignore-bh-eq1} and \eqref{derivatives-palpha}, we can write
	\begin{align*}
		&\Delta_6 =  \frac{1}{h}\bigg|  \int_{(\tau_t^{{h}}-h)\vee t_1}^{t}\int_{\tau^h_s}^{\tau^h_s+h}\int \int \Gamma^h(0,x,\tau^h_s,w) p_\alpha \left(s-\tau^h_s,z-w-b_h(r,w)(s-\tau^h_s)\right)\\
		&\qquad\qquad\qquad\qquad\qquad\qquad\qquad\qquad\qquad\qquad \times (b(s,z)-b_h(r,w))\cdot\nabla p_\alpha (t-s,y-z) \d z \d w \d r \d s \bigg| \\
		&\lesssim   \int_{(\tau_t^{{h}}-h)\vee t_1}^{t}\int \int p_\alpha (\tau_s^h,w-x)p_\alpha (s-\tau^h_s,z-w) |b(s,z)| \frac{p_\alpha (t-s,y-z)}{(t-s)^{\frac{1}{\alpha}}} \d z \d w \d s\\
		&\qquad +  \frac{1}{h}   \int_{(\tau_t^{{h}}-h)\vee t_1}^{t}  \int_{\tau^h_s}^{\tau^h_s+h} \int \int p_\alpha (\tau_s^h,w-x)p_\alpha (s-\tau^h_s,z-w) |b(r,w)| \frac{p_\alpha (t-s,y-z)}{(t-s)^{\frac{1}{\alpha}}} \d z \d w \d r \d s \\
		&=:\Delta_6^1 + \Delta_6^2.
	\end{align*}
	For $\Delta_6^1$, we first use the convolution properties of the stable kernel in $w$ and then apply \eqref{lemma-convo-bulk-2} with $u=(\tau^h_t-h)\vee t_1$, $v=t$, $\beta_1=\frac 1\alpha$ and $\beta_2=0$ to obtain $$\Delta_6^1\lesssim\int_{(\tau_t^{{h}}-h)\vee t_1}^{t} \int p_\alpha (s,z-x) |b(s,z)| \frac{p_\alpha (t-s,y-z)}{(t-s)^{\frac{1}{\alpha}}} \d z  \d s\lesssim p_\alpha (t,y-x)h^{\frac{\gamma}{\alpha}}.$$ 
	
	For $\Delta_6^2$, we use the convolution properties of the stable kernel in $z$ and \eqref{convo-space-sing}
	\begin{align*}
		\Delta_6^2 &\lesssim \frac{1}{h}   \int_{(\tau_t^{{h}}-h)\vee t_1}^{t}  \int_{\tau^h_s}^{\tau^h_s+h}  \int p_\alpha (\tau_s^h,w-x) |b(r,w)| \frac{p_\alpha (t-\tau_s^h,y-w)}{(t-s)^{\frac{1}{\alpha}}}  \d w \d r \d s \\
		&\lesssim p_\alpha (t,y-x)\frac{1}{h}   \int_{(\tau_t^{{h}}-h)\vee t_1}^{t}   \left[ \frac{1}{(\tau_s^h)^{\frac{d}{\alpha p}}} +\frac{1}{(t-\tau_s^h)^{\frac{d}{\alpha p}}} \right]\frac{1}{(t-s)^{\frac{1}{\alpha}}}\int_{\tau^h_s}^{\tau^h_s+h}   \Vert b(r,\cdot)\Vert_{L^p}  \d r\d s\\  
	&\lesssim p_\alpha (t,y-x)h^{-\frac 1q}   \int_{(\tau_t^{{h}}-h)\vee t_1}^{t}   \left[ \frac{1}{(\tau_s^h)^{\frac{d}{\alpha p}}} +\frac{1}{(t-\tau_s^h)^{\frac{d}{\alpha p}}} \right]\frac{1}{(t-s)^{\frac{1}{\alpha}}}\d s  .\end{align*}
	Remarking that $t-\tau_s^h\geq t-s$, that $\tau_s^h \geq h$ and $(t-(\tau_t^{{h}}-h)\vee t_1)\leq 2h$ we get
	\begin{align*}
		\Delta_6^2 	&\lesssim p_\alpha (t,y-x)h^{-\frac{1}{q}}\int_{(\tau_t^{{h}}-h)\vee t_1}^{t}   \left[h^{-\frac{d}{\alpha p}} +\frac{1}{(t-s)^{\frac{d}{\alpha p}}} \right]\frac{1}{(t-s)^{\frac{1}{\alpha}}}\d s  \\
  		&\lesssim p_\alpha (t,y-x)h^{-\frac{1}{q}}\left(h^{-\frac{d}{\alpha p}}(t-(\tau_t^{{h}}-h)\vee t_1)^{1-\frac 1\alpha}+(t-(\tau_t^{{h}}-h)\vee t_1)^{1-\frac 1\alpha-\frac{d}{\alpha p}}\right)\\
		&\lesssim p_\alpha (t,y-x)h^{\frac{\gamma}{\alpha}}.
	\end{align*}
	This is also a valid bound for $|\bar \Delta_6|$ as we only used $|b_h|\lesssim |b|$.\\
	 
	Now that, plugging the above computations for $\Delta_6 $ and \eqref{maj-D2-mainthm}, \eqref{maj-D3-mainthm}, \eqref{MAJ_DELTA_4}, \eqref{maj-D5-mainthm} in \eqref{splitting-error} and using \eqref{derivatives-palpha} for $\Delta_1$, we obtain  
	\begin{align}\label{pre-gronwall}
		&|\Gamma^h (0,x,t,y)-\Gamma(0,x,t,y)| \lesssim p_\alpha (t,y-x)h^{\frac{\gamma}{\alpha}} \nonumber \\ & \qquad\qquad\qquad\qquad +\int_0^t  \int |\Gamma^h(0,x,s,z)-\Gamma(0,x,s,z)| |b(s,z) | \frac{p_\alpha(t-s,y-z)}{(t-s)^{\frac{1}{\alpha}}} \d z \d s.
	\end{align}
	Setting for all $u\in (0,T]$,
	\begin{equation}\label{def-gronwall-function}
		f(u):=\sup_{x,z\in \R^d} \frac{ |\Gamma^h (0,x,u,z)-\Gamma(0,x,u,z)|}{p_\alpha (u,z-x)},
	\end{equation}
we use \eqref{pre-gronwall} then \eqref{p-q-convo} and H\"older's inequality in time to obtain  :
	\begin{align*}
		f(t) &\lesssim h^{\frac{\gamma}{\alpha}} + \sup_{x,y\in \R^d}\frac{1}{p_\alpha (t,y-x)} \int_0^t  \int f(s)p_\alpha (s,z-x) |b(s,z) | \frac{p_\alpha(t-s,y-z)}{(t-s)^{\frac{1}{\alpha}}} \d z \d s\\
		 &\lesssim h^{\frac{\gamma}{\alpha}} + \sup_{x,y\in \R^d}\frac{1}{p_\alpha (t,y-x)} \int_0^t  \Vert p_\alpha (s,\cdot-x) p_\alpha (t-s,y-\cdot)\Vert_{L^{p '}}   \Vert b(s,\cdot) \Vert_{L^p} \frac{f(s)}{(t-s)^{\frac{1}{\alpha}}} \d s\\
		 &\lesssim h^{\frac{\gamma}{\alpha}} +  \left( \int_0^t  \left(  \frac{f(s)}{(t-s)^{\frac{1}{\alpha}}} \left[\frac{1}{s^{\frac{d}{\alpha p}}} + \frac{1}{(t-s)^{\frac{d}{\alpha p}}}\right] \right)^{q'}\d s \right)^{\frac{1}{q'}}.
	\end{align*}
	Up to a convexity inequality, we thus obtain an estimation which permits to conclude by a suitable Gr\"onwall-Volterra lemma
	\begin{align*}
		f(t)^{q'}\lesssim h^{\frac{\gamma q'}{\alpha}} +  \int_0^t f(s)^{q'}   \frac{1}{(t-s)^{\frac{q'}{\alpha}}} \left[\frac{1}{s^{\frac{dq'}{\alpha p}}} + \frac{1}{(t-s)^{\frac{d q'}{\alpha p}}}\right] \d s .
	\end{align*}
	Since $\frac{q'}{\alpha}+\frac{dq'}{\alpha p}<1$, Lemma 2.2 and Example 2.4 \cite{ZhangJFA10} ensure that
	\begin{align*}
		f(t) \lesssim h^{\frac{\gamma}{\alpha}}.
	\end{align*}
	The same reasoning applies for scheme involving $\bar b_h$, which concludes the proof of Theorem \ref{thm-main}.
	
	\section{Proof of Proposition \ref{prop-main-estimates}: Duhamel representation for the density of the schemes and associated controls}\label{subsec-proof-prop1}
	
	\subsection{Duhamel representation for the density of the scheme}
	 Let us first prove \eqref{duhamel-scheme}. Let $t\in (t_k,T]$, $\phi$ be a ${\cal C}^2$ function with compact support and $v(s,y)={\mathbb 1}_{s<t}p_\alpha (t-s,\cdot)\star \phi(y)+{\mathbb 1}_{s=t}\phi(y)$. According to Lemma \ref{lemFK}, $v$ is ${\cal C}^{1,2}$ on $[0,t]\times \R^d$ and satisfies the Feynman-Kac partial differential equation
         $$\forall (s,y)\in[0,t)\times\R^d,\;\partial_s v(s,y)+\mathcal{L}^\alpha v(s,y)=0 .$$Applying Itô's formula between $t_k$ and $t$ to {$v(s,X_s^h)$ where $(X_s^h)_{s\in[t_k,T]}$ denotes the Euler scheme started from $X^h_{t_k}=x$ and evolving according to \eqref{scheme-interpo}}, we obtain :
		\begin{equation*}
			\phi(X_t^h)=v(t_k,x)+M_{t_k,t}^h+\int_{t_k}^t \nabla v(s,X_s^h)\cdot b_h \left(U_{\lfloor \frac{s}{h} \rfloor},X_{\tau_s^h}^h\right) \d s, 
		\end{equation*}
		where $M_{t_k,t}^h=\int_{t_k}^t \int_{\R^d\backslash\{0\}}\Big(v({s,} X_{s^-}^h+x)-v({s,} X_{s^-}^h)\Big) \tilde N({\rm d}s,{\rm d} x)$, in which $\tilde N $ is the compensated Poisson measure associated with $Z$.
		Taking now the expectation (recalling that $(M_{t_k,s}^h)_{s\in [t_k,t]}$ is a martingale) and using Fubini's theorem, we derive
		\begin{equation*}
			\int \phi(y)\Gamma^h (t_k,x,t,y)\d y=v(t_k,x)+\int_{t_k}^t \E_{t_k,x} \left[ \nabla v(s,X_s^h)\cdot b_h \left(U_{\lfloor \frac{s}{h} \rfloor},X_{\tau_s^h}^h\right)\right]\d s.
		\end{equation*}
		Using the definition of $v$, we get
		\begin{align*}
			&\int \phi(y)\Gamma^h (t_k,x,t,y)\d y\\ & \qquad=\int p_\alpha (t-t_k,x-y)\phi(y) \d y+\int \phi(y)\int_{t_k}^t   \E_{t_k,x} \left[ \nabla_y p_\alpha (t-s,X_s^h-y) \cdot b_h \left(U_{\lfloor \frac{s}{h} \rfloor},X_{\tau_s^h}^h\right)\right] \d s\d y.
		\end{align*}
		Since $\phi$ is arbitrary and $p_\alpha(t-s,\cdot)$ is even, we deduce that $\mathrm{d} y$ a.e., 
		\begin{equation}\label{duhamel-inproof}
			\Gamma^h (t_k,x,t,y)= p_\alpha (t-t_k,x-y)-\int_{t_k}^t   \E_{t_k,x} \left[ \nabla_y p_\alpha (t-s,y-X_s^h) \cdot b_h \left(U_{\lfloor \frac{s}{h} \rfloor},X_{\tau_s^h}^h\right)\right] \d s.
		\end{equation}
		We will see later that \eqref{duhamel-inproof} actually holds for all $y\in \R^d$ as a consequence of the H\"older regularity of $\Gamma^h$ in the forward space variable. This concludes the proof of \eqref{duhamel-scheme}.
		
		\subsection{Heat kernel bounds for the scheme}
		We will now prove inequality \eqref{ineq-density-scheme}, upper stable bound for the density of the scheme, in 3 steps. First, we will prove it for $t\in (t_k,t_{k+1}]$, using only the definition of the cutoffed drift and assuming $h<1$. Then, we will prove it between $t_k$ and $t_\ell$, when $t_k-t_\ell$ is small enough at a \textit{macro} scale. We will finally chain the previous estimates to obtain \eqref{ineq-density-scheme} for any time interval $(t_k,t]\subset [0,T]$.\\
		
		\begin{paragraph}{Step 1 : $t\in (t_k,t_{k+1}]$\\[0.5cm]}
			
			Remarking that when $t\in (t_k,t_{k+1}]$, $\forall z\in \R^d,\ \Gamma^h(t_k,x,t,z)=\frac{1}{h}\int_{t_k}^{t_{k+1}} p_\alpha(t-t_k,z-x-(t-t_k)b_h(r,x)) \d r$ (resp. $\bar{\Gamma}^h(t_k,x,t,z)=\frac{1}{h}\int_{t_k}^{t_{k+1}}p_\alpha(t-t_k,z-x-(t-t_k)\bar{b}_h(r,x)) \d r$), we obtain \eqref{ineq-density-scheme} in the case $t\in (t_k,t_{k+1}]$  using \eqref{ignore-bh-eq1} from Lemma \ref{lemma-ignore-bh} to get rid of the drift.
		\end{paragraph}
		
		\begin{paragraph}{Step 2 : $t-t_k$ \textit{small enough}\\[0.5cm]}
			
			Recall that for $j\in \{k,\cdots, \lceil t/h\rceil-1\}$ and $r\in [t_j,t_{j+1}]$, $X_r^h = X_{t_j}^h+(Z_r-Z_{t_j})+b_h(U_j,X_{t_j}^h)(r-t_j)$. Using this and the independence between $X_{t_j}^h$,  $(Z_r-Z_{t_j})$ and $U_j$, we have, starting from the Duhamel representation \eqref{duhamel-scheme},
\begin{align*}
				\Gamma^h& (t_k,x,t,y)=p_\alpha (t-t_k,y-x)-\sum_{j=k}^{\lceil \frac{t}{h}\rceil-1}\int_{t_j}^{t_{j+1}\wedge t} \E_{t_k,x} \left[ \nabla_y p_\alpha (t-r,y-X_r^h) \cdot b_h \left(U_{j},X_{t_j}^h\right)\right] \d r\\
				&=p_\alpha (t-t_k,y-x)\\& \hspace*{7pt}-\!\!\sum_{j=k}^{\lceil \frac{t}{h}\rceil-1}\int_{t_j}^{t_{j+1}\wedge t}\!\!\frac{1}{h} \int_{t_j}^{t_{j+1}} \!\!\!\!\E_{t_k,x} \left[ \nabla_y p_\alpha \left(t-r,y-X_{t_j}^h-(Z_r-Z_{t_j})-b_h(s,X_{t_j}^h)(r-t_j) \right) \cdot b_h \left(s,X_{t_j}^h\right)\right] \d s\d r.
			\end{align*}
			Using Fubini's and Lebesgue's theorems and the convolution property of the stable density, 
			\begin{align*}
				&\E_{t_k,x} \left[ \nabla_y p_\alpha \left(t-r,y-X_{t_j}^h-(Z_r-Z_{t_j})-b_h(s,X_{t_j}^h)(r-t_j)\right)\cdot b_h \left(s,X_{t_j}^h\right)\right] \\
				&\qquad\qquad\qquad= \int \E_{t_k,x} \left[ \nabla_y p_\alpha \left(t-r,y-X_{t_j}^h-z-b_h(s,X_{t_j}^h)(r-t_j)\right) \cdot b_h \left(s,X_{t_j}^h\right)\right]p_\alpha (r-t_j,z) \d z\\
				&\qquad\qquad\qquad= \E_{t_k,x} \left[ \nabla_y  \left( \int p_\alpha \left(t-r,y-X_{t_j}^h-z-b_h(s,X_{t_j}^h)(r-t_j)\right) p_\alpha (r-t_j,z) \d z \right) \cdot b_h \left(s,X_{t_j}^h\right)\right]	\\
				&\qquad\qquad\qquad =\E_{t_k,x} \left[ \nabla_y   p_\alpha \left(t-t_j,y-X_{t_j}^h-b_h(s,X_{t_j}^h)(r-t_j)\right) \cdot b_h \left(s,X_{t_j}^h\right) \right].
			\end{align*}
			Hence
			\begin{align}\label{duhamel-decoupe}
				\Gamma^h \nonumber& (t_k,x,t,y)
				=p_\alpha (t-t_k,y-x)-\int_{t_k}^{t_{k+1}\wedge t}\frac{1}{h} \int_{t_k}^{t_{k+1}} \nabla_y p_\alpha \left(t-t_k,y-x-b_h(s,x)(r-t_k)\right)\cdot b_h \left(s,x\right) \d z \d s\d r\nonumber\\ &\qquad- \sum_{j=k+1}^{\lceil \frac{t}{h}\rceil-1}\int_{t_j}^{t_{j+1}\wedge t}\frac{1}{h} \int_{t_j}^{t_{j+1}} \int \Gamma^h (t_k,x,t_j,z)  \nabla_y p_\alpha \left(t-t_j,y-z-b_h(s,z)(r-t_j)\right)\cdot b_h \left(s,z\right) \d z\d s\d r.
			\end{align}
			Note that we have not used any property related to $b_h$ here, so the same holds with $(\bar \Gamma^h ,\bar b_h)$ in place of $(\Gamma^h,b_h)$.\\
			
			Set for $j\in \{k+1,\cdots,n\} $, $m_{k,j}:=\sup_{x,y\in \R^d}  \frac{\Gamma^h (t_k,x,t_j,y)}{p_\alpha (t_j-t_k,y-x)}$. Observe from the previous one-step part that there exists $C\ge 1$ s.t. $m_{k,j}\le C^{n-k}<+\infty $. The point of step 2 is to make this bound uniform in $n$. Using \eqref{ignore-bh-eq1} to get rid of the \textit{negligible} cutoffed drift, we get, for $n\ge\ell\geq k+1\ge 1$:
			\begin{align*} 
				&\frac{\Gamma^h (t_k,x,t_\ell,y)}{p_\alpha (t_\ell-t_k,y-x)} \lesssim 1+\int_{t_k}^{t_{k+1}}\frac{1}{h(t_\ell-t_k)^{\frac{1}{\alpha}}} \int_{t_k}^{t_{k+1}}  |b_h (s,x)| \d s\d r\\ &\qquad+ \sum_{j=k+1}^{\ell-1}\int_{t_j}^{t_{j+1}}\frac{1}{h} \int_{t_j}^{t_{j+1}} \int m_{k,j}\frac{p_\alpha (t_j-t_k,z-x)}{p_\alpha (t_\ell-t_k,y-x)}\times \frac{1}{(t_\ell-t_j)^{\frac{1}{\alpha}}} p_\alpha \left(t_\ell-t_j,y-z\right) |b_h (s,z)| \d z\d s\d r.
			\end{align*}
			In the first integral, we use the bound $|b_h|\leq h^{-\frac{d}{\alpha p}-\frac{1}{q}}$ (the bound remains valid for $\bar b_h$ since the latter vanishes on the first time step) and in the second we use $t_\ell-t_j\ge t_\ell -s$ for $s\in[t_j,t_{j+1}]$ and then bound $m_{k,j}$ from above:
						\begin{align*} 
				&\frac{\Gamma^h (t_k,x,t_\ell,y)}{p_\alpha (t_\ell-t_k,y-x)}\\
				&\lesssim 1+\frac{h^{1-\frac{d}{\alpha p}-\frac{1}{q}}}{(t_\ell-t_k)^{\frac{1}{\alpha}}}  + \max_{j\in \llbracket k+1,\ell-1\rrbracket} m_{k,j} \int_{t_{k+1}}^{t_{\ell}} \int \frac{1}{(t_\ell-s)^{\frac{1}{\alpha}}} \frac{p_\alpha (s-t_k,z-x)}{p_\alpha (t_\ell-t_k,y-x)}p_\alpha \left(t_\ell-s,y-z\right) |b (s,z)| \d z \d s.
			\end{align*}
			We are now in the right setting to apply \eqref{lemma-convo-bulk-2} (with $u=t_{k+1}, v=t_\ell,\beta_1=1/\alpha, \beta_2=0$), which readily gives
			\begin{align*} 
				\frac{\Gamma^h (t_k,x,t_\ell,y)}{p_\alpha (t_\ell-t_k,y-x)}& \lesssim 1+\frac{h^{1-\frac{d}{\alpha p}-\frac{1}{q}}}{(t_\ell-t_k)^{\frac{1}{\alpha}}}  + \max_{j\in \llbracket k+1,\ell-1\rrbracket} m_{k,j} (t_{\ell}-{t_{k+1}})^{\frac{\gamma}{\alpha}}\\
				&\lesssim 1+(t_\ell-t_k)^{\frac{\gamma}{\alpha}}  + \max_{j\in \llbracket k+1,\ell-1\rrbracket} m_{k,j} (t_{\ell}-{t_{k+1}})^{\frac{\gamma}{\alpha}}\\
				&\lesssim 1+ \max_{j\in \llbracket k+1,\ell-1\rrbracket} m_{k,j} (t_{\ell}-{t_{k+1}})^{\frac{\gamma}{\alpha}}.
			\end{align*}
			Taking the supremum over $(x,y)\in \R^d$ in the l.h.s., and remarking that the r.h.s. is non-decreasing with $\ell$, along with the definition of  $\lesssim$, we get
			$$\max_{j\in \llbracket k+1,\ell \rrbracket} m_{k,j} \leq C+C (t_\ell-t_k)^{\frac{\gamma}{\alpha}}\max_{j\in \llbracket k+1,\ell \rrbracket} m_{k,j}$$
			for some constant $C$ not depending on $h$. Thus, if $C (t_\ell-t_k)^{\frac{\gamma}{\alpha}}< 1$, then 
			\begin{equation}\label{macro-size-bound}
				\max_{j\in \llbracket k+1,l \rrbracket} m_{k,j}\leq \frac{C}{1-C(t_\ell-t_k)^{\frac \gamma\alpha}}.
			\end{equation}
			In particular, it is bounded uniformly in $h$ for $k,\ell$ s.t. $(t_\ell-t_k)<C^{-\frac{\alpha}{\gamma}}$.\\
			As we only used the fact that $|b_h|\lesssim b$ for the main term, which remains true with $\bar b_h$ instead of $b_h$, the same estimates hold for $\bar m_{k,j}:=\sup_{x,y\in \R^d}  \frac{\bar \Gamma^h (t_k,x,t_j,y)}{p_\alpha (t_j-t_k,y-x)}$.
			\begin{remark} Note that in the Gaussian setting, a precise control of the variance was required because of the exponential structure of the Gaussian tails (see \cite{JM23}). In the stable setting, as the tails of the stable kernel are polynomial, these controls are not required.
			\end{remark}
		\end{paragraph}
		
		\begin{paragraph}{Step 3 : chaining the previous estimates \\[0.5cm]}
		In order to obtain the result for any arbitrary time interval, we will now chain  the previous estimates. This will be done in the following way: denote $\theta =C^{-\frac{\alpha}{\gamma}}$ and let us first suppose that $h\le\theta$, which implies that $\tau^h_\theta\ge \frac{\theta}{2}$. Let $t>0$ s.t. $t-t_k>\theta$ and let $J=\lceil \frac{t-t_k}{\tau_\theta^h} \rceil-1\le\frac{2T}{\theta}$. We will first divide $(t_k,t)$ into a main term (over $(t_k,t_k+J\tau_\theta^h)$) composed of $J$ slices of size $\tau_\theta^h$ (and thus on which we can use \eqref{macro-size-bound}) and a remainder term (over $(t_k+J\tau_\theta^h,t)$). This remainder term will then be split into two terms again ($(t_k+J\tau_\theta^h,\tau_t^h)$ and $(\tau_t^h,t)$), in order to account for the fact that $t$ does not necessarily belong to the discretization grid. Over $(t_k+J\tau_\theta^h,\tau_t^h)$, as we work on the grid, we will use \eqref{macro-size-bound} again, and over $(\tau_t^h,t)$, we will use the cutoff and \eqref{one-step-ineq}.\\
			
			With the convention $y_0=x$,
			\begin{align*}
				\Gamma^h (t_k,x,t,y)&=\int_{(\R^d)^J}\prod_{j=1}^{J}\Gamma^h (t_k+(j-1)\tau_\theta^h,y_{j-1},t_k+j\tau_\theta^h,y_j)\Gamma^h(t_k+J\tau_\theta^h,y_J,t,y)\d y_1 ...\d y_J\\
				&\leq \int_{(\R^d)^J}\prod_{j=1}^{J}m_{k+(j-1)\lfloor \frac{\theta}{h}\rfloor,k+j\lfloor \frac{\theta}{h}\rfloor}p_\alpha (\tau_\theta^h,y_j-y_{j-1})\Gamma^h(t_k+J\tau_\theta^h,y_J,t,y)\d y_1 ...\d y_J.
			\end{align*}
			Using the boundedness of $m_{k+(j-1)\lfloor \frac{\theta}{h}\rfloor,k+j\lfloor \frac{\theta}{h}\rfloor}$, we get
			\begin{align*}
				\Gamma^h (t_k,x,t,y)&\lesssim \int_{(\R^d)^J}\prod_{j=1}^{J}p_\alpha (\tau_\theta^h,y_j-y_{j-1})\Gamma^h(t_k+J\tau_\theta^h,y_J,t,y)\d y_1 ...\d y_J\\
				&\lesssim \int p_\alpha (J\tau_\theta^h,y_J-y_0)\Gamma^h(t_k+J\tau_\theta^h,y_J,t,y)\d y_J.
			\end{align*}
			Pay attention that the constants grow exponentially fast with $J$, but $J\le\frac{2T}{\theta}$.
			Remarking that 
			\begin{align*}
				\Gamma^h(t_k+J\tau_\theta^h,y_J,t,y)&= \int \Gamma^h (t_k+J\tau_\theta^h,y_J,\tau_t^h,z)\Gamma^h(\tau_t^h,z,t,y)\d z\\
				&\lesssim m_{k+J \lfloor \frac{\theta}{h}\rfloor,\tau_t^h} \int  p_\alpha (t_k+J\tau_\theta^h,y_J,\tau_t^h,z)p_\alpha(\tau_t^h,z,t,y)\d z\\
				&\lesssim p_\alpha(t_k+J\tau_\theta^h,y_J,t,y),
			\end{align*}
			we obtain, by convolution, 
			\begin{equation*}
				\Gamma^h (t_k,x,t,y)\lesssim p_\alpha (t_k,x,t,y).
                              \end{equation*}
            When $h>\theta$, then $\frac{T}{h}<\frac{T}{\theta}$ and the conclusion remains valid by chaining in a similar way with $\tau^h_\theta$ replaced by $h$ the estimate derived in Step 1.
			The same reasoning applied to $\bar \Gamma^h$ gives	
			\begin{equation*}
				\bar \Gamma^h (t_k,x,t,y)\lesssim p_\alpha (t_k,x,t,y),
			\end{equation*}
			which concludes the proof of \eqref{ineq-density-scheme}.\\
		\end{paragraph}
		
		\subsection{H\"older regularity of $\Gamma^h$ in the forward  variables}
We will establish here the H\"older properties for the density of the scheme stated in  Proposition \ref{prop-main-estimates}. We begin with the forward time
variable and discuss the forward space variable later on. 		
		\subsubsection{H\"older regularity of $\Gamma^h$ in the forward time variable}
			\label{SUB_H_REG_IN_F_TIME}
			
			Let us now prove \eqref{holder-time-gammah}. Let $0\leq k<\ell<n$, $x,y\in \R^d$ and $t\in [t_\ell,t_{\ell+1}]$. \\

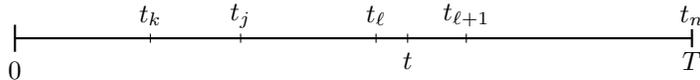
			\begin{centering}
				\begin{tikzpicture}[scale=6]
					\draw[thick] (1.5,0.02)--(1.5,-0.02) node[anchor=north] {};
					\draw[thick] (0,0.03)--(0,-0.03) node[anchor=north] {0};
					
					
					\node at (0.3,0.05) {$t_k$};
					\draw (0.3,0.01)--(0.3,-0.01) {};
					
					\node at (0.5,0.05) {$t_j$};
					\draw (0.5,0.01)--(0.5,-0.01) {};
					
					\node at (0.8,0.05) {$t_\ell$};
					\draw (0.8,0.01)--(0.8,-0.01) {};
					
					\node at (0.87,-0.05) {$t$};
					\draw (0.87,0.01)--(0.87,-0.01) {};
					
					\node at (1,0.05) {$t_{\ell+1}$};
					\draw (1,0.01)--(1,-0.01) {};
					
					\node at (1.5,0.05) {$t_{n}$};
					\draw (1.5,0.01)--(1.5,-0.01) {};
					
					\node at (1.5,-0.05) {$T$};
					
					\draw[thick] (0,0)--(1.5,0);
					\end{tikzpicture}\\[0.5cm]
				\end{centering}
				
				Going back to \eqref{duhamel-decoupe}, we can write:
				\begin{align}
					\Gamma^h& (t_k,x,t_\ell,y)-\Gamma^h (t_k,x,t,y) =p_\alpha (t_\ell-t_k,y-x) - p_\alpha (t-t_k,y-x)\nonumber\\
					&\nonumber-\frac{1}{h}\int_{t_k}^{t_{k+1}} \int_{t_k}^{t_{k+1}} \left[\nabla p_\alpha \left(t_\ell-t_k,w\right)-\nabla p_\alpha \left(t-t_k,w\right) \right]_{w=y-x-b_h(s,x)(r-t_k)}  \cdot b_h \left(s,x\right) \d s\d r\\
					&\nonumber-\frac{\mathbb{1}_{\{\ell \geq k+2\}}}{h} \sum_{j=k+1}^{\ell-1}\int_{t_j}^{t_{j+1}}\int_{t_j}^{t_{j+1}} \int \Gamma^h (t_k,x,t_j,z)  [\nabla p_\alpha \left(t_\ell-t_j,w\right) \\& \nonumber\qquad\qquad\qquad\qquad\qquad\qquad\qquad\qquad\qquad -\nabla p_\alpha \left(t-t_j,w\right) ]_{w=y-z-b_h(s,z)(r-t_j)} \cdot b_h \left(s,z\right) \d z \d s\d r\\
					&\nonumber+\frac{1}{h}\int_{t_\ell}^t \int_{t_\ell}^{t_{\ell+1}} \int \Gamma^h (t_k,x,t_\ell,z)b_h(s,z)\cdot \nabla p_\alpha (t-t_\ell,y-z-b_h(s,z)(r-t_\ell))\d z \d s \d r\\
					&=: \Delta_1 +  \Delta_2 +\Delta_3 +  \Delta_4. \label{splitting-holder-time-gammah}
				\end{align}
				Resp.  $\bar\Gamma^h (t_k,x,t_\ell,y)-\bar\Gamma^h (t_k,x,t,y) =: \bar \Delta_1 +  \bar\Delta_2 +\bar\Delta_3 +\bar  \Delta_4$ for the scheme involving $\bar b_h$.\\
				
				For $\Delta_1$ (which is actually the same as $\bar \Delta_1$), we use \eqref{holder-time-palpha} and $t-t_k\asymp t_\ell-t_k$ then $t-t_\ell<t-t_k$:
				\begin{align}
					|\Delta_1|\lesssim \frac{t-t_\ell}{t-t_k} p_\alpha (t-t_k,y-x) \lesssim \left( \frac{t-t_\ell}{t-t_k}\right)^{\frac{\gamma}{\alpha}} p_\alpha (t-t_k,y-x).\label{maj-holder-time-d1}
				\end{align}
				For $\Delta_2$, let us first bound $\left[\nabla p_\alpha \left(t_\ell-t_k,w\right)-\nabla p_\alpha \left(t-t_k,w\right) \right]_{w=y-x-b_h(s,x)(r-t_k)}$, using again \eqref{holder-time-palpha} along with $t-t_k\asymp t_\ell-t_k$:
				\begin{align*}
					\left|\nabla p_\alpha \left(t_\ell-t_k,w\right)-\nabla p_\alpha \left(t-t_k,w\right) \right|_{w=y-x-b_h(s,x)(r-t_k)} &\lesssim \frac{|t-t_\ell|}{(t-t_k)^{1+\frac{1}{\alpha}}} {p}_\alpha (t-t_k,y-x-b_h(s,x)(r-t_k))
				\end{align*}
				In our current integral, $r-t_k\leq t-t_k$, which means that, using \eqref{ignore-bh-eq1}, we get ${p}_\alpha (t-t_k,y-x-b_h(s,x)(r-t_k)) \lesssim {p}_\alpha (t-t_k,x-y)$. We can thus compute the following bound for $\Delta_2$ (recalling that $|b_h| \lesssim h^{-\frac{d}{\alpha p}-\frac{1}{q}} $):
				\begin{align*}
					|\Delta_2| &= \frac{1}{h}\left|\int_{t_k}^{t_{k+1}} \int_{t_k}^{t_{k+1}} \left[\nabla p_\alpha \left(t_\ell-t_k,w\right)-\nabla p_\alpha \left(t-t_k,w\right) \right]_{w=y-x-b_h(s,x)(r-t_k)}\cdot b_h \left(s,x\right) \d s\d r \right| \\
					&\lesssim\frac{1}{h}\int_{t_k}^{t_{k+1}} \int_{t_k}^{t_{k+1}}\frac{t-t_\ell}{(t-t_k)^{1+\frac{1}{\alpha}}} {p}_\alpha (t-t_k,y-x) |b_h (s,x) |\d s\d r\\
					&\lesssim\frac{t-t_\ell}{(t-t_k)^{1+\frac{1}{\alpha}}} {p}_\alpha (t-t_k,y-x) \int_{t_k}^{t_{k+1}} |b_h (s,x) |\d s\\
					&\lesssim \frac{t-t_\ell}{(t-t_k)^{1+\frac{1}{\alpha}}} {p}_\alpha (t-t_k,y-x) h^{1-\frac{d}{\alpha p}-\frac{1}{q}}.
				\end{align*}
				Using the fact that $t-t_k\geq t_\ell-t_k\geq h$, we get 
				\begin{align}
					|\Delta_2|&\lesssim \frac{t-t_\ell}{t-t_k} {p}_\alpha (t-t_k,y-x) h^{1-\frac{1}{\alpha}-\frac{d}{\alpha p}-\frac{1}{q}}\lesssim\left( \frac{t-t_\ell}{t-t_k} \right)^{\frac{\gamma}{\alpha}}{p}_\alpha (t-t_k,y-x) h^{\frac{\gamma}{\alpha}}.\label{maj-holder-time-d2}
				\end{align}
				For the alternative scheme, we would have used the inequality $|\bar b_h|\lesssim h^{\frac{1}{\alpha}-1}$, which yields 
				$$|\bar \Delta_2|\lesssim \frac{t-t_\ell}{(t-t_k)^{1+\frac{1}{\alpha}}} {p}_\alpha (t-t_k,y-x)\lesssim\left( \frac{t-t_\ell}{t-t_k} \right)^{\frac{\gamma}{\alpha}}{p}_\alpha (t-t_k,y-x).$$
				
				For $\Delta_3$, note that for all $j\in\llbracket k+1,\ell-1 \rrbracket$, denoting $u=t_\ell-t_j$ and $u'=t-t_j$, we can see $u'-u=t-t_\ell\in [0,h]$ as a small perturbation at the scale of $u$ or $u'$.
				This allows us to use \eqref{holder-time-palpha} along with $t-t_j\asymp t_\ell-t_j$ and then \eqref{ignore-bh-eq1}:
				\begin{align*}
					|\nabla p_\alpha \left(t_\ell-t_j,w\right)-\nabla p_\alpha \left(t-t_j,w\right)|_{w=y-z-b_h(s,z)(r-t_j)} &\lesssim \frac{t-t_\ell}{(t-t_j)^{1+\frac{1}{\alpha}}}{p}_\alpha(t-t_j,y-z-b_h(s,z)(r-t_j))\\
					&\lesssim \frac{t-t_\ell}{(t-t_j)^{1+\frac{1}{\alpha}}}{p}_\alpha(t-t_j,y-z).
				\end{align*}
				For the computations on $\Delta_3$, we assume that $\ell\ge k+2$ and introduce an exponent $\gamma_1 \in (\gamma,\alpha]$. Here, we singularize some of the estimates in order to obtain the expected H\"older rate involving $\gamma$. This is somehow a flexibility of the scheme: since we stay away from the final time $t$ for this contribution, we can afford to make non-integrable exponents appear. Those terms will be handled with Lemma \ref{lemma-convo-bulk} (eq. \eqref{lemma-convo-bulk-1}).
                Namely, using $|b_h|\leq b$ and  the stable upper-bound \eqref{ineq-density-scheme}, we get
				\color{black}
				\begin{align*}
					|\Delta_3|&\lesssim \frac{1}{h} \sum_{j=k+1}^{\ell-1} \int_{t_j}^{t_{j+1}}\int_{t_j}^{t_{j+1}} \int  \Gamma^h (t_k,x,t_j,z)  \frac{t-t_\ell}{(t-t_j)^{1+\frac{1}{\alpha}}}{p}_\alpha(t-t_j,y-z)|b (s,z)| \d z \d s\d r\\
					&\lesssim \sum_{j=k+1}^{\ell-1} \int_{t_j}^{t_{j+1}} \frac{(t-t_\ell)^{\frac{\gamma_1}{\alpha}}}{(t-t_j)^{\frac{\gamma_1+1}{\alpha}}}\int {p}_\alpha (t_j-t_k,z-x) {p}_\alpha(t-t_j,y-z)|b (s,z)| \d z \d s\\
					&\lesssim \sum_{j=k+1}^{\ell-2} \int_{t_j}^{t_{j+1}} \frac{(t-t_\ell)^{\frac{\gamma_1}{\alpha}}}{(t-t_j)^{\frac{\gamma_1+1}{\alpha}}}\int {p}_\alpha (t_j-t_k,z-x) {p}_\alpha(t-t_j,y-z)|b (s,z)| \d z \d s\\
					&\qquad + \int_{t_{\ell-1}}^{t_\ell} \frac{(t-t_\ell)^{\frac{\gamma_1}{\alpha}}}{(t-t_{\ell-1})^{\frac{\gamma_1+1}{\alpha}}}\int {p}_\alpha (t_{\ell-1}-t_k,z-x) {p}_\alpha(t-t_{\ell-1},y-z)|b (s,z)| \d z \d s\\
					&=:\Delta_3^1 + \Delta_3^2.
				\end{align*}				
				Assume $\ell \geq k+3$ (otherwise $\Delta_3^1$ vanishes). In $\Delta_3^1$, which only contains non-singular integrals, we now approximate the discrete $(t_j)_{j\in \llbracket k+1,l-2 \rrbracket}$ with $s$ in the corresponding time integrals to apply \eqref{lemma-convo-bulk-1} with $u=t_{k+1},v=t_{\ell-1}, \beta_1=(\gamma_1+1)/\alpha, \beta_2=0 $:
				\begin{align*}
					\Delta_3^1 &\lesssim \mathbb{1}_{\{\ell \geq k+3\}}\int_{t_{k+1}}^{t_{\ell-1}} \frac{(t-t_\ell)^{\frac{\gamma_1}{\alpha}}}{(t-s)^{\frac{\gamma_1+1}{\alpha}}}\int {p}_\alpha (s-t_k,z-x) {p}_\alpha(t-s,z-y)|b (s,z)| \d z \d s\\
					&\lesssim p_\alpha (t-t_k,y-x) (t-t_\ell)^{\frac{\gamma_1}{\alpha}}\left[(t_{\ell-1}-t_{k+1})^{\frac{\gamma-\gamma_1}{\alpha}}+(t-t_{\ell-1})^{\frac{\gamma-\gamma_1}{\alpha}}\right]\\
                    &\lesssim p_\alpha (t-t_k,y-x) (t-t_\ell)^{\frac{\gamma_1}{\alpha}}h^{\frac{\gamma-\gamma_1}{\alpha}},
				\end{align*}
                where, for the last inequality, we used $t_{\ell-1}-t_{k+1}\ge h$ since $\ell-1\ge (k+1)+1$, $t-t_{\ell-1}=t-t_{\ell}+h\geq h$ and $\gamma-\gamma_1<0$.\\
                For $\Delta_3^2$ (last time step), let us first use the convolution estimate \eqref{convo-space-sing}:
                \begin{align*}
                    \Delta_3^2 \lesssim p_\alpha (t-t_k,y-x)\frac{(t-t_\ell)^{\frac{\gamma_1}{\alpha}}}{(t-t_{\ell-1})^{\frac{\gamma_1+1}{\alpha}}} \left[\frac{1}{(t_{\ell-1}-t_k)^{\frac{d}{\alpha p}}} + \frac{1}{(t-t_{\ell-1})^{\frac{d}{\alpha p}}}\right] \int_{t_{\ell-1}}^{t_\ell}  \Vert b (s,\cdot) \Vert_{L^p} \d s.
                \end{align*}
				Using $ \left[\frac{1}{(t_{\ell-1}-t_k)^{\frac{d}{\alpha p}}} + \frac{1}{(t-t_{\ell-1})^{\frac{d}{\alpha p}}}\right]\frac{1}{(t-t_{\ell-1})^{\frac{\gamma_1+1}{\alpha}}}\le h^{-\frac{d}{\alpha p}-\frac{\gamma_1+1}{\alpha}}$ and applying H\"older's inequality to the integral, we obtain:
                \begin{align*}
                    \Delta_3^2 &\lesssim p_\alpha (t-t_k,y-x)(t-t_\ell)^{\frac{\gamma_1}{\alpha}} h^{-\frac{d}{\alpha p}-\frac{\gamma_1+1}{\alpha}} \Vert \mathbb{1}_{(t_{\ell-1},t_\ell)}\Vert_{L^{q'}}\nonumber\\
                    & \lesssim p_\alpha (t-t_k,y-x)(t-t_\ell)^{\frac{\gamma_1}{\alpha}} h^{\frac{\gamma-\gamma_1}{\alpha}}.
                \end{align*}
				Gathering both estimates and recalling that $ t-t_\ell\le h$ and $\gamma-\gamma_1<0 $, we obtain 
                \begin{align}
                    |\Delta_3|
					&\lesssim {p}_\alpha (t-t_k,y-x)(t-t_\ell)^{\frac{\gamma}{\alpha}}.\label{maj-holder-time-d3}
                \end{align}
				Let us now bound $\Delta_4$ in \eqref{splitting-holder-time-gammah}. Recalling that from its definition, $|b_h|\leq |b|$ (resp. $|\bar b_h|\leq |b|$), we can write, using also \eqref{ineq-density-scheme} and \eqref{ignore-bh-eq1}:
				\begin{align*}
					|\Delta_4| 	&\lesssim  \frac{(t-t_\ell)^{1-\frac{1}{\alpha}}}{h} \int_{t_\ell}^{t_{\ell+1}} \int p_\alpha (t_\ell-t_k,z-x)|b(s,z)| p_\alpha (t-t_\ell,z-y)\d z \d s.
				\end{align*}	
				We can now bound $|\Delta_4|$ using \eqref{convo-space-sing}, then $(t_\ell-t_k)^{-1} \leq (t-t_\ell)^{-1}$ and finally $t-t_\ell \leq h$:
				\begin{align}
					|\Delta_4| &\nonumber\lesssim  \frac{(t-t_\ell)^{1-\frac{1}{\alpha}}}{h} \left[ \frac{1}{(t_\ell - t_k)^{\frac{d}{\alpha p}}} + \frac{1}{(t-t_\ell)^{\frac{d}{\alpha p}}}\right] {p}_\alpha (t-t_k,y-x) \int_{t_\ell}^{t_{\ell+1}} \Vert b(s,\cdot) \Vert_{L^{p}} \d s\\
					&\nonumber\lesssim  \frac{(t-t_\ell)^{1-\frac{1}{\alpha}-\frac{d}{\alpha p}}}{h} {p}_\alpha (t-t_k,y-x) \Vert b \Vert_{L^q  -L^p} \Vert \mathbb{1}_{(t_\ell,t_{\ell+1})} \Vert_{L^{q'} }\\
					&\nonumber\lesssim  (t-t_\ell)^{1-\frac{1}{\alpha}-\frac{d}{\alpha p}} h^{1-\frac{1}{q}-1} {p}_\alpha (t-t_k,y-x)\\
					&\nonumber\lesssim (t-t_\ell)^{1-\frac{1}{\alpha}-\frac{d}{\alpha p}-\frac{1}{q}} {p}_\alpha (t-t_k,y-x)\\
					&\lesssim (t-t_\ell)^{\frac{\gamma}{\alpha}} {p}_\alpha (t-t_k,y-x).\label{maj-holder-time-d4}
				\end{align}
				As, for $\Delta_3$ and $\Delta_4$, we only used the fact that $|b_h|\leq |b|$, the same estimations still hold for the alternative scheme involving $\bar b_h$.	Plugging the estimates \eqref{maj-holder-time-d1}-\eqref{maj-holder-time-d4}
				into \eqref{splitting-holder-time-gammah} concludes the proof of \eqref{holder-time-gammah}.
			 
			\subsubsection{H\"older regularity of $\Gamma^h$ in the forward space variable}
				Let us now prove \eqref{holder-space-gammah}. This property is important to prove that any limit point of the law induced by the Euler scheme solves the martingale problem and that its marginals will satisfy heat kernel estimates through a compactness type argument (see Section \ref{SEC_MART_PB} for details).  
				\begin{itemize}
					\item \textbf{Off-diagonal regime:} $|y-y'|\gtrsim(t-t_k)^{1/\alpha}$.\\
					In this case, using the stable upper bound \eqref{ineq-density-scheme}, we only need to write
					\begin{align*}
						|\Gamma^h(t_k,x,t,y') -\Gamma^h(t_k,x,t,y)|&\lesssim \Gamma^h(t_k,x,t,y') +\Gamma^h(t_k,x,t,y)\\
						&\lesssim p_\alpha(t-t_k,y'-x)+p_\alpha(t-t_k,y-x)\\
						&\lesssim \frac{|y-y'|^{\gamma} \wedge (t-t_k)^{\frac{\gamma}{\alpha}}}{(t-t_k)^{\frac{\gamma}{\alpha}}}\left( p_\alpha(t-t_k,y'-x)+p_\alpha(t-t_k,y-x)\right).
					\end{align*}
					\item \textbf{Diagonal regime:} $|y-y'|\lesssim (t-t_k)^{1/\alpha}$. Note that in this setting, $p_\alpha (t-t_k,y-x)\asymp p_\alpha (t-t_k,y'-x)$.\\
					In this case, we go back to \eqref{duhamel-decoupe}, denoting $\ell = \lceil t/h \rceil -1$ (so that $t\in (t_\ell,t_{\ell+1}]$) and we write, similarly to \eqref{splitting-holder-time-gammah}:
					\begin{align}
						&\nonumber\Gamma^h(t_k,x,t,y') -\Gamma^h(t_k,x,t,y)=p_\alpha(t-t_k,y'-x)-p_\alpha(t-t_k,y-x)\\
						&\nonumber-\frac{1}{h}\int_{t_k}^{t_{k+1}\wedge t}\int_{t_k}^{t_{k+1}}  [\nabla p_\alpha (t-t_k,y'-w)- \nabla p_\alpha (t-t_k,y-w)]_{w=x+b_h(s,x)(r-t_k)} \cdot b_h (s,x) \d s\d r\\&\nonumber- \frac{1}{h}\sum_{j=k+1}^{\ell-1}\int_{t_j}^{t_{j+1}} \int_{t_j}^{t_{j+1}} \int \Gamma^h (t_k,x,t_j,z)  \\&\nonumber \qquad\qquad\qquad\qquad [\nabla p_\alpha (t-t_j,y'-w)-\nabla p_\alpha (t-t_j,y-w) ]_{w=z+b_h(s,z)(r-t_j)} \cdot b_h (s,z) \d z\d s\d r\\
						&\nonumber-\frac{{\mathbb{1}}_{\ell\ge k+1}}{h}\int_{t_\ell}^{t} \int_{t_\ell}^{t_{\ell+1}} \int \Gamma^h (t_k,x,t_\ell,z)  \\& \nonumber\qquad\qquad\qquad\qquad [\nabla p_\alpha (t-t_\ell,y'-w)-\nabla p_\alpha (t-t_\ell,y-w) ]_{w=z+b_h(s,z)(r-t_\ell)} \cdot b_h (s,z) \d z\d s\d r\\
						&=:\Delta_1+\Delta_2+\Delta_3+\Delta_4. \label{splitting-holder-space-gammah}
					\end{align}
					Resp. $\bar\Gamma^h(t_k,x,t,y') -\bar\Gamma^h(t_k,x,t,y) =: \bar \Delta_1+\bar\Delta_2+\bar\Delta_3+\bar\Delta_4$ for the scheme with $\bar b_h$.\\
					Those terms will be treated in similar way than for the time sensitivities, up to the fact that we will use the H\"older regularity in space of $p_\alpha$ \eqref{holder-space-palpha} instead of its H\"older regularity in time \eqref{holder-time-palpha}.\\
					
					For $\Delta_1$ (which is the same as $\bar{\Delta}_1$), we directly use \eqref{holder-space-palpha} with $\theta=1$  and the diagonal regime:
					\begin{align}\label{maj-holder-space-d1}
						\Delta_1 &\lesssim \frac{|y-y'|}{(t-t_k)^{\frac{1}{\alpha}}}(p_\alpha(t-t_k,y'-x)+p_\alpha(t-t_k,y-x))\nonumber \\
						&\lesssim \frac{|y-y'|^\gamma}{(t-t_k)^{\frac{\gamma}{\alpha}}}(p_\alpha(t-t_k,y'-x)+p_\alpha(t-t_k,y-x)).
					\end{align}
					
					For $\Delta_2$, we first use \eqref{ignore-bh-eq2} to bound $[\nabla p_\alpha (t-t_k,x+b_h(s,x)(r-t_k)-y')- \nabla p_\alpha (t-t_k,x+b_h(s,x)(r-t_k)-y)]$ and get rid of the drift:
					\begin{align*}
						&|\nabla p_\alpha (t-t_k,y'-x-b_h(s,x)(r-t_k))- \nabla p_\alpha (t-t_k,y-x-b_h(s,x)(r-t_k))|\\
						& \qquad \lesssim\frac{|y-y'|}{(t-t_k)^{\frac{2}{\alpha}}}(p_\alpha(t-t_k,y'-x)+p_\alpha(t-t_k,y-x)).
					\end{align*}
					We then compute, recalling from the definition \eqref{cutoff-1} that $|b_h|\lesssim h^{-\frac{d}{\alpha p}-\frac{1}{q}}$ (resp. $|\bar b_h|\lesssim h^{\frac{1}{\alpha}-1}$), 
					\begin{align*}
						\Delta_2 &\lesssim\frac{|y-y'|}{(t-t_k)^{\frac{2}{\alpha}}} (p_\alpha(t-t_k,y'-x)+p_\alpha(t-t_k,y-x)) \int_{t_k}^{t_{k+1}\wedge t}    |b_h (s,x)| \d s\\
						&\lesssim\frac{|y-y'|}{(t-t_k)^{\frac{1}{\alpha}}}(p_\alpha(t-t_k,y'-x)+p_\alpha(t-t_k,y-x)) \frac{(t_{k+1}\wedge t-t_k)h^{-\frac{d}{\alpha p}-\frac{1}{q}}}{(t-t_{k})^{\frac{1}{\alpha}}}\\
&\lesssim\frac{|y-y'|}{(t-t_k)^{\frac{1}{\alpha}}}(p_\alpha(t-t_k,y'-x)+p_\alpha(t-t_k,y-x)) (t_{k+1}\wedge t-t_k)^{1-\frac{1}{\alpha}}h^{-\frac{d}{\alpha p}-\frac{1}{q}}\\						\mathrm{resp.} \qquad \bar \Delta_2 & \lesssim\frac{|y-y'|}{(t-t_k)^{\frac{1}{\alpha}}}(p_\alpha(t-t_k,y'-x)+p_\alpha(t-t_k,y-x)) (t_{k+1}\wedge t-t_k)^{1-\frac{1}{\alpha}}h^{\frac{1}{\alpha}-1}.
					\end{align*}
					In our current diagonal regime, we can write $\frac{|y-y'|}{(t-t_k)^{\frac{1}{\alpha}}} \lesssim \frac{|y-y'|^\gamma}{(t-t_k)^{\frac{\gamma}{\alpha}}}$, which, along with $t_{k+1}\wedge t-t_k\le h$, yields
					\begin{align}\label{maj-holder-space-d2}
						\Delta_2 \lesssim \frac{|y-y'|^\gamma}{(t-t_k)^{\frac{\gamma}{\alpha}}}(p_\alpha(t-t_k,y'-x)+p_\alpha(t-t_k,y-x))h^{\frac{\gamma}{\alpha}}.
					\end{align}
					$$				\mathrm{resp.} \qquad \bar \Delta_2 \lesssim \frac{|y-y'|^\gamma}{(t-t_k)^{\frac{\gamma}{\alpha}}}(p_\alpha(t-t_k,y'-x)+p_\alpha(t-t_k,y-x)).$$
					
					For $\Delta_3$, we first use \eqref{ignore-bh-eq2} and \eqref{ineq-density-scheme} to write, for $\gamma_1 \in (\gamma,1]$
					\begin{align*}
						\Delta_3& \lesssim \sum_{j=k+1}^{\ell-1}\int_{t_j}^{t_{j+1}} \int p_\alpha(t_j-t_k,z-x)\left( \frac{|y-y'|^{\gamma_1}}{(t-t_j)^{\frac{\gamma_1}{\alpha}}}\wedge 1\right)\frac{1}{(t-t_j)^{\frac{1}{\alpha}}}\\ & \qquad\qquad\qquad\qquad\qquad\times \left(p_\alpha(t-t_j,y-z)+p_\alpha(t-t_j,y'-z)\right) |b_h (s,z)| \d z\d s.
					\end{align*}
					Then, we will proceed differently depending on whether, at the current time $t_j$, the spatial difference $|y-y'| $ we are interested in is in the diagonal regime 
                                        w.r.t. the corresponding time scale $t-t_j $. To this end, let us split between what we call \textit{meso-scale} diagonal and off-diagonal regimes (respectively $|y-y'|< (t-t_j)^{\frac{1}{\alpha}}$ and $|y-y'|\geq (t-t_j)^{\frac{1}{\alpha}}$). This meso-scale dichotomy did not appear in the proof of the H\"older time-regularity of $\Gamma^h$. It does now because of technical reasons:  we need to retrieve the loss induced by the introduction of $\gamma_1\in (\gamma,\alpha]$, which can only be done in the mesoscopic diagonal case.\\

					Let us point out that $|y-y'|\leq (t-t_j)^{\frac{1}{\alpha}} \ssi j\leq \frac{t-|y-y'|^\alpha}{h}$. Set $j_{\max}:= \left\lfloor \frac{t-|y-y'|^\alpha}{h} \right\rfloor \wedge (\ell-2)$. We recall that, when $t_j$ is close to $t$, a \textit{local} off-diagonal regime might appear. With the previous notations it will precisely be the case from $t_{j_{\max}+1} $ to $t_\ell $ whenever $j_{\max}<\ell-2$ . We can thus write
					\begin{align*}
						\Delta_3 \lesssim&{\mathbb{1}}_{j_{\max}=\ell-2}\frac{|y-y'|^\gamma}{(t-t_{\ell-1})^{\frac{1+\gamma}{\alpha}}}\int_{t_{\ell-1}}^{t_{\ell}} \int p_\alpha(t_{\ell-1}-t_k,z-x)\\&\phantom{+{\mathbb{1}}_{j_{\max}=\ell-2}\frac{|y-y'|^\gamma}{(t-t_{\ell-1})^{\frac{1+\gamma}{\alpha}}}\int_{t_{\ell-1}}^{t_{\ell}} \int} \left(p_\alpha(t-t_{\ell-1},y-z)+p_\alpha(t-t_{\ell-1},y'-z)\right)|b_h (s,z)| \d z \d s\\
                        &\qquad +{\mathbb{1}}_{j_{\max}<\ell-2}\sum_{j=j_{\max}+1}^{\ell-1} \int_{t_j}^{t_{j+1}} \int p_\alpha(t_j-t_k,z-x) \frac{1}{(t-t_j)^{\frac{1}{\alpha}}}\\ & \qquad\qquad\qquad\qquad\qquad\times\left(p_\alpha(t-t_j,y-z)+p_\alpha(t-t_j,y'-z)\right)|b_h (s,z)| \d z \d s\\
						&\qquad + \sum_{j=k+1}^{j_{\max}}\int_{t_j}^{t_{j+1}}  \int p_\alpha(t_j-t_k,z-x) \frac{|y-y'|^{\gamma_1}}{(t-t_j)^{\frac{1+\gamma_1}{\alpha}}} \\ & \qquad\qquad\qquad\qquad\qquad\times   \left(p_\alpha(t-t_j,y-z)+p_\alpha(t-t_j,y'-z)\right) |b_h (s,z)|  \d z \d s\\
						=:&\Delta_3^{EDGE} +\Delta_3^{OD}+\Delta_3^{D}.
					\end{align*}
                    For the first term, we first use $|b_h|\lesssim |b|$ and \eqref{convo-space-sing}, then H\"older's inequality for the time integral and last $t-t_{\ell-1}=t-t_\ell+h\geq h$ to obtain:
                    \begin{align*}
                        \Delta_3^{EDGE}\lesssim & \frac{p_\alpha(t-t_k,y-x)+p_\alpha(t-t_k,y'-x)}{(t-t_{\ell-1})^{\frac{1}{\alpha}}}\left[\frac{1}{(t_{\ell-1}-t_k)^{\frac{d}{\alpha p}}} + \frac{1}{(t-t_{\ell-1})^{\frac{d}{\alpha p}}}\right] \int_{t_{\ell-1}}^{t_{\ell}}\Vert b (s,\cdot)\Vert_{L^p}  \d s\\
                        &\times \frac{|y-y'|^\gamma}{(t-t_{\ell-1})^{\frac \gamma\alpha}}\\
                        \lesssim &\left(p_\alpha(t-t_k,y-x)+p_\alpha(t-t_k,y'-x)\right) \frac{1}{(t-t_{\ell-1})^{\frac{1}{\alpha}}}\left[\frac{1}{(t_{\ell-1}-t_k)^{\frac{d}{\alpha p}}} + \frac{1}{(t-t_{\ell-1})^{\frac{d}{\alpha p}}}\right]h^{1-\frac{1}{q}}\\
       &                 \times \frac{|y-y'|^\gamma}{(t-t_{\ell-1})^{\frac \gamma\alpha}}\\& \lesssim \left(p_\alpha(t-t_k,y-x)+p_\alpha(t-t_k,y'-x)\right)|y-y'|^\gamma.
                    \end{align*}
                   Next, note that in the integrals appearing in $\Delta_3^{OD}$ and $\Delta_3^D$, we can use $t-s\le t-t_j$ for $s\in[t_j,t_{j+1}]$. Together with $|\bar{b}_h |\leq |b|$, this yields
					\begin{align*}
						\Delta_3^{OD}+\Delta_3^{D}& \lesssim  \mathbb{1}_{j_{\max}<\ell-2}\int_{t_{j_{\max}+1}}^{t_{\ell}} \int p_\alpha(s-t_k,z-x) \frac{1}{(t-s)^{\frac{1}{\alpha}}}\\ & \qquad\qquad\qquad\qquad\qquad \times\left(p_\alpha(t-s,y-z)+p_\alpha(t-s,y'-z)\right)|b (s,z)| \d z \d s\ \\
						&\qquad + \int_{t_{k+1}}^{t_{j_{\max}}}  \int p_\alpha(s-t_k,z-x) \frac{|y-y'|^{\gamma_1}}{(t-s)^{\frac{1+\gamma_1}{\alpha}}} \\ & \qquad\qquad\qquad\qquad\qquad\qquad\times   \left(p_\alpha(t-s,y-z)+p_\alpha(t-s,y'-z)\right) |b (s,z)|  \d z \d s.
					\end{align*}
					For $\Delta_3^{OD}$, we simply use \eqref{lemma-convo-bulk-2}, with $u=t_{j_{\max}+1},v=t_{\ell},\beta_1=1/\alpha,\beta_2=0 $. \\
                    For $\Delta_3^{D}$,
                    \begin{itemize}
                        \item if $\frac 12 (t_{j_{\max}}-t_{k+1})>t-t_{j_{\max}} $, we use \eqref{lemma-convo-bulk-1} with $u=t_{k+1}, v=t_{j_{\max}},\beta_1=(\gamma_1+1)/\alpha,\beta_2=0$.
                        \item if $\frac 12 (t_{j_{\max}}-t_{k+1})\leq t-t_{j_{\max}} $, we use the bound $(t-s)^{-\frac{\gamma_1}\alpha}\le (t-t_{j_{\max}})^{-\frac{\gamma_1}\alpha}$ for $s\le t_{j_{\max}}$ then apply \eqref{lemma-convo-bulk-2} with $u=t_{k+1}, v=t_{j_{\max}},\beta_1=1/\alpha,\beta_2=0 $.
                    \end{itemize}
                    This yields \begin{align*}
						\Delta_3^{OD}+\Delta_3^{D} \lesssim& \left(p_\alpha(t-t_k,y-x)+p_\alpha(t-t_k,y'-x)\right)\left[ (t-t_{j_{\max}+1})^{\frac{\gamma}{\alpha}} \mathbb{1}_{{j_{\max}}<\ell-2}\right.\\
						&+ |y-y'|^{\gamma_1}\big(\{(t_{j_{\max}}-t_{k+1})^{\frac{\gamma-\gamma_1}{\alpha}} +(t-t_{j_{\max}})^{\frac{\gamma-\gamma_1}{\alpha}}\}\mathbb{1}_{\frac 12 (t_{j_{\max}}-t_{k+1})>t-t_{j_{\max}}}\\
						&\left.+(t-t_{j_{\max}})^{-\frac{\gamma_1}\alpha}(t_{j_{\max}}-t_{k+1})^{\frac \gamma\alpha} \mathbb{1}_{\frac 12 (t_{j_{\max}}-t_{k+1})\le t-t_{j_{\max}}} \big) \right]\\
						\lesssim & \left(p_\alpha(t-t_k,y-x)+p_\alpha(t-t_k,y'-x)\right)\left[ (t-t_{j_{\max}+1})^{\frac{\gamma}{\alpha}} \mathbb{1}_{{j_{\max}}<\ell-2}\right.\\
						&\left.+|y-y'|^{\gamma_1} (t-t_{j_{\max}})^{\frac{\gamma-\gamma_1}\alpha}\right].
					\end{align*}
					Since $(t-t_{j_{\max}+1})\leq |y-y'|^\alpha$ if ${j_{\max}}<\ell-2 $ and $(t-t_{j_{\max}})\geq |y-y'|^\alpha$, we obtain
					\begin{align}\label{maj-holder-space-d3}
						\Delta_3^{OD}+\Delta_3^{D}&\lesssim \left(p_\alpha(t-t_k,y-x)+p_\alpha(t-t_k,y'-x)\right) |y-y'|^{\gamma}.
					\end{align}
					Finally, for $\Delta_4$, we suppose that $\ell\ge k+1$ since otherwise this term vanishes.  Using again \eqref{ignore-bh-eq2}, we get
					\begin{align*}
						&|\nabla_{y'} p_\alpha (t-t_\ell,y'-w)-\nabla_y p_\alpha (t-t_\ell,y-w) |_{w=z+b_h(s,z)(r-t_\ell)}\\ &\qquad\qquad\qquad\qquad\qquad\lesssim \frac{|y-y'|^\gamma}{(t-t_\ell)^{\frac{1+\gamma}{\alpha}}}(p_\alpha(t-t_\ell,y'-z)+p_\alpha(t-t_\ell,y-z)),
					\end{align*}
					yielding, along with \eqref{ineq-density-scheme} and $|\bar b_h |\leq |b|$,
					\begin{align*}
						\Delta_4 &\lesssim 	\frac{|y-y'|^\gamma}{(t-t_\ell)^{\frac{1+\gamma}{\alpha}}}\int_{t_\ell}^{t} \int p_\alpha (t_\ell - t_k,z-x)  (p_\alpha(t-t_\ell,y'-z)+p_\alpha(t-t_\ell,y-z)) |b (s,z)| \d z \d s.
					\end{align*}
					Let $$d_4:\{y,y'\}\ni\mathfrak{y}\mapsto\int_{t_\ell}^t \int p_\alpha (t_\ell - t_k,z-x) p_\alpha (t-t_\ell,\mathfrak{y}-z)|b (s,z) | \d z \d s,$$ so that $|\Delta_4|\lesssim\frac{|y-y'|^\gamma}{(t-t_\ell)^{\frac{1+\gamma}{\alpha}}} (d_4(y)+d_4(y'))$. Let us then bound $d_4$ using the convolution inequality \eqref{convo-space-sing}. For $\mathfrak{y} \in \{y,y'\}$,
					\begin{align*}
						d_4 (\mathfrak{y})&\lesssim \int_{t_\ell}^{t} \int p_\alpha (t_\ell-t_k,z-x) p_\alpha(t-t_\ell,\mathfrak{y}-z) |b (s,z)| \d z \d s\\
						&\lesssim 	\left[\frac{1}{(t_\ell-t_k)^{\frac{d}{\alpha p}}} +\frac{1}{(t-t_\ell)^{\frac{d}{\alpha p}}}  \right] p_\alpha (t-t_k,\mathfrak{y}-x)\int_{t_\ell}^{t} \Vert b (s,\cdot) \Vert_{L^{p}}\d s\\
						&\lesssim p_\alpha (t-t_k,\mathfrak{y}-x) \Vert b \Vert_{L^q-L^p} \Vert \mathbb{1}_{(t_\ell,t)}\Vert_{L^{q'}}   	\left[\frac{1}{(t_\ell-t_k)^{\frac{d}{\alpha p}}} +\frac{1}{(t-t_\ell)^{\frac{d}{\alpha p}}}  \right]\\
						&\lesssim p_\alpha (t-t_k,\mathfrak{y}-x) (t-t_\ell)^{1-\frac{1}{q}-\frac{d}{\alpha p}},
					\end{align*}
					where, for the last inequality, we used the fact that $t_\ell-t_k\geq t-t_\ell$. Plugging this into $|\Delta_4|$ yields
					\begin{align}\label{maj-holder-space-d4}
						|\Delta_4| & \nonumber \lesssim (p_\alpha (t-t_k,y-x)+p_\alpha (t-t_k,y'-x)) |y-y'|^\gamma (t-t_\ell)^{1-\frac{1}{q}-\frac{d}{\alpha p}-\frac{1}{\alpha}-\frac{\gamma}{\alpha}}\\
						& \lesssim (p_\alpha (t-t_k,y-x)+p_\alpha (t-t_k,y'-x)) |y-y'|^\gamma.
					\end{align}
					The estimates for $\Delta_3$ and $\Delta_4$ remain valid for $\bar \Delta_3$ and $\bar \Delta_4$  since we only used $|b_h|\lesssim |b|$. Plugging the estimations 
\eqref{maj-holder-space-d1}-\eqref{maj-holder-space-d4} into \eqref{splitting-holder-space-gammah} concludes the proof of \eqref{holder-space-gammah} and of Proposition \ref{prop-main-estimates}.
				\end{itemize}

	\section{Proof of existence of a unique weak solution and heat kernel estimates for the SDE \eqref{sde} (Theorem \ref{thm-wp-diffusion})}
	\subsection{Uniqueness of solutions to the Duhamel formulation \eqref{duhamel-diffusion} satisfying the estimation \eqref{ineq-density-diffusion}}\label{UNIQUE_MARG_LAW}Assume that $\Gamma_1$ and $\Gamma_2$ both satisfy the estimation \eqref{ineq-density-diffusion} and the Duhamel formula. Then
		$$\mu_x (t) := \sup_{y\in \R^d} \frac{|\Gamma_1 (0,x,t,y)-\Gamma_2(0,x,t,y)|}{p_\alpha (t,y-x)}.$$
is bounded on $(0,T]$ and we can write for all $(t,y)\in (0,T]\times \R^d$:
		\begin{align*}
			\Gamma_1 (0,x,t,y)-\Gamma_2(0,x,t,y) &= \int_0^t  \int b(r,z)\cdot \nabla_y p_\alpha (t-r,y-z) \left[ \Gamma_2 (0,x,r,z)-\Gamma_1 (0,x,r,z)\right]\d z \d r.
		\end{align*}
We deduce that for $(t,y)\in(0,T]\times\R^d$, 
		\begin{align*}
			\left| \frac{\Gamma_1 (0,x,t,y)-\Gamma_2(0,x,t,y)}{p_\alpha (t,y-x)} \right| &\leq \frac{1}{p_\alpha (t,y-x)}\int_0^t \int |b(r,z)| | \nabla_y p_\alpha (t-r,y-z) | p_\alpha (r,z-x) \mu_x (r)\d z \d r.
		\end{align*}
	Using \eqref{derivatives-palpha} and \eqref{convo-space-sing}, we get:
		\begin{align*}
			\left| \frac{\Gamma_1 (0,x,t,y)-\Gamma_2(0,x,t,y)}{p_\alpha (t,y-x)} \right| &\leq \int_0^t \frac{ \Vert b(r,\cdot)\Vert_{L^p} }{(t-r)^{\frac{1}{\alpha}}} \left[\frac{1}{r^{\frac{d}{\alpha p}}} + \frac{1}{(t-r)^{\frac{d}{\alpha p}}}\right] \mu_x (r) \d r.
		\end{align*}
		Taking the supremum over $y\in \R^d$ on the l.h.s. and applying H\"older's inequality in time,
                we get like in the last step of the proof of Theorem \ref{thm-main}\begin{align*}
	\forall t\in(0,T],\;	\mu_x(t)^{q'}\lesssim  \int_0^t   \frac{\mu_x(r)^{q'}}{(t-r)^{\frac{q'}{\alpha}}} \left[\frac{1}{r^{\frac{dq'}{\alpha p}}} + \frac{1}{(t-r)^{\frac{d q'}{\alpha p}}}\right] \d r .
	\end{align*}
Since $\frac{q'}{\alpha}+\frac{dq'}{\alpha p}<1$, Lemma 2.2 and Example 2.4 \cite{ZhangJFA10} ensure that 
$\forall t\in (0,T]$, $\mu_x(t)=0$, from which we immediately deduce $\Gamma_1=\Gamma_2$.
	\label{SEC_MART_PB}
	\subsection{Tightness of the laws $P^h$ of $((X^h_s)_{s\in [0,T]})_h$ and $\bar P^h $ of  $((\bar X^h_s)_{s\in [0,T]})_h$}
        \label{subsec-tightness}
        Let ${\mathfrak{B}}^\eta:=\E_{0,x} \left[\int_0^T \left| b_h(U_{\lfloor\frac{s}{h}\rfloor},X_{\tau_s^h}^h)\right|^\eta \d s \right]$, where $\eta>1 $ is chosen sufficiently close to $1$ in order that, under \eqref{serrin}, $p/\eta > 1$, $q/\eta> 1 $ and $ \eta (d/p+\alpha/q)<\alpha$.
        Using $|b_h|\le Bh^{-\frac{d}{\alpha p}-\frac 1 q}$ on the first time step then $|b_h|\le |b|$, \eqref{ineq-density-scheme}, H\"older's inequality and \eqref{spatial-moments} with $\delta=0 $, we obtain
    \begin{align*}
			{\mathfrak{B}}^\eta& \lesssim h^{1-\frac{d}{\alpha p}-\frac 1 q}+\sum_{k=1}^{n-1} \int_{t_k}^{t_{k+1}}\frac{1}{h}\int_{t_k}^{t_{k+1}} \int \left| b(r,y)\right|^\eta \Gamma^h (0,x,t_k,y) \d y\d r \d s \\
			& \lesssim  h^{1-\frac{d}{\alpha p}-\frac 1 q}+ \sum_{k=1}^{n-1} \int_{t_k}^{t_{k+1}} \int \left| b(r,y)\right|^\eta p_\alpha (t_k,y-x) \d y\d r \\
			& \lesssim  h^{1-\frac{d}{\alpha p}-\frac 1 q}+\sum_{k=1}^{n-1} \int_{t_k}^{t_{k+1}}  \Vert b(r,\cdot)^\eta \Vert_{L^{\frac{p}{\eta}}} \Vert p_\alpha (t_k,\cdot-x)\Vert_{L^{(\frac{p}{\eta})'}} \d r \\
			& \lesssim  h^{1-\frac{d}{\alpha p}-\frac 1 q}+ \sum_{k=1}^{n-1}  \int_{t_k}^{t_{k+1}}  t_k^{-\frac{d\eta }{ \alpha p}}  \Vert b(r,\cdot) \Vert_{L^{p}}^\eta\d r.
		\end{align*}    
We then write
		$$\int_{t_k}^{t_{k+1}}  t_k^{-\frac{d \eta}{\alpha p}}  \Vert b(r,\cdot) \Vert_{L^{p}}^\eta\d r  \leq \int_{t_k}^{t_{k+1}}  r^{-\frac{d \eta }{\alpha p}}  \Vert b(r,\cdot) \Vert_{L^{p}}^\eta\d r,$$
		and use a H\"older inequality to obtain from the condition $\eta(d/p+\alpha/q)<\alpha $ which ensures that $-(\frac{q}{\eta})'\frac{d\eta}{\alpha p}>-1$,
		\begin{align}\label{control-b}
			{\mathfrak{B}}^\eta 
                  &\lesssim   \left( \int_{t_1}^{t_{n}}  r^{-(\frac{q}{\eta})'\frac{d\eta}{\alpha p}}  \d  r\nonumber \right)^{\frac{1}{(\frac{q}{\eta})'}}\times \Vert r \mapsto \Vert b(r,\cdot) \Vert_{L^{p}}^\eta \Vert_{L^{\frac{q}{\eta}}}\nonumber\\
			&\lesssim h^{1-\frac{d}{\alpha p}-\frac 1 q}+T^{1-\eta(\frac{d}{ \alpha p}+\frac{1}{q} )}\|b\|_{L^q-L^p}^\eta.
		\end{align}
                In the same way, $$\bar{\mathfrak{B}}^\eta:=\E_{0,x} \left[\int_0^T \left| \bar b_h(U_{\lfloor\frac{s}{h}\rfloor},\bar X_{\tau_s^h}^h)\right|^\eta \d s \right]=\E_{0,x} \left[\int_h^T \left|  \bar  b_h(U_{\lfloor\frac{s}{h}\rfloor},\bar  X_{\tau_s^h}^h)\right|^\eta \d s\right]\lesssim T^{1-\eta(\frac{d}{ \alpha p}+\frac{1}{q} )}\|b\|_{L^q-L^p}^\eta.$$
 By H\"older's inequality, we have
		\begin{align}\label{control-B}
			\forall 0\leq u \leq t \leq T, \left|\int_{u}^{t} b_h(U_{\lfloor\frac{s}{h}\rfloor},X_{\tau_s^h}^h) \d s\right| \lesssim (t-u)^{\frac{\eta-1}{\eta}} \left(\int_0^T \left| b_h(U_{\lfloor\frac{s}{h}\rfloor},X_{\tau_s^h}^h)\right|^\eta \d s \right)^{\frac{1}{\eta}}.
		\end{align}
and the same estimation holds with $(b_h,X^h)$ replaced by $(\bar b_h,\bar X^h)$. 
		Since by \eqref{spatial-moments} applied with $\delta=1$, $\E[|Z_t-Z_u|]\lesssim (t-u)^{\frac 1\alpha}$, setting $\zeta=\left(1-\frac 1\eta\right)\wedge \frac 1\alpha>0$, we deduce that \begin{equation}
                    \forall 0\leq u \leq t \leq T,\;\E [|X_t^h-X_u^h|]+\E [|\bar X_t^h-\bar X_u^h|]\lesssim (t-u)^\zeta\label{estimom},
                \end{equation}
        which ensures the tightness of the laws $P^h$ of $X^h$ and $\bar P^h$ of $(\bar X^h)$ on the space $\mathcal D([0,T],\R^d)$ of càdlàg functions endowed with the Skorokhod topology (see Proposition 34.9 from \cite{Bas11} for example). Let $\big(\xi_s\big)_{s\in [0,T]}$ denote the canonical process on this space.\\
		We may then extract a subsequence, still denoted by $(P^h)$ (resp. $(\bar P^h)$), such that $P^h$ (resp. $(\bar P^h)$) weakly converges to some limit probability $P$ on $\mathcal D([0,T],\R^d)$ as $h\rightarrow 0$. For $u,t\in[0,T]$ outside the at most countable set $\{s\in (0,T]:P(|\xi_s-\xi_{s-}|>0)>0\}$, the law of  $({X}^h_u,{X}^h_t)$ (resp. $(\bar{X}^h_u,\bar{X}^h_t)$) converges to $P\circ(\xi_u,\xi_t)^{-1}$ so that \eqref{estimom} combined with the right-continuity of sample-paths ensures that $\sup_{0\le u<t\le T}(t-u)^{-\zeta}\int_{{\cal D}([0,T],\R^d)}|\xi_t-\xi_u|P(\d\xi)<\infty$. As a consequence \begin{equation}
                  \{s\in (0,T]:P(|\xi_s-\xi_{s-}|>0)>0\}=\emptyset\label{tempsdisc}\end{equation} and for each $t\in(0,T]$, the distribution $\Gamma^h(0,x,t,y) \d y$ of $X^h_t$ (resp. $\bar\Gamma^h(0,x,t,y)\d y$ of $\bar X^h_t$) converges weakly to $P_t=P\circ \xi_t^{-1}$. By \eqref{ineq-density-scheme}
and \eqref{holder-space-gammah}, the Ascoli-Arzelà theorem ensures that we can extract a further subsequence such that $y\mapsto \Gamma^h(0,x,t,y)$ (resp. $y\mapsto \bar\Gamma^h(0,x,t,y)$) converges uniformly on the compact subsets of $\R^d$ to some limit $y\mapsto \Gamma(0,x,t,y)$ so that $P_t(\d y)=\Gamma(0,x,t,y)\d y$. Taking the limit $h\to 0$ into \eqref{ineq-density-scheme}
and \eqref{holder-space-gammah} ensures that $\Gamma$ satisfies \eqref{ineq-density-diffusion} and \eqref{holder-space-gamma}.

We are next going to prove that the limit probability measure $P$ solves the following martingale problem.

\begin{definition}[Martingale Problem]\label{DEF_MART_PB}
A probability measure $ P$  on the space $\mathcal D([0,T],\R^d)$ of càdlàg functions with time-marginals $(P_t)_{t\in[0,T]}$, solves the  martingale problem related to $b\cdot \nabla +\mathcal L^\alpha$ and $x\in \R^d$ if :
\begin{itemize}
\item[(i)] $ P_0=\delta_x$,
\item[(ii)] for a.a. $t\in(0,T]$, $P_t(\mathrm{d}y)=\rho(t,y)\d y$ for some $\rho \in L^{q'}((0,T],L^{p'}(\R^d))$,
\item[(iii)] for all $\mathcal C^{1,2}$ function $f$ on $[0,T]\times\mathbb R^d$ bounded together with its derivatives
 , the process
\begin{equation}\label{M}
\Bigg\{M_t^f=f(t,\xi_t)-f(0,\xi_0)-\int_0^t \Big(\textcolor{black}{(\partial_s +\mathcal L^\alpha)f(s,\xi_s)}+ b (s,\xi_s) \cdot \nabla
 f(s,\xi_s)\Big)\d s\Bigg\}_{0\le t\le T},\tag{${\rm M}$}
\end{equation}
is a $ P$ martingale. 
\end{itemize}
\end{definition}
Let us point out that, in the current singular drift setting, condition $(ii)$ which guarantees that $$\int_{\mathcal D([0,T],\R^d)}\int_0^T|b(s,\xi_s)|\d s P(\mathrm{d}\xi)<\infty$$ is somehow the minimal one required for  	all the terms in \eqref{M} to be well defined. 

Before checking that the limit probability measure $P$ solves the martingale problem, let us prove that this implies that $\Gamma$ solves \eqref{duhamel-diffusion}, which concludes the proof of Theorem \ref{thm-wp-diffusion} (in fact, for this purpose, it would be enough to check that the limit probability measure associated with either the schemes $X^h$ or the schemes $\bar X^h$ solves the martingale problem).
	     Let $t\in (0,T]$ and $\phi:\R^d\to\R$ be a ${\cal C}^4$ function with compact support. Choosing $f(s,z)={\mathbb 1}_{[0,t)}(s)p_\alpha(t-s,\cdot)\star\phi(z)+{\mathbb 1}_{[t,T]}(s)(\phi(z)-(s-t){\cal L}^\alpha\phi(z))$ which, \textcolor{black}{according to Lemma \ref{lemFK},} satisfies $(\partial_s+{\cal L}^\alpha)f(s,z)=0$ for $(s,z)\in[0,t]\times\R^d$ and writing the centering  of  $M^f_t$ (introduced in Definition \ref{DEF_MART_PB}) under $P$, we obtain that
        \begin{align*}
         \int_{\R^d}\phi(y)\Gamma(0,x,t,y)\d y=\int_{\R^d}\phi(y)p_{\alpha}(t,x-y)\d y+\int_0^t\int_{\R^d}\Gamma(0,x,s,z) b(s,z)\cdot\nabla_zf(s,z)\d z \d s.
        \end{align*}
        Using \eqref{ineq-density-diffusion} and \eqref{derivatives-palpha} to justify the use of Fubini's theorem and the fact that for $s\in(0,T]$, $p_\alpha(s,\cdot)$ is an even function, we deduce that
  \begin{align*}
   \int_{\R^d}\phi(y)\Gamma(0,x,t,y)\d y=\int_{\R^d}\phi(y)\left(p_{\alpha}(t,y-x)-\int_0^t\Gamma(0,x,s,z) b(s,z)\cdot\nabla_yp_\alpha(s,y-z)\d s\right)\d y
  \end{align*}
  Since $\phi$ is arbitrary, we conclude that $(0,T]\times\R^d\ni(t,y)\mapsto \Gamma(0,x,t,y)$ satisfies \eqref{duhamel-diffusion}.
  \subsection{Any limit point solves the martingale problem}
  Let us now prove that the limit point $P$ solves the martingale problem associated with \eqref{sde} and introduced in Definition \ref{DEF_MART_PB}.  Since for each $h$, $X^h_0=x=\bar{X}^h_0$, one has $P_0=\delta_x$. Moreover, for $t\in (0,T]$, $P_t(dy)=\Gamma(0,x,t,y)dy$ with $\Gamma$ satisfying \eqref{ineq-density-diffusion}. By \eqref{spatial-moments} applied with $\delta=0$, $\|\Gamma(0,x,t,\cdot)\|_{L^{p'}}\le C t^{-\frac{d}{\alpha p}}$ where the right-hand side belongs to $L^{q'}([0,T])$ since $q'\frac{d}{\alpha p}<1$ by \eqref{serrin}. As a consequence, $\Gamma(0,x,\cdot,\cdot)\in L^{q'}((0,T],L^{p'}(\R^d))$. Therefore properties $(i)$ and $(ii)$ in Definition \ref{DEF_MART_PB} hold.
  
	 Let $f:[0,T]\times\mathbb R^d\to\R$ be $\mathcal C^{1,2}$ and bounded together with its derivatives, $\psi:(\R^d)^p \rightarrow \R$ be continuous and bounded, $0\leq s_1\leq...\leq s_p < u \leq t \leq T$ with $u>0$ and $F:\mathcal D([0,T],\R^d)\rightarrow \R$ be  defined by
	\begin{equation}\label{DEF_F}
	F(\xi):=\left(f(t,\xi_t)-f(u,\xi_u) - \int_u^t \left[(\partial_s+\mathcal{L}^\alpha)f (s,\xi_s)+b(s,\xi_s)\cdot \nabla f(s,\xi_s)\right] \d s\right)\psi (\xi_{s_1},...,\xi_{s_p}).
	\end{equation}
	In order to prove that $P$ satisfies $(iii)$ in Definition \ref{DEF_MART_PB}, we will show that {$\int_{{\mathcal D}([0,T],\R^d)}F(\xi)P(\d \xi)=0$.} 
        
	\begin{paragraph}{Proof of $\lim_{h\rightarrow 0} \E [F( X^h)]=0=\lim_{h\rightarrow 0} \E [F(\bar X^h)]$. \\[0.5cm]}
		Using Itô's formula, we can write
		$$f(t, X_t^h)-f(u, X_u^h) =  M_{t}^h - M_{u}^h+ \int_u^t \nabla f(s, X_s^h)\cdot  b_h (U_{\lfloor s/h\rfloor}, X_{\tau_s^h}^h) \d s+\int_u^t  (\partial_s+\mathcal{L}^\alpha)f(s, X_s^h) \d s,$$
		with $ M_{s}^h=\int_0^s \int_{\R^d\backslash\{0\}}\Big(f(r, X_{r^-}^h+x)-f(r, X_{r^-}^h)\Big) \tilde N(\d r,\d x)$ where $\tilde N $ is the compensated Poisson measure associated with $Z$.
		Since $ M^h$ is a martingale, taking expectations, we get:
		\begin{align}
			\E [F( X^h)]&=\E \left[\left(\int_u^t \left( b_h (U_{\lfloor s/h\rfloor}, X_{\tau_s^h}^h) -b(s, X_s^h)\right) \cdot \nabla f(s, X_s^h) \d s\right) \psi ( X_{s_1}^h,..., X_{s_p}^h)\right]\notag\\
			&=\E \left[\left(\int_u^t  b_h (U_{\lfloor s/h\rfloor}, X_{\tau_s^h}^h) \cdot \left( \nabla f(s, X_s^h)-\nabla f(s, X_{\tau_s^h}^h)\right) \d s\right) \psi ( X_{s_1}^h,..., X_{s_p}^h)\right]\notag\\
			&\qquad + \E \left[\left(\int_u^t \left( b_h (U_{\lfloor s/h\rfloor}, X_{\tau_s^h}^h) -  b_h(s, X_{\tau_s^h}^h)\right) \cdot \nabla f(s, X_{\tau_s^h}^h) \d s\right) \psi ( X_{s_1}^h,..., X_{s_p}^h)\right]\notag\\
			&\qquad +\E \left[\left(\int_u^t \left( b_h (s, X_{\tau_s^h}^h) -b(s, X_{\tau_s^h}^h)\right) \cdot \nabla f(s, X_{\tau_s^h}^h) \d s\right) \psi ( X_{s_1}^h,..., X_{s_p}^h)\right]\notag\\
			&\qquad +\E \left[\left(\int_u^t \left(b (s, X_{\tau_s^h}^h)\cdot\nabla f(s, X_{\tau_s^h}^h)-b(s, X_s^h)\nabla f(s, X_s^h)\right)   \d s\right) \psi ( X_{s_1}^h,..., X_{s_p}^h)\right]\notag\\
			&=: \Delta_1+\Delta_2+\Delta_3+\Delta_4.\label{DECOUP_DELTA_MART}
		\end{align}
In the same way, $\E [F(\bar X^h)]=\bar\Delta_1+\bar\Delta_2+\bar\Delta_3+\bar\Delta_4$ where $\bar\Delta_i$ is defined like $\Delta_i$ with $(X^h,b_h)$ replaced by $(\bar X^h,\bar b_h)$ for $i\in\{1,\cdots,4\}$.

		For $\Delta_1$, we first write, using $|b_h|\lesssim h^{-\frac{d}{\alpha p}-\frac 1q}  $ and conditioning w.r.t. 
		${\mathcal F}_{\tau_s^h}=\sigma\big(X^h_u,\ 0\le u\le \tau_s^h \big) $,
		\begin{align*}
			\E_{ {\mathcal F}_{\tau_s^h}} |\nabla f(s, X_s^h) -\nabla f(s, X_{\tau_s^h}^h)|& \leq \Vert \nabla^2 f \Vert_{L^\infty}\E_{ {\mathcal F}_{\tau_s^h}} | X_{s}^h- X_{\tau_s}^h|\\
			& \leq \Vert \nabla^2 f \Vert_{L^\infty}\E_{ {\mathcal F}_{\tau_s^h}} \left[ \int_{\tau_s^h}^{s}  |  b_h (U_{\lfloor r/h\rfloor}, X_{\tau_r^h}^h) | \d r + |Z_s-Z_{\tau_s^h}| \right]\\
			& \lesssim \Vert \nabla^2 f \Vert_{L^\infty} \left[ \int_{\tau_s^h}^{s}  h^{-\frac{d}{\alpha p}-\frac 1q}dr + (s-\tau_s^h)^{\frac{1}{\alpha}}\right]\\
			& \lesssim\Vert \nabla^2 f \Vert_{L^\infty} h^{\frac 1\alpha}.
		\end{align*}
		Using this bound along with $|b_h|\le |b|$ and \eqref{control-b}, we can compute
		\begin{align*}
                  |\Delta_1| \lesssim \Vert \psi \Vert_{L^\infty}  \Vert \nabla^2 f \Vert_{L^\infty} h^{\frac 1\alpha}\E\left[\int_{0}^T |b(U_{\lfloor s/h\rfloor}, X_{\tau_s^h}^h)|\d s \right]	\lesssim \Vert \psi \Vert_{L^\infty}  \Vert \nabla^2 f \Vert_{L^\infty} h^{\frac 1\alpha}.	\end{align*}
                The same bound holds for $|\bar\Delta_1|$ since the larger cutoff $|\bar b_h|\lesssim h^{\frac{1}{\alpha}-1}$ does not deteriorates the estimation of $\E_{\bar {\mathcal F}_{\tau_s^h}} |\nabla f(s,\bar X_s^h) -\nabla f(s,\bar X_{\tau_s^h}^h)|$ where $\bar {\mathcal F}_{\tau_s^h}=\sigma\big(\bar X_u^h,\ 0\le u\le \tau_s^h \big) $.
                      
		For $\Delta_2$, supposing that $h$ is small enough to ensure that $\tau^h_u<\tau^h_t$, we split the time integral into three terms: a main term over $(\tau_u^h+h,\tau_t^h)$ which matches the time grid, and two terms around the edges, over $(u,\tau_u^h+h)$ and $(\tau_t^h,t)$ respectively. For the main term, we will use the following cancellation:
		$$\E \left[\int_{\tau_u^h+h}^{\tau_t^h} \left( b_h (U_{\lfloor s/h\rfloor}, X_{\tau_s^h}^h) - b_h(s, X_{\tau_s^h}^h)\right) \cdot \nabla f(s, X_{\tau_s^h}^h) \d s \; \bigg| \;  {\mathcal F}_{\tau_u^h+h} \right]=0.$$
		For the other two terms, we use that $| b_h (U_{\lfloor s/h\rfloor}, X_{\tau_s^h}^h) - b_h(s, X_{\tau_s^h}^h)|\leq | b_h (U_{\lfloor s/h\rfloor}, X_{\tau_s^h}^h)|+|  b_h(s, X_{\tau_s^h}^h)|\lesssim h^{-\frac{d}{\alpha p}-\frac 1q}$ and the inequalities $\tau_u^h+h-u\leq h$ and $t-\tau_t^h\leq h$ to write
		\begin{align*}
			|\Delta_2|&\lesssim \Vert \psi \Vert_{L^\infty} \Vert \nabla f\Vert_{L^{\infty}} \E \left[  \int_{[u,\tau_u^h+h] \cup[\tau_t^h,t]} \left(| b_h (U_{\lfloor s/h\rfloor}, X_{\tau_s^h}^h) |+| b_h(s, X_{\tau_s^h}^h)| \right) \d s\right]\\
			& \lesssim\Vert \psi \Vert_{L^\infty} \Vert \nabla f \Vert_{L^{\infty}} h^{1-\frac{d}{\alpha p}-\frac 1q} .
		\end{align*}
In the same way, $|\bar \Delta_2|\lesssim\Vert \psi \Vert_{L^\infty} \Vert \nabla f \Vert_{L^{\infty}} h^{\frac{1}{\alpha}}$.

      When $p=q=\infty$, $\Delta_3$ vanishes. Otherwise, applying \eqref{maj-delta2-main} with $\lambda=\eta$ where $\eta>1$ is such that \eqref{control-b} holds, we get
 		\begin{align*}
			|\Delta_3| \leq  \Vert\psi \Vert_{L^\infty}\Vert \nabla f \Vert_{L^\infty} h^{(\frac d{\alpha p}+\frac 1 q)(\eta-1)}  \E \left[\int_u^t |b(s, X_{\tau_s^h}^h)|^\eta \d s \right] \lesssim \Vert\psi \Vert_{L^\infty}\Vert \nabla f \Vert_{L^\infty} h^{(\frac d{\alpha p}+\frac 1 q)(\eta-1)}.
		\end{align*}
 Since $|b-\bar b^h|\leq |b|^\eta B^{1-\eta} h^{(1-\frac{1}{\alpha })(\eta-1)}$, we obtain in the same way that $|\bar\Delta_3|\lesssim \Vert\psi \Vert_{L^\infty}\Vert \nabla f \Vert_{L^\infty} h^{(1-\frac 1\alpha)(\eta-1)}$.
 For $\Delta_4$,  we have
 \begin{align*}
   &\left|\E_{{\cal F}_{\tau^h_u-h}}\left[\int_u^t \left(b (s, X_{\tau_s^h}^h)\cdot\nabla f(s, X_{\tau_s^h}^h)-b(s, X_s^h)\nabla f(s, X_s^h)\right)   \d s\right]\right|\\&\le \int_u^t \int \left| b (s,z)\cdot\nabla  f (z) \right| \left|  \Gamma^h (\tau^h_u-h,X^h_{\tau^h_u-h},\tau_s^h,z) - \Gamma^h (\tau^h_u-h,X^h_{\tau^h_u-h},s,z)\right|  \d z \d s .
 \end{align*}
 Assuming w.l.o.g. that $h $ is small enough to have $\tau_u^h-h\ge s_p $, 
 we deduce that:
		\begin{align*}
			|\Delta_4 |&\lesssim \Vert\psi \Vert_{L^\infty} \int_u^t \int \int \left| b (s,z)\cdot\nabla  f (z) \right| \left|  \Gamma^h (\tau^h_u-h,y,\tau_s^h,z) - \Gamma^h (\tau^h_u-h,y,s,z)\right|  \Gamma^h (0,x,\tau_u^h-h,y) \d z \d y \d s .
		\end{align*}
		Then, we use the H\"older regularity \eqref{holder-time-gammah} of $ \Gamma^h$ in the forward time variable:
		\begin{align*}
			|\Delta_4 |&\lesssim \Vert\psi \Vert_{L^\infty} \Vert \nabla  f \Vert_{L^{\infty}}\int_u^t \int \int |b (s,z)| \frac{(s-\tau_s^h)^\frac{\gamma}{\alpha}}{(\tau_s^h-\tau_u^h+h)^{\frac{\gamma}{\alpha}}} p_\alpha (s-\tau_u^h+h,z-y) \Gamma^h (0,x,\tau_u^h-h,y) \d z \d y \d s.
		\end{align*}
		Since $s-\tau_s^h \leq h$ and $\tau_s^h-\tau_u^h+h>s-u$, we get, using \eqref{ineq-density-scheme}, H\"older's inequality and \eqref{spatial-moments}
		\begin{align*}
			|\Delta_4 |&\lesssim \Vert\psi \Vert_{L^\infty} \Vert \nabla  f \Vert_{L^{\infty}} h^\frac{\gamma}{\alpha} \int_u^t \int \int \frac{|b (s,z)| }{(s-u)^{\frac{\gamma}{\alpha}}} p_\alpha (s-\tau_u^h+h,z-y)  p_\alpha (\tau_u^h-h,y-x)  \d z \d y \d s\\
			&\lesssim \Vert\psi \Vert_{L^\infty} \Vert \nabla  f \Vert_{L^{\infty}} h^\frac{\gamma}{\alpha} \int_u^t \int  \frac{|b (s,z)| }{(s-u)^{\frac{\gamma}{\alpha}}} p_\alpha (s,z-x) \d z \d s\\
			&\lesssim \Vert\psi \Vert_{L^\infty} \Vert \nabla  f \Vert_{L^{\infty}} h^\frac{\gamma}{\alpha} \int_u^t  \frac{\Vert b (s,\cdot)\Vert_{L^{p}} }{(s-u)^{\frac{\gamma}{\alpha}}s^{\frac{d}{\alpha p}}}\d s.
		\end{align*}
		Finally, using H\"older's inequality in time and $\frac 1{q'}-\frac{\gamma}{\alpha}-\frac d{\alpha p}=\frac 1\alpha $, we obtain
		\begin{align*}
			|\Delta_4| \lesssim \Vert\psi \Vert_{L^\infty} \Vert \nabla  f \Vert_{L^{\infty}} \Vert b \Vert_{L^{q}-L^p }h^\frac{\gamma}{\alpha} (t-u)^{\frac{1}{\alpha}}.
		\end{align*}
        The same estimation holds for $|\bar\Delta_4|$.
		Putting together the previous estimates on $(\Delta_i)_{i\in \{1,\cdots,4\}} $ in \eqref{DECOUP_DELTA_MART} \textcolor{black}{we obtain} $\lim_{h\rightarrow 0} \E \left[F(X^h)\right] = 0$. In the same way, $\lim_{h\rightarrow 0} \E \left[F(\bar X^h)\right] = 0$
		\end{paragraph}
	
		\begin{paragraph}{$P$ solves the martingale problem.\\[0.5cm]}
	        In this paragraph, we only consider the case when $P$ is the limit of the laws of the schemes $X^h$ since the argument is exactly the same when $P$ is the limit of the laws of the schemes $\bar X^h$. The lack of continuity of the functional $F$ on ${\mathcal D}([0,T],\R^d)$ prevents from deducing immediately that $\int_{{\mathcal D}([0,T],\R^d)}F(\xi)P(\mathrm{d} \xi)=0$. Let us first suppose that $p<\infty$ and set $\tilde q=q{\mathbb 1}_{q<\infty}+\frac{\alpha p +1}{(\alpha-1)p-d}{\mathbb 1}_{q<\infty}$. We have $\frac{d}{p}+\frac{\alpha}{\tilde q}<\alpha-1$.  We introduce for $\varepsilon\in(0,1]$,  a smooth and bounded function $b_\varepsilon $  such that $\lim_{\varepsilon\to 0}\|b_\varepsilon-b\|_{L^{\tilde q}-L^{p}}=0$.
			The functional $F_\varepsilon$ defined like $F$ in \eqref{DEF_F}, but with $b_\varepsilon$ replacing $b$ is bounded. According to \eqref{tempsdisc}, for fixed $\varepsilon\in (0,1]$, $P$ gives full weight to continuity points of $F_\varepsilon$ and since $\lim_{h\to 0}\E[F( X^h)]=0$, we have
			$$\int_{{\mathcal D}([0,T],\R^d)}F_\varepsilon(\xi)P(\mathrm{d} \xi)=\lim_{h\to 0}\E[F_\varepsilon( X^h)]=\lim_{h\to 0}\E[F_\varepsilon( X^h)-F( X^h)].$$
			We deduce that
			$$\left|\int_{{\mathcal D}([0,T],\R^d)}F(\xi)P(\mathrm{d} \xi)\right|\le \limsup_{\varepsilon\to 0}\int_{{\mathcal D}([0,T],\R^d)}|F(\xi)-F_\varepsilon(\xi)|P(\mathrm{d} \xi)+\limsup_{\varepsilon\to 0}\limsup_{h\to 0}\E[|F_\varepsilon( X^h)-F( X^h)|].$$
			One has, using \eqref{ineq-density-scheme}, then H\"older's inequality in space together with \eqref{spatial-moments} applied with $(\ell,\delta)=(p,0)$ and last H\"older's inequality in time,
			\begin{align*}
				\E[|F_\varepsilon( X^{h})-F( X^{h})|]&\le \|\psi\|_{L^\infty}\|\nabla f\|_{L^\infty}\int_u^t\E[|b_\varepsilon(s, X^{h}_s)-b(s, X^{h}_s)
				|]\d s\notag\\
				&\lesssim \|\psi\|_{L^\infty}\|\nabla f\|_{L^\infty}\int_u^t\int |b_\varepsilon(s,y)-b(s,y)|p_\alpha(s,y-x) \d y \d s \notag\\
				&\lesssim \|\psi\|_{L^\infty}\|\nabla f\|_{L^\infty}\int_u^t\frac{\|b_\varepsilon(s,.) -b(s,.)\|_{L^{ p}}}{s^{\frac{d}{\alpha  p}}} \d s\\
				&\lesssim \|\psi\|_{L^\infty}\|\nabla f\|_{L^\infty}\|b_\varepsilon-b\|_{L^{\tilde q}-L^{p}}t^{1-\left(\frac 1 {\tilde q}+\frac d{ \alpha  p}\right)}.
			\end{align*}
	        Since the same estimation holds for $\int_{{\mathcal D}([0,T],\R^d)}|F(\xi)-F_\varepsilon(\xi)|P(\mathrm{d} \xi)$, because the heat kernel estimates hold as well for the limit point, we conclude that $\int_{{\mathcal D}([0,T],\R^d)}F(\xi)P(\mathrm{d}\xi)$ $=0$. Taking $ f,\psi,u,s_1,\ldots,s_p,t$ in countable dense subsets, we deduce that $P$ satisfies  $(iii)$ in Definition \ref{DEF_MART_PB}.\\
	                      
	        Let us now deal with the case $p=\infty$. We set $(\tilde p,\tilde q)=(\frac{dq+1}{(\alpha-1)q-\alpha},q){\mathbb 1}_{q<\infty}+(\frac{3d}{\alpha-1},\frac{3\alpha}{\alpha-1}){\mathbb 1}_{q=\infty}$. We have $\frac{d}{\tilde p}+\frac{\alpha}{\tilde q}<\alpha-1$.  We introduce for $\varepsilon\in(0,1]$,  a smooth and bounded function $b_\varepsilon $  such that $\|b_\varepsilon\|_{L^{q}-L^{\infty}}\le 2\|b\|_{L^{q}-L^{\infty}}$ and, for each $K\in\N^*$, setting $b_\varepsilon^K(t,x)=\mathbb{1}_{[-K,K]^d}(x)b_\varepsilon(t,x)$ and $b^K(t,x)=\mathbb{1}_{[-K,K]^d}(x)b(t,x)$, we  have $\lim_{\varepsilon\to 0}\|b^K_\varepsilon-b^K\|_{L^{\tilde q}-L^{\tilde p}}=0$. The above reasoning when $p<\infty$ remains valid once we now bound $\E[|F_\varepsilon( X^{h})-F( X^{h})|]$ from above by
	    	\begin{align*}
				&\|\psi\|_{L^\infty}\|\nabla f\|_{L^\infty}\int_u^t\E[|b^K_\varepsilon(s, X^{h}_s)-b^K(s, X^{h}_s)
				|+(\|b_\varepsilon(s,\cdot)\|_{L^\infty}+\|b(s,\cdot)\|_{L^\infty}){\mathbb 1}_{|X^h_s|\ge K}]\d s\notag\\
				&\lesssim \|\psi\|_{L^\infty}\|\nabla f\|_{L^\infty}\left(\|b^K_\varepsilon-b^K\|_{L^{\tilde q}-L^{\tilde p}}t^{1-\left(\frac 1 {\tilde q}+\frac d{ \alpha  \tilde p}\right)}+\|b\|_{L^q-L^\infty}\left(\int_u^t(\P(|X^h_s|\ge K))^{q'}\d s\right)^{1/q'}\right).
			\end{align*}                  
			According to \eqref{estimom}, $\int_u^t(\P(|X^h_s|\ge K))^{q'}\d s$ can be made arbitrarily small uniformly in $h$ for $K$ large enough while  for fixed $K$, $\|b^K_\varepsilon-b^K\|_{L^{\tilde q}-L^{\tilde p}}$ goes to $0$ with $\varepsilon$. This concludes the proof. 
	
		\end{paragraph}
	\subsection{Uniqueness of the solution to the martingale problem}
		For this paragraph, we assume $p,q<\infty$ (otherwise, we can proceed in a similar way to the previous paragraph to mollify the drift). Let $(b_m)_{m\in \N}$ denote a sequence of bounded smooth approximating functions s.t. $\Vert b-b_m\Vert_{L^q-L^p}\rightarrow 0$ as $m\rightarrow \infty$. We study the mollified equation
		\begin{equation}\label{pde-mol}
			(\partial_s + \mathcal{L}^\alpha + b_m\cdot\nabla)u_m(s,x)=f(s,x), \qquad (s,x)\in [0,t)\times \R^d, u_m(t,\cdot)=0.
		\end{equation}
		It is well known that for a smooth compactly supported $f$, \eqref{pde-mol} has a unique smooth bounded classical solution (see \cite{MP14}). Furthermore, the following Schauder estimates (whose proofs are postponed to Appendix \ref{subsec-schauder-proofs}) hold:
		\begin{lemma}[Schauder]\label{lemma-schauder} Let $f:[0,T]\times\R^d \to\R$ be ${\cal C}^{1,2}$ with compact support and $(u_m)_{m\in \N}$ denote the sequence of classical solutions to the mollified PDEs \eqref{pde-mol}. Then, for all $\xi \in [0,(\gamma+1)/\alpha)$, for all $0\leq s\leq s'\leq t$, for all $x\in \R^d$, and for all $m\in \N$,
			\begin{align}
				\Vert \nabla u_m \Vert_{L^\infty} &\lesssim  \Vert f \Vert_{L^\infty}\label{gradient-bound},\\
				|u_m(s',x)-u_m(s,x)|&\lesssim |s'-s|^{\xi} \Vert f \Vert_{L^\infty}.\label{time-holder}
			\end{align}
		\end{lemma}	
		Let $P^1$ and $P^2$ be solutions of the martingale problem associated with $b\cdot\nabla + \mathcal{L}^\alpha$ and $x\in \R^d$ in the sense of Definition \ref{DEF_MART_PB}. Let $f$ be a smooth bounded function. For all $m\in \N$, denote $u_m\in \mathcal{C}^1([0,T],\mathcal{C}^2(\R^d,\R))$ the classical solution to the Cauchy problem associated with \eqref{pde-mol} with source term $f$. For $i\in \{ 1,2\}$, 
		\begin{equation}\label{uniqueness-1}
			\left\{M_s^{u_m}={u_m}(s,\xi_s)-{u_m}(0,x)-\int_0^s (\partial_r +\mathcal L^\alpha+b\cdot\nabla){u_m}(r,\xi_r)\d r\right\}_{0\le s\le t}
		\end{equation}
		is a $P^i$-martingale. Equations \eqref{time-holder} and \eqref{gradient-bound} allow us to apply the Ascoli-Arzel\`a theorem to $(u_m)$: let  $(u_{m_k})_k$ be a subsequence of $(u_m)_m$ which converges uniformly on every compact subset of $[0,t]\times \R^d$ to some $u_\infty$.
		Now, taking $\rho^i \in L^{q'}((0,T],L^{p'}(\R^d))$ such that $P_t^i(\mathrm{d}y)=\rho^i(t,y)\d y$, taking expectations in \eqref{uniqueness-1} and using the Fubini theorem, we have, when $s\rightarrow t$,
		\begin{equation}\label{uniqueness-expectation}
			\E^{P^i}\left[ \int_0^t f(r,\xi_r) \d r\right]=-{u_{m_k}}(0,x)+ \int_0^t \int  (b_{m_k}-b)(r,z) \cdot\nabla{u_{m_k}}(r,z)\rho^i (r,z)\d z\d r.
		\end{equation}
		Since $\rho^i\in L^{q'}((0,T],L^{p'}(\R^d))$, using a H\"older inequality in space and then one in time along with the fact that $\Vert \nabla u_{m_k} \Vert_{L^\infty}$ is bounded uniformly in $k$ (from Equation \eqref{gradient-bound}), we obtain
		\begin{align*}
			\left|\int_0^t \int  (b_{m_k}-b)(r,z) \cdot\nabla{u_{m_k}}(r,z)\rho^i (r,z)\d z\d r \right|\lesssim \Vert \nabla u_{m_k} \Vert_{L^\infty} \Vert b-b_{m_k} \Vert_{L^q-L^p} \Vert \rho^i \Vert_{L^{q'}-L^{p'}} \underset{k\rightarrow \infty}{\longrightarrow} 0.
		\end{align*}
		Thus, taking the limit as $k$ goes to $\infty$ in \eqref{uniqueness-expectation}, we obtain
		\begin{equation}
			\E^{P^1}\left[ \int_0^t f(r,\xi_r) \d r\right]=-{u_\infty}(0,x)=\E^{P^2}\left[ \int_0^t f(r,\xi_r) \d r\right],
		\end{equation}
		which readily gives $P^1=P^2$ (see e.g. Theorem 4.2 in \cite{EK86}).

\appendix
\section{Proof of the technical lemmas involving the stable density}\label{APP_HK}
\subsection{Proof of Lemma \ref{lemma-stable-sensitivities} (Stable Sensitivities)}
	
	Item \eqref{derivatives-palpha} directly follows from Section 2 in \cite{Kol00}. Let us prove \eqref{holder-space-palpha}.
	\begin{itemize}
		\item \textbf{Diagonal case:} $|x-x'|\leq u^{1/\alpha}$. Since we are looking at a small perturbation in the space variable, it makes sense to use a Taylor expansion:
		\begin{align*}
			\left|\nabla_x^\zeta p_\alpha (u,x)-\nabla_x^\zeta p_\alpha (u,x')\right| &= \left| \int_0^1 \nabla_x\nabla_x^{\zeta} p_\alpha (u,x'+(x-x')\lambda) \cdot (x-x') \d \lambda\right|\\
			&\lesssim  \frac{|x-x'|}{u^{\frac{1+|\zeta|}{\alpha}}}\int_0^1 \bar{p}_\alpha (u,x'+(x-x')\lambda) \d \lambda,
		\end{align*}
		using \eqref{derivatives-palpha} and $p_\alpha \asymp \bar{p}_\alpha$ (see \eqref{ARONSON_STABLE} and \eqref{DEF_P_BAR}) for the last inequality.
		Up to a modification of the underlying constant,
		\begin{align*}
			\bar{p}_\alpha (u,x'+(x-x')\lambda)
			&\lesssim u^{-\frac{d}{\alpha}}\left(2+\frac{|x'+\lambda (x-x')|}{u^{\frac{1}{\alpha}}} \right)^{-d-\alpha}
			\lesssim u^{-\frac{d}{\alpha}}\left(2-\frac{|x-x'|}{u^{\frac 1\alpha}}+\frac{|x'|}{u^{\frac{1}{\alpha}}} \right)^{-d-\alpha} \\
			& \lesssim u^{-\frac{d}{\alpha}}\left(1 +\frac{|x'|}{u^{\frac{1}{\alpha}}} \right)^{-d-\alpha} \lesssim \bar{p}_\alpha (u,x').
		\end{align*}
		We conclude the proof in the diagonal case noting that for all $\theta\in (0,1] $, $\frac{|x-x'|}{u^{\frac{1}{\alpha}}}\le \big(\frac{|x-x'|}{u^{\frac{1}{\alpha}}}\big)^\theta$.
		\item \textbf{Off-diagonal case:} $|x-x'|\geq u^{1/\alpha}$. In this case, a Taylor expansion in space is not relevant. We simply use the fact that $1=\frac{|x-x'|^\theta}{u^{\frac{\theta}{\alpha}}}\wedge 1$ and \eqref{derivatives-palpha}:
		\begin{align*}
			\left|\nabla_x^\zeta p_\alpha (u,x)-\nabla_x^\zeta p_\alpha (u,x')\right| &\leq  \left(\frac{|x-x'|^\theta}{u^{\frac{\theta}{\alpha}}} \wedge 1 \right)\left( \left| \nabla_x^\zeta p_\alpha (u,x)\right|+\left|\nabla_x^\zeta p_\alpha (u,x') \right| \right) \\
			&\lesssim\left( \frac{|x-x'|^\theta}{u^{\frac{\theta}{\alpha}}}\wedge 1 \right)\frac{1}{u^{\frac{|\zeta|}{\alpha}}}\left(p_\alpha (u,x)+p_\alpha (u,x')\right).
		\end{align*}
	\end{itemize}
	This concludes the proof of \eqref{holder-space-palpha}.\\
	
	Let us now prove \eqref{holder-time-palpha}. Let $0<u\le u'\le T$. Assume first $|u-u'|\leq \frac{u}{2}$.
	\begin{align*}
		\left|\nabla_x^\zeta p_\alpha (u,x)-\nabla_x^\zeta p_\alpha (u',x)\right| &= \left| \int_0^1  \partial_t \nabla_x^\zeta p_\alpha (u+(u'-u)\lambda,x) (u'-u)\d \lambda \right|\\
		&\lesssim \int_0^1 \frac{1}{(u+(u'-u)\lambda)^{1+\frac{|\zeta|}{\alpha}}}\bar{p}_\alpha(u+(u'-u)\lambda,x) |u'-u| \d \lambda \\
		&\lesssim \frac{|u-u'|}{u^{1+\frac{|\zeta|}{\alpha}}} \int_0^1\bar{p}_\alpha(u+(u'-u)\lambda,x) \d \lambda, 
	\end{align*}
recalling that $ u'\ge u$ for the last inequality. We now discuss in function of the position of the spatial variable $x$ w.r.t. the current time $u$.
	\begin{itemize}
		\item \textbf{Diagonal case:} $|x|\leq u^{1/\alpha}$. Then, 
		\begin{align*}
			\bar{p}_\alpha (u+ (u'-u)\lambda,x) \lesssim  (u+(u'-u)\lambda)^{-\frac{d}{\alpha}} \lesssim u^{-\frac{d}{\alpha}}\asymp \bar{p}_\alpha (u,x) \asymp {p}_\alpha (u,x).
		\end{align*}
		\item\textbf{ Off-diagonal case:}  $|x|\geq u^{1/\alpha}$. 
		\begin{align*}
			\bar{p}_\alpha (u+ (u'-u)\lambda,x) \lesssim \frac{u+(u'-u)\lambda }{|x|^{d+\alpha}} \lesssim \frac{u}{|x|^{d+\alpha}} \asymp \bar{p}_\alpha(u,x)\asymp {p}_\alpha (u,x).
		\end{align*}
	\end{itemize}
	Note that the condition $|u'-u|\le \frac u 2 $ is actually needed only for the second above inequality. Namely, it ensures that the term $\lambda(u'-u) $ has the same magnitude than $u$ (otherwise the previous expansions are useless and the estimation is direct as discussed below).
	
	In turn, we obtain
	\begin{align*}
		\left|\nabla_x^\zeta p_\alpha (u,x)-\nabla_x^\zeta p_\alpha (u',x)\right| \lesssim \frac{|u-u'|}{u^{1+\frac{|\zeta|}{\alpha}}}   \bar{p}_\alpha(u,z) \lesssim \frac{|u-u'|^\theta}{u^{\theta+\frac{|\zeta|}{\alpha}}}   {p}_\alpha(u,z),
	\end{align*}
	for all $\theta\in (0,1] $.
	
	In the case $|u-u'|\geq \frac{u}{2}$, we simply write using \eqref{derivatives-palpha}
	\begin{align*}
		\left|\nabla_x^\zeta p_\alpha (u,x)-\nabla_x^\zeta p_\alpha (u',x)\right|&\leq \Big(2\frac{|u-u'|}{u}\Big)^\theta\left( |\nabla_x^\zeta p_\alpha (u,x)|+ |\nabla_x^\zeta p_\alpha (u',x)|\right)\\
		&\lesssim\frac{|u-u'|^\theta}{u^{\theta+\frac{|\zeta|}{\alpha}}}  \left(p_\alpha (u,x)+p_\alpha (u',x)\right),
	\end{align*}
	which concludes the proof of \eqref{holder-time-palpha}.

        Let us now prove \eqref{spatial-moments}. Using $\bar{p}_\alpha \asymp p_\alpha$, we can write
	\begin{align*}
		\Vert p_\alpha (u,\cdot)|\cdot|^\zeta \Vert_{L^{\ell '}}^{\ell '} &=\int_{\R^d} p_\alpha (u,y)^{\ell '} |y|^{\zeta \ell '} \d y \lesssim \int_{\R^d} \frac{1}{u^{\frac{d \ell '}{\alpha}}}\times \frac{1}{\left(1+\frac{|y|}{u^{\frac{1}{\alpha}}} \right)^{\ell ' (d+\alpha)} }|y|^{\zeta \ell '} \d y.
	\end{align*}
	Set $z=yu^{-\frac{1}{\alpha}}$:
	\begin{align*}
		\Vert p_\alpha (u,\cdot)|\cdot|^\zeta \Vert_{L^{\ell '}}^{\ell '} &\lesssim u^{\frac{d}{\alpha}(1-\ell ')+\frac{\zeta \ell '}{\alpha}} \int_{\R^d} \frac{1}{\left(1+|z|\right)^{\ell ' (d+\alpha)} }|z|^{\zeta \ell '} \d z,
	\end{align*}
	which converges whenever $\ell ' (\zeta -d-\alpha)+d-1<-1 \ssi \zeta<d+\alpha-\frac{d}{\ell'}$, in which case we obtain
	$$\Vert p_\alpha (u,\cdot)|\cdot|^\zeta \Vert_{L^{\ell '}}\lesssim u^{-\frac{d }{\alpha\ell}+\frac{\zeta }{\alpha}}.$$
	Let us now prove the convolution part (\ref{p-q-convo}). Denote
	\begin{align}
		\mathfrak{I} &:= \Vert \bar{p}_\alpha (t-u,\cdot-y) \bar{p}_\alpha (u-s,x-\cdot)\Vert_{L^{\ell'}}^{\ell'}\notag\\
		&\lesssim \int \frac{1}{(t-u)^{\frac{d\ell'}{\alpha}}} \times\frac{1}{\left( 1+\frac{|z-y|}{(t-u)^{\frac{1}{\alpha}}}\right)^{(d+\alpha)\ell'}}\times\frac{1}{(u-s)^{\frac{d\ell'}{\alpha}}} \times\frac{1}{\left( 1+\frac{|x-z|}{(u-s)^{\frac{1}{\alpha}}}\right)^{(d+\alpha)\ell'}}\d z.\label{DEF_I_FOR_CONV}
	\end{align}
		We now discuss in function of the magnitude of the distance $|x-y| $ w.r.t. to the global time scale $t-s$.

	\begin{itemize}
		\item \textbf{Diagonal case:} $|x-y|<(t-s)^{1/\alpha}$\\
		In this case, either $(t-u)\geq \frac{1}{2} (t-s)$ or $(u-s)\geq \frac{1}{2} (t-s)$, we can then use the  global diagonal bound in \eqref{DEF_I_FOR_CONV} for the corresponding density.
		\begin{itemize}
			\item If $(t-u)\geq \frac{1}{2} (t-s)$,
			\begin{align*}
				\mathfrak{I} 
				&\lesssim \frac{1}{(t-u)^{\frac{d\ell'}{\alpha}}}\times  \frac{1}{(u-s)^{\frac{d}{\alpha}(\ell'-1)}} \int \frac{1}{(u-s)^{\frac{d}{\alpha}}}\times \frac{1}{\left( 1+\frac{|x-z|}{(u-s)^{\frac{1}{\alpha}}}\right)^{(d+\alpha)\ell'}}\d z
				\lesssim \frac{1}{(t-u)^{\frac{d\ell'}{\alpha}}} \times \frac{1}{(u-s)^{\frac{d}{\alpha}(\ell'-1)}}.
			\end{align*}
			Since $(t-u)\geq \frac{1}{2} (t-s)$,
			$\frac{1}{(t-u)^{\frac{d}{\alpha}}} \lesssim \frac{1}{(t-s)^{\frac{d}{\alpha}}} \asymp \bar{p}_\alpha  (t-s,y-x),$
			and
			\begin{equation*}
				\mathfrak{I}\lesssim \bar{p}_\alpha (t-s,x-y)^{\ell'} \frac{1}{(u-s)^{\frac{d}{\alpha}(\ell'-1)}}.
			\end{equation*}
			\item If $(u-s)\geq \frac{1}{2} (t-s)$, we readily obtain by symmetry
			\begin{equation*}
				\mathfrak{I}\lesssim \bar{p}_\alpha (t-s,x-y)^{\ell'} \frac{1}{(t-u)^{\frac{d}{\alpha}(\ell'-1)}}.
			\end{equation*}
		\end{itemize}
		\item \textbf{Off-diagonal case:} $|x-y|\geq (t-s)^{1/\alpha}$\\
		In this case, either $|x-z|\geq\frac{1}{2}|x-y|$ or $|z-y|>\frac{1}{2}|x-y|$, i.e. one of the two contributions in $\mathfrak{I}$ is in the off-diagonal regime, allowing us to use \eqref{off-diag}. In this cas we split the upper-bound for $\mathfrak{I}$ in \eqref{DEF_I_FOR_CONV} as follows:
		\begin{align*}
\mathfrak I\lesssim \int \frac{1}{(t-u)^{\frac{d\ell'}{\alpha}}}\times \frac{1}{\left( 1+\frac{|z-y|}{(t-u)^{\frac{1}{\alpha}}}\right)^{(d+\alpha)\ell'}}\times \frac{1}{(u-s)^{\frac{d\ell'}{\alpha}}}\times \frac{(\mathbb{1}_{|x-z|\geq\frac{1}{2}|x-y|}+\mathbb{1}_{|z-y|\geq\frac{1}{2}|x-y|})}{\left( 1+\frac{|x-z|}{(u-s)^{\frac{1}{\alpha}}}\right)^{(d+\alpha)\ell'}}\d z=:\mathfrak I_1+\mathfrak I_2.
		\end{align*}
		\begin{itemize}
			\item For $\mathfrak I_1 $,  $|x-z|\geq\frac{1}{2}|x-y|>\frac{1}{2} (t-s)^{1/\alpha}$, we get 
			\begin{align*}
				\mathfrak{I}_1 &\lesssim \frac{1}{(u-s)^{\frac{d\ell'}{\alpha}}}\times \frac{1}{\left( 1+\frac{|x-y|}{(u-s)^{\frac{1}{\alpha}}}\right)^{(d+\alpha)\ell'}} \int \frac{1}{(t-u)^{\frac{d\ell'}{\alpha}}} \times\frac{1}{\left( 1+\frac{|z-y|}{(t-u)^{\frac{1}{\alpha}}}\right)^{(d+\alpha)\ell'}} \mathbb{1}_{|x-z|\geq\frac{1}{2}|x-y|} \d z\\
				&\lesssim \bar{p}_\alpha (u-s,x-y)^{\ell'} \frac{1}{(t-u)^{\frac{d}{\alpha}(\ell'-1)}}\int \frac{1}{(t-u)^{\frac{d}{\alpha}}}\times \frac{1}{\left( 1+\frac{|z-y|}{(t-u)^{\frac{1}{\alpha}}}\right)^{(d+\alpha)\ell'}}\mathbb{1}_{|x-z|\geq\frac{1}{2}|x-y|} \d z 
			\end{align*}
			
			Since $|x-y|>(u-s)^{1/\alpha}$,
			$\bar{p}_\alpha (u-s,x-y) \asymp \frac{u-s}{|x-y|^{d+\alpha}} \leq \frac{t-s}{|x-y|^{d+\alpha}} \asymp \bar{p}_\alpha (t-s,x-y),$
			and 
			\begin{equation*}
				\mathfrak{I}_1\lesssim \bar{p}_\alpha (t-s,x-y)^{\ell'} \frac{1}{(t-u)^{\frac{d}{\alpha}(\ell'-1)}}.
			\end{equation*}
			\item For $\mathfrak I_2 $, $|z-y|>\frac{1}{2}|x-y|>\frac{1}{2} (t-s)^{1/\alpha}$, the same computations give, when swapping the roles of $|x-z|$ and $|y-z|$,
			\begin{equation}
				\mathfrak{I}_2\lesssim \bar{p}_\alpha (t-s,x-y)^{\ell'} \frac{1}{(u-s)^{\frac{d}{\alpha}(\ell'-1)}}.
			\end{equation}
		\end{itemize}
	\end{itemize}
	
	In each case, we have established that
	\begin{align*}
		\Vert \bar{p}_\alpha (t-u,\cdot-y) \bar{p}_\alpha (u-s,x-\cdot)\Vert_{L^{\ell'}} &= \mathfrak{I}^{\frac{1}{\ell'}} \lesssim \left[ \frac{1}{(t-u)^{\frac{d}{\alpha}\frac{\ell'-1}{\ell'}}} +\frac{1}{(u-s)^{\frac{d}{\alpha}\frac{\ell'-1}{\ell'}}}\right] \bar{p}_\alpha (t-s,x-y) \nonumber \\
		& \lesssim \left[ \frac{1}{(t-u)^{\frac{d}{\alpha  \ell}}} +\frac{1}{(u-s)^{\frac{d}{\alpha  \ell}}}\right] \bar{p}_\alpha(t-s,x-y),
	\end{align*}
	which concludes the proof of \eqref{p-q-convo}.\\
	Equation \eqref{convo-space-sing} then eventually follows from \eqref{p-q-convo} and H\"older's inequality.
\subsection{Proof of Lemma \ref{lemFK} (Feynman-Kac partial differential equation)}
Recall that for $u>0, z\in \R^d$:
	\begin{align*}
p_\alpha(u,z)&=\frac{1}{(2\pi)^d}\int_{\R^d} \exp(u \Phi_{\mu,\alpha}(\zeta))\exp(-i\zeta \cdot z) \d\zeta,\\
\Phi_{\mu,\alpha}(\zeta)&=\int_{\R_+}\int_{\mathbb{S}^{d-1}}\left[ \exp(i \zeta \cdot \rho \xi) - 1\right]\mu(\mathrm{d} \xi)\frac{\d \rho}{\rho^{1+\alpha}} =-C_{\alpha,d}\int_{\mathbb S^{d-1}} |\zeta\cdot \xi|^\alpha\mu({\mathrm{d}}\xi),\ C_{\alpha,d}>0,
\end{align*}
being the Khinchin exponent associated with the operator $\mathcal L^\alpha$.
	It is thus direct to see from the non-degeneracy assumption \eqref{minmajleb} that there exists $c>0$ s.t. $\forall \zeta\in\R^d,\;\int_{\mathbb S^{d-1}} |\zeta\cdot \xi|^\alpha\mu({\rm d}\xi)\ge c|\zeta|^\alpha$ so that $ \exp(u\Phi_{\mu,\alpha}(\zeta))\le \exp(- cC_{\alpha,d}u|\zeta|^\alpha)$. We deduce that $p_\alpha $ is smooth on $\R_+^*\times \R^d $ and
	\begin{align*}
          \partial_u p_\alpha(u,z)
                                    =\frac{1}{(2\pi)^d}\int_{\R^d} \Phi_{\mu,\alpha}(\zeta) \exp(u \Phi_{\mu,\alpha}(\zeta))\exp(-i\zeta \cdot z) \d\zeta.\end{align*}
        Since, by symmetry of the measure $\mu$ and Fubini's theorem,
\begin{align*}
  \Phi_{\mu,\alpha}(\zeta) \exp(u \Phi_{\mu,\alpha}(\zeta))&=\int_{\R^d}\int_{\R_+}\int_{\mathbb{S}^{d-1}}\left[ \exp(-i \zeta \cdot \rho \xi) - 1\right]\mu(\mathrm{d} \xi)\frac{\d \rho}{\rho^{1+\alpha}}\int_{\R^d}\exp(i\zeta \cdot x)p_\alpha(u,x)\d x\\&=\int_{\R_+}\int_{\mathbb{S}^{d-1}}\int_{\R^d}\left(\exp(i \zeta \cdot (x- \rho \xi))-\exp(i \zeta \cdot x)\right)p_\alpha(u,x)\d x\mu(\mathrm{d} \xi)\frac{\d \rho}{\rho^{1+\alpha}}
\\&=\int_{\R_+}\int_{\mathbb{S}^{d-1}}\int_{\R^d}\exp(i \zeta \cdot x)\left(p_\alpha(u,x+\rho\xi)-p_\alpha(u,x)\right)\d x\mu(\mathrm{d} \xi)\frac{\d \rho}{\rho^{1+\alpha}}\\&=\int_{\R^d}\exp(i\zeta \cdot x)\mathcal L^\alpha p_\alpha(u,x)\d x,
\end{align*}
one has $\partial_up_\alpha(u,z)=\mathcal L^\alpha p_\alpha(u,z)$.\\
        
        The fact that $v$ solves the Feynman-Kac partial differential equation on $[0,t)\times\R^d$ is easily deduced using \eqref{ARONSON_STABLE} and \eqref{derivatives-palpha} to apply Lebesgue's and Fubini's theorems. Last, for $s\in [0,t)$, 
        \begin{align*}
        |v(s,y)-\phi(y)|&=\int_{\R^d}\left|\phi\left(y-(t-s)^{\frac{1}{\alpha}}z\right)-\phi(y)\right|(t-s)^{\frac{d}{\alpha}}p_\alpha(t-s,(t-s)^{\frac{1}{\alpha}}z)\d z\\
        &=\int_{\R^d}\left|\phi\left(y-(t-s)^{\frac{1}{\alpha}}z\right)-\phi(y)\right|p_\alpha(1,z)\d z,
        \end{align*}
 and Lebesgue's theorem ensures that the right-hand side converges to $0$ as $s$ goes to $t$.\\
	
We refer to \cite{Kol00} (in particular the introduction, Proposition 2.5 and Section 3) for additional details and properties about the density $p_\alpha$.
\subsection{Proof of Lemma \ref{lemma-convo-bulk}: stable time-space convolutions with Lebesgue function}
	Let us first use \eqref{convo-space-sing} with $\ell=p$ and then H\"older's inequality in time:
	\begin{align*}
		I_{\beta_1,\beta_2}(u,v):=\int_{u}^{v}\int &p_\alpha (r,z-x)|b(r,z)|p_\alpha (t-r,y-z) \frac{1}{(t-r)^{\beta_1}}\frac{1}{r^{\beta_2}} \d z\d r\\
		&\lesssim p_\alpha (t,y-x) \int_{u}^{v} \Vert b(r,\cdot)\Vert_{L^p} \left[\frac{1}{r^{\frac{d}{\alpha p}}} + \frac{1}{(t-r)^{\frac{d}{\alpha p}}}\right]\frac{1}{(t-r)^{\beta_1}}\frac{1}{r^{\beta_2}} \d r\\
		&\lesssim p_\alpha (t,y-x)\left( \int_{u}^{v} \left[\frac{1}{r^{\frac{dq' }{\alpha p}}} + \frac{1}{(t-r)^{\frac{dq'}{\alpha p}}}\right]\frac{1}{(t-r)^{q' \beta_1}}\frac{1}{r^{q'\beta_2}} \d r\right)^{\frac{1}{q'}}\\
		&=:p_\alpha (t,y-x)S_{\beta_1,\beta_2}(u,v).
	\end{align*}
	Set $\lambda =\frac{r-u}{v-u} \iff r=u+\lambda (v-u)$, then 
	\begin{align*}
	&S_{\beta_1,\beta_2}(u,v)^{q'}\\
	\lesssim & (v-u)\int_{0}^{1} \left[\frac{1}{(u+\lambda (v-u))^{\frac{dq' }{\alpha p}}} + \frac{1}{(t-u-\lambda (v-u))^{\frac{dq'}{\alpha p}}}\right]\frac{1}{(t-u-\lambda(v-u))^{q' \beta_1}}\frac{1}{(u+\lambda(v-u))^{q'\beta_2}} \d \lambda.
	\end{align*}
	Assume first that $q' \left(\frac{d}{\alpha p}+\beta_i\right)<1,\ i\in \{1,2\}$ (integrable case). Then,
	\begin{align*}
	&S_{\beta_1,\beta_2}(u,v)^{q'}\\
	\lesssim & (v-u)^{1-\frac{dq' }{\alpha p}-q'(\beta_1+\beta_2)}\int_{0}^{1} \left[\frac{1}{\lambda ^{\frac{dq' }{\alpha p}}} + \frac{1}{(1-\lambda)^{\frac{dq'}{\alpha p}}}\right]\frac{1}{(1-\lambda)^{q' \beta_1}}\frac{1}{\lambda^{q'\beta_2}} \d \lambda\lesssim (v-u)^{1-\frac{dq' }{\alpha p}-q'(\beta_1+\beta_2)},
	\end{align*}
and
$$I_{\beta_1,\beta_2}(u,v)\lesssim p_\alpha(t,y-x)(v-u)^{1-\frac 1q-\frac{d }{\alpha p}-(\beta_1+\beta_2)} $$
which, recalling \eqref{gap}, gives \eqref{lemma-convo-bulk-2}.	

Let us now consider the case $q' \left(\frac{d}{\alpha p}+\beta_1\right)>1, q' \left(\frac{d}{\alpha p}+\beta_2\right)<1$ (singular case) with $v<t$. We then write:
	\begin{align*}
	&S_{\beta_1,\beta_2}(u,v)^{q'}\\
	\lesssim & (v-u)^{1-q'(\frac{d}{\alpha p}+\beta_1+\beta_2)}\int_{0}^{1} \left[\frac{1}{\lambda^{q'(\frac{d}{\alpha p}+\beta_2)}(\frac{t-u}{v-u}-\lambda)^{q'\beta_1}} + \frac{1}{(\frac{t-u}{v-u}-\lambda )^{q'(\frac{d}{\alpha p}+\beta_1)}\lambda^{q'\beta_2}}\right]\d \lambda\\
	\lesssim & (v-u)^{1-q'(\frac{d}{\alpha p}+\beta_1+\beta_2)}\Big(\int_0^{\frac 12}\frac{1}{\lambda^{q'(\frac d{\alpha p}+\beta_2)}}\d \lambda+\int_{\frac 12}^1\frac{1}{(\frac{t-u}{v-u}-\lambda )^{q'(\frac{d}{\alpha p}+\beta_1)}}\d \lambda \Big)\\
	\lesssim & (v-u)^{1-q'(\frac{d}{\alpha p}+\beta_1+\beta_2)}\Big(1+(\frac{t-u}{v-u}-1 )^{1-q'(\frac{d}{\alpha p}+\beta_1)}\Big)\\
	\lesssim&  (v-u)^{1-q'(\frac{d}{\alpha p}+\beta_1+\beta_2)}+(v-u)^{-q'\beta_2}(t-v)^{1-q'(\frac{d}{\alpha p}+\beta_1)}.
	\end{align*}
Hence, in the divergent case we have established
$$I_{\beta_1,\beta_2}(u,v)\lesssim p_\alpha(t,y-x)\Big((v-u)^{1-\frac 1q-\frac{d }{\alpha p}-(\beta_1+\beta_2)}+(v-u)^{-\beta_2}(t-v)^{1-\frac 1q-(\frac{d}{\alpha p}+\beta_1)}\Big), $$
which precisely gives \eqref{lemma-convo-bulk-1}. This concludes the proof of  Lemma \ref{lemma-convo-bulk}.

\subsection{Proof of Lemma \ref{lemma-ignore-bh} (About the cutoff on a one-step transition)}
  Using the fact that $p_\alpha \asymp \bar{p}_\alpha$ and $\left|y-sb_h(r,x)\right| \geq |y|-s|b_h(r,x)|$, we get for $0\le s\le \min (u,h) $,
	\begin{align}
		\bar{p}_\alpha (u,y-sb_h(r,x)) 
		&\lesssim \frac{1}{u^{\frac{d}{\alpha}}} \left( 2+\frac{|y-sb_h(r,x)|}{u^{\frac{1}{\alpha}}} \right)^{-(d+\alpha)} 
		\lesssim \frac{1}{u^{\frac{d}{\alpha}}} \left( 2-\frac{s}{u^{\frac{1}{\alpha}}}| b_h(r,x)|+ \frac{|y|}{u^{\frac{1}{\alpha}}} \right)^{-(d+\alpha)}\nonumber\\
		&\lesssim \frac{1}{u^{\frac{d}{\alpha}}} \left( 2-s^{1-\frac{1}{\alpha}}h^{-\frac{d}{\alpha p}-\frac{1}{q}} + \frac{|y|}{u^{\frac{1}{\alpha}}} \right)^{-(d+\alpha)}
		\lesssim \frac{1}{u^{\frac{d}{\alpha}}} \left( 2-h^{\frac{\gamma}{\alpha}} + \frac{|y|}{u^{\frac{1}{\alpha}}} \right)^{-(d+\alpha)}\nonumber\\
		&\lesssim \frac{1}{u^{\frac{d}{\alpha}}} \left( 1 + \frac{|y|}{u^{\frac{1}{\alpha}}} \right)^{-(d+\alpha)}
		\lesssim \bar{p}_\alpha(u,y)  \label{one-step-ineq},
	\end{align}	 
	provided $h<1$ for the last inequality, which we can assume w.l.o.g.\\
	In the case of $\bar{b}_h$, we derive similarly,
	\begin{align*}
		\bar{p}_\alpha (u,y-s \bar{b}_h(r,x)) 
		&\lesssim \frac{1}{u^{\frac{d}{\alpha}}} \left( 2- \frac{s}{u^{\frac{1}{\alpha}}}| \bar{b}_h(r,x)|+ \frac{|y|}{u^{\frac{1}{\alpha}}} \right)^{-(d+\alpha)}
		\lesssim \frac{1}{u^{\frac{d}{\alpha}}} \left( 2- h^0+ \frac{|y|}{u^{\frac{1}{\alpha}}} \right)^{-(d+\alpha)}\\
		&\lesssim \frac{1}{u^{\frac{d}{\alpha}}} \left( 1+ \frac{|y|}{u^{\frac{1}{\alpha}}} \right)^{-(d+\alpha)}=\bar{p}_\alpha (u,y).
	\end{align*}	
	This proves \eqref{ignore-bh-eq1} for $|\zeta|=0$. For $0<|\zeta|\leq 1$, one simply needs to apply \eqref{derivatives-palpha} beforehands.	
	For the proof of \eqref{ignore-bh-eq2}, it is enough to apply \eqref{holder-space-palpha} to $\left| \nabla^\zeta p_\alpha (u,y-s\mathfrak b_h(r,x))-\nabla^\zeta p_\alpha (u,y'-s\mathfrak b_h(r,x)) \right| $, where $ \mathfrak b_h\in \{ b_h,\bar b_h\}$, which yields
	\begin{align*}
		&\left| \nabla^\zeta p_\alpha (u,y-s\mathfrak b_h(r,x))-\nabla^\zeta p_\alpha (u,y'-s\mathfrak b_h(r,x)) \right|\\
		& \qquad\qquad\qquad\qquad\lesssim \left(\frac{|y-y'|^\delta}{u^{\frac{\delta}{\alpha}}}\wedge 1\right) \frac{1}{u^{\frac{|\zeta|}{\alpha}}}\left(p_\alpha (u,y-s\mathfrak b_h(r,x))+p_\alpha (u,y'-s\mathfrak b_h(r,x)) \right),
	\end{align*}
	for all $\delta\in (0,1] $ and then use \eqref{ignore-bh-eq1} to get rid of the drift in the previous equation.
	
	\subsection{Proof of Lemma \ref{lemma-schauder}: Schauder estimates for the mollified PDE \eqref{pde-mol}}\label{subsec-schauder-proofs}
	
	Let $m\in\N$ and $u_m$ denote the classical solution to $\eqref{pde-mol}$. For $s\in (0,t],x\in \R^d$, computing $u_m(t,x)+(Z_t-Z_s)-u_m(s,x)$ by Itô's formula and taking expectations, we obtain
	\begin{align}\label{duhamel-extended}
		u_m(s,x) &= - \int_s^t \int f(r,y)p_\alpha (r-s,y-x)\d y \d r\nonumber\\
		&\qquad + \int_s^t \int b_m(r,y)\cdot \nabla u_m(r,y) p_\alpha (r-s,y-x)\d y \d r\nonumber\\
		&=: I_1(s,x) + I_2(s,x).
	\end{align}
	Let us first prove the gradient bound \eqref{gradient-bound}. For $I_1$, using that $f$ is bounded on $[0,T]\times\R^d$ along with \eqref{derivatives-palpha}, we get
	\begin{align*}
		|\nabla I_1 (s,x)| & \leq \int_s^t \int |f(r,y)||\nabla_x p_\alpha (r-s,y-x)|\d y \d r \\
		&\lesssim \Vert f\Vert_{L^\infty-L^\infty}\int_s^t \int  p_\alpha (r-s,y-x)\frac{1}{(r-s)^{\frac{1}{\alpha}}}\d y \d r \\
		&\lesssim \Vert f\Vert_{L^\infty-L^\infty}\int_s^t  \frac{1}{(r-s)^{\frac{1}{\alpha}}}\d r \\
		&\lesssim \Vert f\Vert_{L^\infty-L^\infty}.
	\end{align*}
	For $I_2$, let us first note that due to standard Schauder estimates (see \cite{MP14}), we already know that $\nabla u_m(r,\cdot)$ is bounded (although not necessarily uniformly in $m$) for all $r\in [0,t]$. We can thus write, using a H\"older inequality, then \eqref{spatial-moments}, and finally a H\"older inequality in time,
	\begin{align*}
		|\nabla I_2 (s,x)| & \leq \int_s^t \Vert \nabla u_m (r,\cdot) \Vert_{L^\infty}\int |b_m(r,y)|  p_\alpha (r-s,y-x)\d y \d r\\
		& \leq \int_s^t \Vert \nabla u_m (r,\cdot) \Vert_{L^\infty} \Vert b_m(r,\cdot) \Vert_{L^p}  \frac{1}{(r-s)^{\frac{d}{\alpha p}}} \d r \\
		& \leq \Vert b_m \Vert_{L^q-L^p} \left(\int_s^t \Vert \nabla u_m (r,\cdot) \Vert_{L^\infty}^{q'}  \frac{1}{(r-s)^{\frac{d q'}{\alpha p}}} \d r \right)^{\frac{1}{q'}}.
	\end{align*}
	Gathering the previous estimates, we have 
	\begin{align*}
		\Vert\nabla u_m (s,\cdot)\Vert_{L^\infty}^{q'} & \lesssim\Vert f\Vert_{L^\infty-L^\infty}^{q'}+ \Vert b_m \Vert_{L^q-L^p} ^{q'}\int_s^t \Vert \nabla u_m (r,\cdot) \Vert_{L^\infty}^{q'}  \frac{1}{(r-s)^{\frac{d q'}{\alpha p}}} \d r.
	\end{align*}
	Since $\frac{dq'}{\alpha p}<1$, using Lemma 2.2 and Example 2.4 \cite{ZhangJFA10}, we deduce \eqref{gradient-bound}.\\
	Let us now prove \eqref{time-holder}. Using the previous notations, we can write
	\begin{align*}
		|u_m (s',x) -u_m (s,x)| &\leq  | I_1 (s',x) -  I_1 (s,x)| + | I_2 (s',x) - I_2 (s,x)|.
	\end{align*}
	For the first term, using \eqref{holder-time-palpha} with $\theta = \xi$, for any $\xi \in [0,(\gamma+1)/\alpha)$ and $\zeta=0$, we readily have
	\begin{align*}
		| I_1 (s',x) - I_1 (s,x)|  \lesssim |s'-s|^{\xi}  \Vert f \Vert_{L^\infty-L^\infty}.
	\end{align*}
	For the second term, using \eqref{holder-time-palpha} with $\theta =\xi$  and $\zeta=0$, as well as \eqref{derivatives-palpha}, we can write
	\begin{align*}
		| I_2 (s',x) - I_2 (s,x)| &\lesssim \int_{s'}^t \int |b_m(r,y)\cdot \nabla u_m(r,y)| \frac{|{s'}-s|^{\xi}}{(r-{s'})^{\xi}}(p_\alpha (r-{s'},y-x)+p_\alpha (r-s,y-x))\d y \d r\\
		&\qquad + \int_s^{s'}\int |b_m(r,y)\cdot \nabla u_m(r,y)| p_\alpha (r-s,y-x)\d y \d r.
	\end{align*}
	Using then a H\"older inequality in space, \eqref{spatial-moments}, a H\"older inequality in time and the previously established boundedness of $\nabla u_m$, we get
	\begin{align*}
		| I_2 (s',x) -  I_2 (s,x)| &\lesssim \Vert b_m \Vert_{L^q-L^p} \Vert \nabla u_m \Vert_{L^\infty} |{s'}-s|^{\xi}\left(\int_{s'}^t  \frac{1}{(r-{s'})^{q'\xi}}\left[\frac{1}{(r-s)^{\frac{dq'}{\alpha p}}} + \frac{1}{(r-{s'})^{\frac{dq'}{\alpha p}}}\right] \d r\right)^{\frac{1}{q'}}\\
		&\qquad + \Vert b_m \Vert_{L^q-L^p} \Vert \nabla u_m \Vert_{L^\infty} \left(\int_s^{s'} \frac{1}{(r-s)^{\frac{dq'}{\alpha p}}}\d r\right)^{\frac{1}{q'}}.
	\end{align*}
	Notice now that
	$$q'\left(\xi + \frac{d}{\alpha p}\right)<q'\left(\frac{\gamma}{\alpha} +\frac 1\alpha+ \frac{d}{\alpha p}\right) = 1,$$
	which concludes the proof of \eqref{time-holder}.
	\color{black}

	\bibliographystyle{alpha}
	\bibliography{ART}

\begin{thebibliography}{CdRM22}

\bibitem[Bas11]{Bas11}
Richard Bass.
\newblock {\em Stochastic Processes}.
\newblock Cambridge University Press, 2011.

\bibitem[BJ22]{BJ20}
Oumaima Bencheikh and Benjamin Jourdain.
\newblock Convergence in total variation of the euler--maruyama scheme applied
  to diffusion processes with measurable drift coefficient and additive noise.
\newblock {\em SIAM Journal on Numerical Analysis}, 60(4):1701--1740, 2022.

\bibitem[CdRM22]{CdRM22}
Paul-\'{E}ric Chaudru~de Raynal and St{\'e}phane Menozzi.
\newblock {On multidimensional stable-driven stochastic differential equations
  with {B}esov drift}.
\newblock {\em Electronic Journal of Probability}, 27(none):1 -- 52, 2022.

\bibitem[CHZ20]{CHZ20}
Zhen-Qing Chen, Zimo Hao, and Xicheng Zhang.
\newblock H\"{o}lder regularity and gradient estimates for {SDE}s driven by
  cylindrical {$\alpha$}-stable processes.
\newblock {\em Electron. J. Probab.}, 25:Paper No. 137, 23, 2020.

\bibitem[CZZ21]{CZZ21}
Zhen-Qing Chen, Xicheng Zhang, and Guohuan Zhao.
\newblock Supercritical {SDE}s driven by multiplicative stable-like {L}\'{e}vy
  processes.
\newblock {\em Trans. Amer. Math. Soc.}, 374(11):7621--7655, 2021.

\bibitem[DGI22]{DGI22}
Tiziano {De Angelis}, Maximilien Germain, and Elena Issoglio.
\newblock A numerical scheme for stochastic differential equations with
  distributional drift.
\newblock {\em Stochastic Processes and their Applications}, 154:55--90, 2022.

\bibitem[EK86]{EK86}
Stewart Ethier and Thomas Kurtz.
\newblock {\em Markov Processes: Characterization and Convergence}.
\newblock John Wiley and Sons, New York, 1986.

\bibitem[Esc06]{Esc06}
Carlos Escudero.
\newblock The fractional {K}eller-{S}egel model.
\newblock {\em Nonlinearity}, 19(12):2909--2918, 2006.

\bibitem[Fit23]{Fit23}
Mathis Fitoussi.
\newblock Heat kernel estimates for stable-driven {SDE}s with distributional
  drift.
\newblock {\em Potential Analysis}, 2023.

\bibitem[Fri64]{AF64}
Avner Friedman.
\newblock {\em Partial Differential Equations of Parabolic Type}.
\newblock Dover Publications, 1964.

\bibitem[Hol22]{Hol22}
Teodor Holland.
\newblock A note on the weak rate of convergence for the {E}uler-{M}aruyama
  scheme with {H}\"older drift.
\newblock {\em arXiv, 2206.12830}, 2022.

\bibitem[JM23]{JM23}
Benjamin Jourdain and St{\'e}phane Menozzi.
\newblock {Convergence Rate of the Euler-Maruyama Scheme Applied to Diffusion
  Processes with $L^q$ -- $L^\rho$ Drift Coefficient and Additive Noise}.
\newblock {\em Annals of Applied Probability}, 2023.

\bibitem[KK18]{KK18}
Victoria Knopova and Alexei Kulik.
\newblock Parametrix construction of the transition probability density of the
  solution to an {SDE} driven by {$\alpha$}-stable noise.
\newblock {\em Ann. Inst. Henri Poincar\'{e} Probab. Stat.}, 54(1):100--140,
  2018.

\bibitem[Kol00]{Kol00}
Vassili Kolokoltsov.
\newblock Symmetric stable laws and stable-like jump-diffusions.
\newblock {\em Proc. London Math. Soc. (3)}, 80(3):725--768, 2000.

\bibitem[KR05]{KR05}
Nicolai~V. Krylov and Michael R\"{o}ckner.
\newblock Strong solutions of stochastic equations with singular time dependent
  drift.
\newblock {\em Probab. Theory Related Fields}, 131(2):154--196, 2005.

\bibitem[Kul19]{Kul19}
Alexei Kulik.
\newblock On weak uniqueness and distributional properties of a solution to an
  {SDE} with {$\alpha$}-stable noise.
\newblock {\em Stochastic Process. Appl.}, 129(2):473--506, 2019.

\bibitem[L{\^e}20]{Le20}
Khoa L{\^e}.
\newblock {A stochastic sewing lemma and applications}.
\newblock {\em Electronic Journal of Probability}, 25(none):1 -- 55, 2020.

\bibitem[LL22]{LL22}
Khoa L\^e and Chengcheng Ling.
\newblock Taming singular stochastic differential equations: A numerical
  method, 2022.

\bibitem[MP91]{MP91}
Remigijus Mikulevicius and Eckhard Platen.
\newblock Rate of convergence of the {E}uler approximation for diffusion
  processes.
\newblock {\em Mathematische Nachrichten}, 1991.

\bibitem[MP14]{MP14}
Remigijus Mikulevicius and Henrikas Pragarauskas.
\newblock On the {C}auchy problem for integro-differential operators in
  {H}\"{o}lder classes and the uniqueness of the martingale problem.
\newblock {\em Potential Anal.}, 40(4):539--563, 2014.

\bibitem[MS12]{MS12}
Mark~M. Meerschaert and Alla Sikorskii.
\newblock {\em Stochastic models for fractional calculus}, volume~43 of {\em De
  Gruyter Studies in Mathematics}.
\newblock Walter de Gruyter \& Co., Berlin, 2012.

\bibitem[MZ22]{MZ22}
St\'ephane Menozzi and Xicheng Zhang.
\newblock Heat kernel of supercritical nonlocal operators with unbounded
  drifts.
\newblock {\em Journal de l{\textquoteright}\'Ecole polytechnique {\textemdash}
  Math\'ematiques}, 9:537--579, 2022.

\bibitem[Por94]{Por94}
N.~I. Portenko.
\newblock Some perturbations of drift-type for symmetric stable processes.
\newblock {\em Random Oper. Stochastic Equations}, 2(3):211--224, 1994.

\bibitem[PP95]{PP95}
S.~I. Podolynny and N.~I. Portenko.
\newblock On multidimensional stable processes with locally unbounded drift.
\newblock {\em Random Oper. Stochastic Equations}, 3(2):113--124, 1995.

\bibitem[PvZ22]{PvZ22}
Nicolas Perkowski and Willem van Zuijlen.
\newblock Quantitative heat-kernel estimates for diffusions with distributional
  drift.
\newblock {\em Potential Analysis}, 2022.
\newblock https://doi.org/10.1007/s11118-021-09984-3.

\bibitem[Sat99]{Sat99}
Ken-iti Sato.
\newblock {\em L\'{e}vy Processes and Infinitely divisible Distributions}.
\newblock Cambridge University Press, 1999.

\bibitem[TT90]{TT90}
D.~Talay and L.~Tubaro.
\newblock Expansion of the global error for numerical schemes solving
  sto\-chastic differential equations.
\newblock {\em Stoch. Anal. and App.}, 8-4:94--120, 1990.

\bibitem[TTW74]{TTW74}
Hiroshi Tanaka, Masaaki Tsuchiya, and Shinzo Watanabe.
\newblock {Perturbation of drift-type for L\'{e}vy processes}.
\newblock {\em Journal of Mathematics of Kyoto University}, 14(1):73 -- 92,
  1974.

\bibitem[Wat07]{Wa07}
Toshiro Watanabe.
\newblock Asymptotic estimates of multi-dimensional stable densities and their
  applications.
\newblock {\em Transactions of the American Mathematical Society},
  359(6):2851--2879, 2007.

\bibitem[XZ20]{XZ20}
Longjie Xie and Xicheng Zhang.
\newblock Ergodicity of stochastic differential equations with jumps and
  singular coefficients.
\newblock {\em Ann. Inst. Henri Poincar\'{e} Probab. Stat.}, 56(1):175--229,
  2020.

\bibitem[Zha10]{ZhangJFA10}
Xicheng Zhang.
\newblock Stochastic {V}olterra equations in {B}anach spaces and stochastic
  partial differential equation.
\newblock {\em J. Funct. Anal.}, 258(4):1361--1425, 2010.

\end{thebibliography}

\end{document}